\documentclass[]{amsart}

\usepackage{amscd}
\usepackage{amsfonts}
\usepackage{amssymb}
\usepackage{amsmath}
\numberwithin{equation}{section}

\newtheorem{pro}{Proposition}[section]
\newtheorem{lemma}[pro]{Lemma}
\newtheorem{theorem}[pro]{Theorem}
\newtheorem{definition}[pro]{Definition}
\newtheorem{remark}[pro]{Remark}

\newtheorem{corollary}[pro]{Corollary}

\newcommand{\fg}{\frak g}

\newcommand{\fp}{\frak p}

\newcommand{\fs}{\frak s}

\newcommand{\Ll}{\mathcal{L}}
\newcommand{\al}{\alpha}
\newcommand{\be}{\beta}
\newcommand{\ga}{\gamma}

\newcommand{\f}{\varphi}

\newcommand{\se}{\theta}
\newcommand{\om}{\omega}
\newcommand{\we}{\wedge}
\newcommand{\lam}{\lambda}
\newcommand{\sig}{\sigma}
\newcommand{\ra}{{\rightarrow}}
\newcommand{\lra}{{\longrightarrow}}

\newcommand{\ie}{{i.e$.$\,}}
\newcommand{\cf}{{cf$.$\,}}

\newcommand{\CC}{{\mathbb C}}

\newcommand{\HH}{{\mathbb H}}

\newcommand{\PP}{{\mathbb P}}
\newcommand{\RR}{{\mathbb R}}

\newcommand{\ZZ}{{\mathbb Z}}

\newcommand{\GL}{\mathop{\rm GL}\nolimits}

\newcommand{\PSp}{\mathop{\rm PSp}\nolimits}

\newcommand{\Null}{\mathop{\rm Null}\nolimits}
\newcommand{\Sim}{\mathop{\rm Sim}\nolimits}
\newcommand{\Aut}{\mathop{\rm Aut}\nolimits}

\newcommand{\dev}{\mathop{\rm dev}\nolimits}

\begin{document}
\baselineskip 13pt
\pagestyle{myheadings}
\thispagestyle{empty}
\setcounter{page}{1}

\title[\ ]
{Pseudo-conformal quaternionic $CR$ structure
on $(4n+3)$-dimensional manifolds}

\author[\ ]{Dmitri ALEKSEEVSKY and Yoshinobu KAMISHIMA}
\address{Department of Mathematics, Hull University,
GEN. Cottingham Road HU6 7RX, Hull, England}
\email{d.v.alekseevsky@maths.hull.ac.uk}
\address{Department of Mathematics, Tokyo Metropolitan 
University,Minami-Ohsawa 1-1,
Hachioji, Tokyo 192-0397, Japan}
\email{kami@comp.metro-u.ac.jp}

\date{\today}
\date{\today}

\begin{abstract} We study a geometric structure on a $(4n+3)$-dimensional smooth manifold $M$ which is an integrable, nondegenerate codimension 
$3$-subbundle $\mathcal D$ on $M$ whose fiber supports
the structure of $4n$-dimensional quaternionic vector space.
It is thought of as a generalization of the quaternionic $CR$ structure. 
In order to study this geometric structure on $M$,
we single out an  ${\mathfrak s}{\mathfrak p}(1)$-valued 
$1$-form $\om$ locally on a neighborhood $U$ of
$M$ such that ${\rm Null}\om=\mathcal D|U$.
We shall construct the invariants on the pair $(M,\om)$
whose vanishing implies that $M$ is uniformized with respect to
a finite dimensional flat quaternionic $CR$ geometry.
The invariants obtained on $(4n+3)$-$M$ have the same formula as
the curvature tensor of quaternionic (indefinite) K\"ahler manifolds.
From this viewpoint,
we shall exhibit a quaternionic analogue of Chern-Moser's
$CR$ structure.
\end{abstract}

\maketitle

 
\setcounter{section}{0}
\section*{Introduction}\label{INT}
The Weyl curvature tensor is a conformal invariant
of Riemannian manifolds and the Chern-Moser curvature tensor
is a $CR$ invariant on strictly pseudo-convex $CR$-manifolds.
 A geometric significance of the vanishing of these curvature tensors
 is the appearance of the finite dimensional Lie group $\mathcal G$
with homogeneous space $X$.
The geometry $(\mathcal G, X)$ is known as 
conformally flat geometry $({\rm PO}(n+1,1), S^n)$,
spherical $CR$-geometry $({\rm PU}(n+1,1), S^{2n+1})$
respectively.
Similarly, the complete simply connected
quaternionic $(n+1)$-dimensional
quaternionic hyperbolic space 
${\HH}^{n+1}_{\bf H}$ with the group of isometries 
${\rm PSp}(n+1,1)$ has the natural compactification
homeomorphic to a $(4n+4)$-ball endowed with an extended smooth action of  
${\rm PSp}(n+1,1)$.
When the boundary sphere $S^{4n+3}$ of the ball
is viewed as the real hypersurface
in the quaternionic projective space ${\HH}{\PP}^{n+1}$,
the elements of ${\rm PSp}(n+1,1)$ act as quaternionic projective 
transformations of $S^{4n+3}$.
Since the action of ${\rm PSp}(n+1,1)$ is transitive
on $S^{4n+3}$, we obtain a flat (spherical) quaternionic $CR$
geometry $({\rm PSp}(n+1,1),S^{4n+3})$. (Compare \cite{KA}.)
Combined with the above two geometries, this exhibits {\it parabolic geometry}
on the boundary of the compactification of rank one-symmetric
space of noncompact type.
(See \cite{BS},\cite{CM},\cite{WE},\cite{KA1}.)
This observation naturally leads us to the problem
of the existence of a geometric structure
on a $(4n+3)$-dimensional manifold $M$ and the
geometric invariant whose vanishing implies that
$M$ is locally equivalent to the flat
quaternionic $CR$ manifold $S^{4n+3}$.
For this purpose we shall introduce a notion of
 pseudo-conformal quaternionic $CR$ structure 
$(\mathcal D,\{\om_\al\}_{\al=1,2,3})$
on a $(4n+3)$-dimensional manifold $M$. 
First of all we recall a pseudo-conformal quaternionic
structure which is discussed in \cite{AK1}.
Contrary to the nondegenerate $CR$ structure,
the almost complex structure on $\mathcal D$
is not assumed to be integrable.
However, by the requirement of our structure equations,
we can prove the integrability of
quaternionic structure in $\S \ref{threeintegrable}$:
\vskip0.2cm
\par\noindent{\bf Theorem A. }  
{\em Each almost complex structure $\bar J_\al$ of the
quaternionic $CR$ structure is integrable on the codimension $1$
contact subbundle $\mathop{\Null} \om_\al$ $(\al=1,2,3)$.
 }
\vskip0.2cm
A $(4n+3)$-dimensional complete
simply connected quaternionic pseudo-Riemannian space form 
$\Sigma^{3+4p,4q}_{\HH}$ of type $(3+4p,4q)$ with constant curvature $1$
has been introduced in \cite{KU1} $(p+q=n)$.
In $\S \ref{MCR}$, we show that this is a model space of nondegenerate
quaternionic $CR$ structure.
There exists a canonical pseudo-Riemannian metric $g$
associated to the nondegenerate pseudo-conformal quaternionic $CR$
structure. Then we see in  $\S \ref{principal}$
that the integrability of three
almost complex structures $\{\bar J_\al\}_{\al=1,2,3}$
is equivalent with the
condition that $(M,g)$ is a pseudo-Sasakian $3$-structure.
(Compare \cite{AK_1}.)
Using the pseudo-Riemannian connection of
the pseudo-Sasakian $3$-structure, we can define a quaternionic $CR$ curvature tensor (\cf $\S \ref{EP}$).
Based on this curvature tensor, we shall establish
a curvature tensor
$T$ which is invariant under
the equivalence of pseudo-conformal quaternionic $CR$ structures
in $\S \ref{EP1}$. The $(4n+3)$-dimensional manifold $S^{3+4p,4q}$ introduced
in $\S \ref{existPQCR}$ is a pseudo-conformal quaternionic $CR$ manifold
 with all vanishing pseudo-conformal $QCR$ curvature tensor. The model
 $S^{3+4p,4q}$ contains $\Sigma^{3+4p,4q}_{\HH}$ as an open dense subset.
In $\S \ref{CI}$, we shall describe an explicit formula of our tensor $T$
(cf. Theorem \ref{formula}).
\vskip0.2cm
\par\noindent{\bf Theorem B.}\  
{\em
Let $T=({T}^i_{jk\ell})$ be the fourth-order curvature tensor
on a nondegenerate pseudo-conformal $QCR$ manifold $M$ 
in dimension $4n+3$ $(n\geq 0)$.
If $n\geq 2$, then
$T=({T}^i_{jk\ell})\in \mathcal R_0({\rm Sp}(p,q)\cdot{\rm Sp}(1))$
which has the formula:
\begin{equation*}\label{Tensor}
\begin{split}
{T}^i_{jk\ell}&=R^i_{jk\ell}-\Bigl\{(g_{j\ell}\delta^{i}_k
-g_{jk}\delta^{i}_\ell)+
\Bigl[{I}_{j\ell}{I}^{i}_{k}-{I}_{jk}{I}_{i}^{\ell}+
2{I}^{i}_{j}{I}_{k\ell}\Bigr.\Bigr.\\
& \ \ \ \ \ \ \Bigl.\Bigl.+{J}_{j\ell}{J}^{i}_{k}-{J}_{jk}{J}_{i}^{\ell}+
2{J}^{i}_{j}{J}_{k\ell}
+{K}_{j\ell}{K}^{i}_{k}-{K}_{jk}{K}^{i}_{\ell}+
2{K}^{i}_{j}{K}_{k\ell} \Bigr]\Bigr\}.
\end{split}
\end{equation*}
When $n=1$,
$T=({W}^i_{jk\ell})\in \mathcal R_0({\rm SO}(4))$ which
has the same formula as the Weyl conformal curvature tensor.
When $n=0$, there exists the fourth-order curvature tensor
$TW$ on $M$ which has the same formula as the
Weyl-Schouten tensor.
}
\vskip0.2cm

We shall prove that
the vanishing of curvature tensor $T$ on 
a nondegenerate pseudo-conformal quaternionic $CR$ manifold $M$ of type 
$(3+4p,4q)$
gives rise to a {\em uniformization} with respect to the
flat (spherical) pseudo-conformal quaternionic $CR$
geometry in $\S \ref {VC}$, see
Theorem \ref{uni1}.
(Compare \cite{KU} for uniformization in general.)

\vskip0.2cm
\par\noindent{\bf Theorem C.}\ \
{\em 
\begin{itemize}
\item[(i)]
If $M$ is a $(4n+3)$-dimensional 
 nondegenerate pseudo-conformal quaternionic $CR$
manifold of type $(3+4p,4q)$
$(p+q=n\geq 1)$ whose curvature tensor $T$ vanishes, then
$M$ is uniformized over $S^{3+4p,4q}$
with respect to the group ${\rm PSp}(p+1,q+1)$.
\item[(ii)]
If $M$ is a $3$-dimensional pseudo-conformal quaternionic $CR$
manifold whose
curvature tensor $TW$ vanishes, then
$M$ is conformally flat $($\ie locally modelled on $S^3$
with respect to the group ${\rm PSp}(1,1)$ $)$.
\end{itemize}

}
\vskip0.2cm
In particular, if $p=n,q=0$, then
$S^{3+4n,0}=S^{4n+3}$ is the
positive-definite flat (spherical) quaternionic  $CR$
and our pseudo-conformal quaternionic $CR$ geometry
is the spherical quaternionic $CR$ geometry 
$({\rm PSp}(n+1,1),S^{4n+3})$ as explained in the beginning
of Introduction.
By Theorem C,
if $M$ is a flat (spherical) pseudo-conformal positive definite
quaternionic 
$CR$ manifold, then
there exists a developing
map $\mathop{\dev}:\tilde M\ra S^{4n+3}=S^{4n+3,0}$
from the universal covering space $\tilde M$.
It is an immersion preserving the pseudo-conformal quaternionic $CR$
structure such that $\mathop{\dev}^*\om_0=
\lam\cdot \tilde\om\cdot \bar\lam$ where $\tilde \om$ is
the lift of $\om$ to $\tilde M$. As the global case,
positive-definite pseudo-conformal quaternionic $CR$
manifolds contain the class of
$3$-Sasakian manifolds.
(Refer to \cite{BL},\cite{BGM},\cite{TA},\cite{TA1} for Sasakian structure.)
However, we emphasize that the converse is not true.
We shall recall two typical classes of compact (spherical)
pseudo-conformal quaternionic $CR$ manifolds but not Sasakian $3$-manifolds
\cite{KA};
one is a quaternionic Heisenberg manifold ${\mathcal M}/\Gamma$.
Some finite cover of ${\mathcal M}/\Gamma$ is a Heisenberg nilmanifold
which is a principal $3$-torus bundle over
the flat quaternionic $n$-torus $T^n_{\HH}$. (Compare $\S \ref{pcHei}$.)
Another manifold is a pseudo-Riemannian standard space form 
$\Sigma^{3,4n}_{\HH}/\Gamma$ of constant negative curvature of
type $(4n,3)$. It is a compact quotient of the
homogeneous space $\Sigma^{3,4n}_{\HH}={\rm Sp}(1,n)/{\rm Sp}(n)$.
Some finite cover
of $\Sigma^{3,4n}_{\HH}/\Gamma$ is a principal $S^3$-bundle over the
quaternionic hyperbolic space form ${\HH}_{\HH}^n/\Gamma^*$.
Those manifolds are not positive-definite compact $3$-Sasakian manifolds.
(Compare \cite{KA}, \cite{KA2} more generally.)

When a geometric structure is either contact structure or complex contact
structure, it is known that the first Stiefel-Whitney class
or the first Chern class is the obstruction to the existence
of global $1$-forms representing their  strutures respectively.
As a concluding remark to the pseudo-conformal quaternionic
structure but not necessarily pseudo-conformal quaternionic $CR$
structure, we verify that 
the obstruction relates to the first Pontrjagin class $p_1(M)$
of a $(4n+3)$-dimensional pseudo-conformal quaternionic
manifold $M$ $(n\geq 1)$.
In $\S \ref{q-bundle}$, we prove that
the following relation of the first Pontrjagin 
classes. (See Theorem \ref{pont}.)

\vskip0.2cm
\par\noindent{\bf Theorem D.}\ \
{\em  
Let $(M,\mathcal D)$ be a $(4n+3)$-dimensional
pseudo-conformal quaternionic manifold.
Then  the first Pontrjagin classes of $M$ and the bundle $L=TM/\mathcal D$ 
has the relation that $\displaystyle 2p_1(M)=(n+2)p_1(L)$.
Moreover, if $M$ is simply connected, then
the following are equivalent.
\begin{enumerate}
\item $2p_1(M)=0$. In particular, the first rational Pontrjagin class vanishes. 
\item There exists a global 
${\rm Im}\HH$-valued
$1$-form $\omega$ on $M$ which represents a pseudo-conformal
quaternionic structure $\mathcal D$. In particular,
there exists a hypercomplex structure $\{I,J,K\}$ on $\mathcal D$.
\end{enumerate}

}
\par\ \par

\section{Pseudo-conformal quaternionic $CR$ structure}\label{QCCS}
When $\HH$ denotes the field of quaternions, the Lie algebra
${\mathfrak {sp}}(1)$ of ${\rm Sp}(1)$ is identified with ${\rm Im}{\HH}=
{\RR}\mbox{\boldmath$i$}+{\RR}\mbox{\boldmath$j$}+{\RR}\mbox{\boldmath$k$}$.
Let $M$ be a $(4n+3)$-dimensional smooth manifold $M$.
A $4n$-dimensional orientable subbundle $\mathcal D$ equipped with
a quaternionic structure $Q$ is called  a {\em pseudo-conformal quaternionic}
structure on $M$ if it satisfies that
\begin{equation}\label{qcc}
\begin{split}
(i)\ & \mathcal D\cup [\mathcal D,\mathcal D]=TM.\\
(ii)\ & \mbox{The $3$-dimensional quotient bundle 
$TM/\mathcal D$ at any point}\\
 \ \ \  & \mbox{is isomorphic to the Lie algebra ${\rm Im}{\HH}$}.\\
(iii)\ & \mbox{There exists a ${\rm Im}\HH$-valued $1$-form 
$\omega=\omega_1\mbox{\boldmath$i$}+\omega_2\mbox{\boldmath$j$}+
\omega_3\mbox{\boldmath$k$}$ }\\
\ \ \ \  & \mbox{locally defined on a neighborhood of $M$ such that}\\
\ \ \ \  & \mbox{ $\displaystyle \mathcal D=\mathop{\Null} \omega=\mathop{\cap}_{\al=1}^{3}\mathop{\Null} \omega_\al$ and $d\om_\al |\mathcal D$
is nondegenerate.}\\ 
\ \ \ \ & \mbox{Here each $\omega_\al$ is a real valued $1$-form $(\al=1,2,3)$.}\\
(iv)\ & \mbox{The endomorphism $J_\ga=(d\om_\be|\mathcal D)^{-1}\circ 
(d\om_{\al}|\mathcal D):\mathcal D
\ra \mathcal D$ }\\
\ \ \ & \mbox{constitutes the quaternionic structure $Q$ on $\mathcal D$:}\\
\ \ \ & {J_\ga}^2=-1,\ J_\al J_\be=J_\ga=-J_\be J_\al,\ (\ga=1,2,3)\ etc.
\end{split}\end{equation}

\par\ \par

We note the following from the condition $(iv)$.

\begin{lemma}\label{symmetric-bi}
Put $\sigma_\al=(d\om_\al|\mathcal D)$ on $\mathcal D$. 
There is the following equality:
\[
\sigma_1(J_1X,Y)=\sigma_2(J_2X,Y)=\sigma_3(J_3X,Y)\
(\forall\ X,Y\in \mathcal D).
\]Moreover, the form 
\begin{equation}\label{1trans}
g^{\mathcal D}=\sigma_\al\circ J_\al
\end{equation}is a nondegenerate
$Q$-invariant symmetric bilinear form on $\mathcal D$ \ $(\al=1,2,3)$,
\ie $\displaystyle g^{\mathcal D}(X,Y)=g^{\mathcal D}(J_\al X,J_\al Y)$,
$g^{\mathcal D}(X,J_\al Y)=\sigma_\al(X,Y)$, etc.,
$(\al=1,2,3)$.
\end{lemma}

\begin{proof}
By $(iv)$, it follows that
\begin{equation}\label{JJ}
\begin{split}
\sigma_\al(J_\al X,Y)&=\sigma_\al(J_\be(J_\ga X),Y)=\sigma_\ga(J_\ga X,Y)\\
&=\sigma_\ga(J_\al (J_\be X),Y)=\sigma_\be(J_\be X,Y).
\end{split}
\end{equation} 
Put $g^{\mathcal D}(X,Y)=\sigma_\al(J_\al X,Y)$ for $X,Y\in \mathcal D$
$(\al=1,2,3)$, which is nondegenerate by $(iii)$. 
As $-J_\be=\sigma_\ga^{-1}\circ \sigma_{\al}$ by $(iv)$,
calculate
that $\displaystyle
g^{\mathcal D}(Y,X)=-\sigma_\al(X,J_\al Y)=
\sigma_\ga(J_\be X,J_\al Y)=-\sigma_\be(Y,J_\be X)
=g^{\mathcal D}(X,Y)$.
It folows that
$\displaystyle
g^{\mathcal D}(X,Y)=\sigma_\al(J_\al X,Y)=\sigma_\al(J_\al(J_\al Y),J_\al X)
=g^{\mathcal D}(J_\al Y,J_\al X)$.
\end{proof}
In general, there is no canonical choice of 
$\omega$ which annihilates $\mathcal D$.
The fiber of the quotient bundle $TM/\mathcal D$ is isomorphic to 
${\rm Im}\ {\HH}$ by $\om$ on a neighborhood $U$ by $(ii)$.
The coordinate change of the fiber $\HH$
is described as $v\ra \lam\cdot v\cdot\mu$ for 
some nonzero scalars $\mu,\nu\in \HH$.
If $\om'$ is another $1$-form such that  $\mathop{\Null}\om'=\mathcal D$
on a neighborhood $U'$, then it follows that 
 $\om'=\lam\cdot \om\cdot \mu$ for some functions $\lam,\mu$ locally 
defined on $U\cap U'$. This can be rewritten as 
$\om'=u\cdot a\cdot \om\cdot b$
where $a,b$ are functions with valued in ${\rm Sp}(1)$
and $u$ is a
positive function. Since $\bar \om'=-\om'$, it follows that
$a\cdot \om\cdot b=\bar b\cdot \om\cdot\bar a$, 
\ie $(\bar b a)\cdot \om\cdot (\bar b a)=\om$. 
As $\om:T(U\cap U')\ra {\rm Im}\ {\HH}$ is surjective,
$\bar b a$ centralizes ${\rm Im}\ {\HH}$ so that $\bar ba\in \RR$.
Hence, $b=\pm {\bar a}$.
As we may assume that $\mathcal D$ is orientable, 
$\om'$ is uniquely determined
by
\begin{equation}\label{intersection1}
\om'=u\cdot a\cdot \om\cdot\bar a\ \mbox{for some functions}\
a\in{\rm Sp}(1),u>0\ \mbox{on}\ U\cap U'.
\end{equation}

We must show that the definition of \eqref{qcc} does not depend on the
choice of any form $\om'$ of \eqref{intersection1} locally
conjugate to $\om$, 
\ie it satisfies $(iii)$, $(iv)$.
Differentiate the equation \eqref{intersection1} which yields
$\displaystyle d\om'=u\cdot a\cdot d\om\cdot\bar a$
on $\mathcal D|U\cap U'$.
Let $\displaystyle A= (a_{ij})\in {\rm SO}(3)$
be the matrix function determined by
\begin{equation}\label{barA}
{\rm Ad}_{a}\left(\begin{array}{c}
\mbox{\boldmath$i$}\\
\mbox{\boldmath$j$}\\
\mbox{\boldmath$k$}\end{array}\right)
=a\left(\begin{array}{c}
\mbox{\boldmath$i$}\\
\mbox{\boldmath$j$}\\
\mbox{\boldmath$k$}\end{array}\right)\bar a
=A\left(\begin{array}{c}
\mbox{\boldmath$i$}\\
\mbox{\boldmath$j$}\\
\mbox{\boldmath$k$}\end{array}\right).
\end{equation}
When $\om$ is replaced by $\om'$,
a new quaternionic structure on $\mathcal D$ is obtained 
as
\begin{equation}\label{two-quaternionic}
\left(\begin{array}{c}
{J_1'}\\
{J_2'}\\
{J_3'}
\end{array}\right)
={}^tA\left(\begin{array}{c}
J_1\\
J_2\\
J_3
\end{array}\right).
\end{equation}
Note from \eqref{intersection1} that
\begin{equation}\label{intersection2}
\begin{split}
({\om'}_1, {\om'}_2, {\om'}_3)
=(\om_1, \om_2, \om_3)u\cdot A
=u(\mathop{\sum}_{\be=1}^3 a_{\be 1}\om_\be,
\mathop{\sum}_{\be=1}^3 a_{\be2}\om_\be,
\mathop{\sum}_{\be=1}^3 a_{\be3}\om_\be).
\end{split}
\end{equation} Differentiate
\eqref{intersection2} and restricting to $\mathcal D$,
\begin{equation*}\begin{split}
d{\om'}_\al(X,Y)&=u\mathop{\sum}_\be a_{\be\al}d\om_\be(X,Y)\
(\mbox{Lemma \ref{symmetric-bi}}) \\
&=-u(a_{1\al}g^{\mathcal D}(J_1X,Y)+a_{2\al}g^{\mathcal D}(J_2X,Y)
+a_{3\al}g^{\mathcal D}(J_3X,Y))\\
&=-ug^{\mathcal D}((a_{1\al}J_1+a_{2\al}J_2+a_{3\al}J_3)X,Y)
=-ug^{\mathcal D}({J'}_\al X,Y).
\end{split}
\end{equation*}It follows 

\begin{equation}\label{3coincide}
d{\om'}_\al({J'}_\al X,Y)=ug^{\mathcal D}(X,Y)\ (\al=1,2,3).
\end{equation}
In particular, $d\om'_\al|\mathcal D$ is nondegenerate, proving $(iii)$.
Put $\sigma'_\al=d\om'_\al|\mathcal D$. As in $(iv)$,
 the endomorphism is defined by the rule:
$\displaystyle I'_\ga=(\sigma'_\be|\mathcal D)^{-1}\circ 
(\sigma'_{\al}|\mathcal D)$, \ie 
$\displaystyle \sigma'_\be(I'_\ga X,Y)=\sigma'_{\al}(X,Y) \ (\forall X,Y\in\mathcal D.)$. Then we show that $\{I'_\al\}_{\al=1,2,3}$ 
coincides with $\{J'_\al\}_{\al=1,2,3}$  on $\mathcal D$.
For this, as $\sigma'_\al(X,Y)=-ug^{\mathcal D}(J'_\al X,Y)$ by \eqref{3coincide},
it follows that $\sigma'_\be(I'_\ga X,Y)=-ug^{\mathcal D}(J'_\be (I'_\ga X),Y)$ 
and the above equality implies that
$J'_\be (I'_\ga X)=J'_\al X$ $(\forall X\in\mathcal D)$. Hence,
$I'_\ga = -J'_\be J'_\al=J'_\ga$.
This proves $(iv)$.

The nondegenerate bilinear form $g^{\mathcal D}$ is locally defined on $\mathcal D|U$ of signature $(4p,4q)$ with $4p$-times 
positive sign and $4q$-times negative sign such that
$p+q=n$.
As above put ${g'}^{\mathcal D}(X,Y)=d{\om'}_\al({J'}_\al X,Y)$
\ $(X,Y\in\mathcal D)$. We have
\begin{corollary}\label{conf-metric}
If $\om'=u\bar a\cdot \om\cdot a$ on $U\cap U'$, then
\[
{g'}^{\mathcal D}=u\cdot g^{\mathcal D}.
\]
In particular, the signature $(p,q)$ is constant
on $U\cap U'$ (and hence everywhere on $M$)
under the  change $\om'=u\bar a\cdot \om\cdot a$.
\end{corollary}

Consider locally the structure equation:
\begin{equation}\label{integraleq1}
\rho_\al=d\om_\al+2\om_\be\we\om_\ga
\end{equation}where $(\al,\be,\ga)\sim (1,2,3)$ up to cyclic permutation.
If the skew symmetric $2$-form $\rho_\al$ satisfies that
\begin{equation}\label{localflow}
\mathop{\Null} \rho_1=\mathop{\Null} \rho_2=\mathop{\Null} \rho_3,
\end{equation}then the pair $(\om, Q)$ is a
 {\em local quaternionic $CR$ structure} on $M$. See \cite{BI}, \cite{AK_1}.
If the quaternionic $CR$ structure is globally
defined on $M$, \ie there exists a ${\rm Im}\HH$-valued $1$-form $\om$ defined entirely on $M$, then it coincides with 
the pseudo-Sasakian $3$-structure of $M$,
see $\S \ref{PS3}$.
 
We introduce the following definition  due to the manner of Libermann \cite{LI}.
\begin{definition}\label{pc-qcr}
The pair $(\mathcal D, Q)$ on $M$ is said to be 
a {\em pseudo-conformal quaternionic $CR$} structure 
if there exists locally a $1$-form $\eta$ with $\mathop{\Null}\eta=\mathcal D$
on a neighborhood $U$ of $M$ such that $\eta$ is conjugate to a 
quaternionic $CR$ structure on $U$. Namely there exists a
${\rm Im}\HH$-valued $1$-form $\om$ representing
the quaternionic $CR$ structure of $U$ for which
$\eta=\lam\cdot \om\cdot \bar\lam$ where
$\lam:U\ra \HH$ is a function and
$\bar\lam $ is the conjugate of the quaternion.
\end{definition}

\section{Quaternionic $CR$ structure}\label{defqcc}
Suppose that $\om$ is a quaternionic $CR$ structure on a neighborhood of $M$.
Let $\rho_\al=d\om_\al+2\om_\be\we\om_\ga$ be as in
\eqref{integraleq1}.
Put $V=\mathop{\Null} \rho_\al$ $(\al=1,2,3)$ (\cf \eqref{localflow}).
Since ${\rm dim} \mathcal D=4n$,
let $\{v_1,v_2,v_3\}$  be  a basis of $V$.
Put $\om_i(v_j)=a_{ij}$. As $\om_1\we\om_2\we\om_3|V\neq 0$,
the $3\times 3$-matrix $(a_{ij})$ is nonsingular. Put
$b_{ij}={}^t(a_{ij})^{-1}$ and $\xi_j=\sum b_{jk}v_k$.
Then  $\om_\al(\xi_\be)=\delta_{\al\be}$ and locally,
\begin{equation}\label{delta} 
V=\{\xi_\al,\al=1,2,3\}
\end{equation}

\begin{lemma}\label{property1}
Let $\Ll$ be the Lie derivative.
Then, 
$\Ll_{\xi_\al}(\mathcal D)=\mathcal D$\ \ $(\al =1,2,3)$.
\end{lemma}
\begin{proof}
For $X\in \mathcal D$,
 $\om_\be({\Ll}_{\xi_\al}(X))=\om_\be([\xi_\al,X])$.
As
\[
0=\rho_\be(\xi_\al,X)=d\om_\be(\xi_\al,X)+2\om_\ga\we\om_\al
(\xi_\al,X)=\frac 12(-\om_\be([\xi_\al,X]),
\]we have $\om_\be([\xi_\al,X])=0$ for $\be=1,2,3$.
Hence, $\displaystyle {\Ll}_{\xi_\al}(X)\in \mathcal D=
\mathop{\cap}_{\be=1}^{3}{\rm Null}\ \om_\be$.

\end{proof}
We prove also that $\Ll_\xi V=V$ for $\xi\in V$.
\begin{lemma}\label{so3}
The distribution $V$ is integrable. The vector fields 
$\xi_\al$ determined by \eqref{delta} generates the Lie algebra isomorphic
to $\mathfrak {so}(3)$, \ie $[\xi_\al,\xi_\be]=2\xi_\ga$.
\ $(\al,\be,\ga)\sim(1,2,3)$.
\end{lemma}

\begin{proof}
By \eqref{delta}, note that
\begin{equation}\label{localflow1}
\begin{split}
V=\{\xi\in TM\ |\ \rho_1(\xi,v)=\rho_2(\xi,v)=\rho_3(\xi,v)=0, \ \forall v\in
TM \}=\{\xi_\al\ ;\al=1,2,3\}.
\end{split}
\end{equation}
Since $\displaystyle 0=\rho_\al(\xi_\be,\xi_\ga)=
\frac 12(-\om_\al([\xi_\be,\xi_\ga])+2)$,
it follows that $[\xi_\be,\xi_\ga]-2\xi_\al\in {\rm Null}\ \om_\al$.
Applying $\displaystyle\rho_\be(\xi_\be,\xi_\ga)
=\frac 12(-\om_\be([\xi_\be,\xi_\ga])+0)=0$, it yields
also that $[\xi_\be,\xi_\ga]-2\xi_\al\in {\rm Null}\ \om_\be$. 
Similarly as $\rho_\ga(\xi_\be,\xi_\ga)=0$, we obtain
$\displaystyle [\xi_\be,\xi_\ga]-2\xi_\al\in\mathop{\cap}_{\be=1}^{3}
{\rm Null}\ \om_\be=\mathcal D$ for $\al=1,2,3$.
As $\displaystyle \rho_\al([\xi_\be,\xi_\ga]-2\xi_\al,v)=
\rho_\al([\xi_\be,\xi_\ga],v)$ for arbitrary $v\in \mathcal D$,
By the definition of $\rho_\al$,
calculate
\begin{equation*}
\begin{split}
\rho_\al([\xi_\be,\xi_\ga],v)&=
-\frac 12\om_\be([[\xi_\be,\xi_\ga],v])\\
&=\frac 12(\om_\be([[\xi_\ga,v],\xi_\be])+\om_\be([[\xi_\be,v],
\xi_\ga]))\ \mbox{(by Jacobi identity)}\\
&=0\ \mbox{(by Lemma \ref{property1})}.
\end{split}
\end{equation*}
Since $\rho_\al$ is nondegenerate on $\mathcal D$ by $(iii)$,
$[\xi_\be,\xi_\ga]=2\xi_\al$ $(\al=1,2,3)$.
Hence, such a Lie algebra $V$ is locally isomorphic to the Lie algebra of
${\rm SO}(3)$.  
\end{proof}
 
We collect the properties of $\om_\al,\rho_\al,J_\al,g^{\mathcal D}$.
(Compare \cite{AK_1}.)
\begin{lemma}\label{properties}
Up to cyclic permutation of $(\al,\be,\ga)\sim
(1,2,3)$, the following properties hold.

\begin{itemize}
\item[(1)] $\Ll_{\xi_\al}\om_\al =0$, $\Ll_{\xi_\al}\om_\be=\om_\ga
=-\Ll_{\xi_\be}\om_\al$.
\item[(2)] $\Ll_{\xi_\al}\rho_\al =0$, $\Ll_{\xi_\al}\rho_\be=\rho_\ga
=-\Ll_{\xi_\be}\rho_\al$.
\item[(3)] $\Ll_{\xi_\al}J_\al =0$, $\Ll_{\xi_\al}J_\be=J_\ga
=-\Ll_{\xi_\be}J_\al$.
\item[(4)] $\Ll_{\xi_{\al}}g^{\mathcal D}=0$.
\end{itemize}

\end{lemma}
 \begin{proof}
(1). First note that
$\iota_{\xi_{\al}}\om_\al(x)=\om_\al(\xi_\al)=1$ $(x\in M)$,
$\iota_{\xi_{\al}}(\om_\be\we\om_\ga)(X)=
\om_\be\we\om_\ga(\xi_{\al},X)=0$ $(\al\neq \be,\ga)$,
and $\iota_{\xi_\al}\rho_\al(X)=\rho_\al(\xi_\al,X)=0$ 
by \eqref{localflow1}.

\begin{equation}\begin{split}
\Ll_{\xi_{\al}}\om_\al=(d\iota_{\xi_{\al}}+\iota_{\xi_{\al}}d)\om_\al
&= \iota_{\xi_{\al}}d\om_\al
=\iota_{\xi_{\al}}(-2\om_\be\we\om_\ga+\rho_\al)\ 
\mbox{by}\ \eqref {integraleq1}\\
&=-2\iota_{\xi_{\al}}(\om_\be\we\om_\ga)+\iota_{\xi_\al}\rho_\al=0,
\end{split}\end{equation}
Next, 
\[
\Ll_{\xi_{\al}}\om_\be=\iota_{\xi_{\al}}d\om_\be=
\iota_{\xi_{\al}}(-2\om_\ga\we\om_\al+\rho_\be)
=-2\iota_{\xi_{\al}}(\om_\ga\we\om_\al),\]
while 
$-2\iota_{\xi_{\al}}(\om_\ga\we\om_\al)(v)=0$ for
$v\ \not\in\mathop{\Null}\ \om_\ga$ and
$-2\iota_{\xi_{\al}}(\om_\ga\we\om_\al)(\xi_\ga)=1$.
Hence $\Ll_{\xi_{\al}}\om_\be=\om_\ga$.

\noindent (2).  
\begin{equation}
\begin{split}
\Ll_{\xi_\al}\rho_\be &=\Ll_{\xi_\al}(d\om_\be+2\om_\ga\we\om_\al)\\
&=(d\iota_{\xi_{\al}}+\iota_{\xi_{\al}}d)d\om_\be+
2\Ll_{\xi_{\al}}(\om_\ga\we\om_\al)\\
&=d\iota_{\xi_{\al}}d\om_\be+2\Ll_{\xi_{\al}}\om_\ga\we\om_\al+
2\om_\ga\we\Ll_{\xi_{\al}}\om_\al\\
&=d(\Ll_{\xi_\al}-d\iota_{\xi_{\al}})\om_\be+2\Ll_{\xi_{\al}}\om_\ga\we\om_\al
\ \  (\mbox{by}\ (1))\\
&=d(\Ll_{\xi_\al}\om_\be)-2\Ll_{\xi_{\ga}}\om_\al\we\om_\al
=d\om_\ga-2\om_\be\we\om_\al\\
&=d\om_\ga+2\om_\al\we\om_\be=\rho_\ga.
\end{split}
\end{equation}Similarly,
\begin{equation}
\begin{split}
\Ll_{\xi_\al}\rho_\al &=\Ll_{\xi_\al}(d\om_\al+2\om_\be\we\om_\ga)\\
&=d\iota_{\xi_{\al}}d\om_\al+2\Ll_{\xi_{\al}}\om_\be\we\om_\ga+
2\om_\be\we\Ll_{\xi_{\al}}\om_\ga\\
&=d(\Ll_{\xi_\al}-d\iota_{\xi_{\al}})\om_\al+2\om_\ga\we\om_\ga
+2\om_\be\we (-\om_\be)\\
&=d\Ll_{\xi_\al}\om_\al=0\ \ (\mbox{by}\ (1)).
\end{split}
\end{equation}

\noindent (3).  
As $\Ll_{\xi_\al}\rho_\al=0$ by property $(2)$,
\begin{equation*}\begin{split}
0&=(\Ll_{\xi_\al}\rho_\al)(J_\be X,Y)\\
&=\Ll_{\xi_\al}(\sig_\al(J_\be X,Y))-\sig_\al(\Ll_{\xi_\al}(J_\be X),Y)
-\sig_\al(J_\be X,\Ll_{\xi_\al}Y).
\end{split}
\end{equation*}
Noting that
$J_\be=\sig_\al^{-1}\circ \sig_{\ga}$ by Lemma \ref{symmetric-bi},
we have
\begin{equation}\begin{split}
&\sigma_\al((\Ll_{\xi_\al}J_\be) X,Y)
=\sigma_\al(\Ll_{\xi_\al}(J_\be X),Y)-\sigma_\al(J_\be\Ll_{\xi_\al}(X),Y)\\
&=\Ll_{\xi_\al}(\sigma_\al(J_\be X,Y))-\sigma_\al(J_\be X,\Ll_{\xi_\al} Y)
-\sigma_\al(J_\be\Ll_{\xi_\al}X,Y)\\
&=(\Ll_{\xi_\al}\sigma_\ga)(X,Y)=-\sigma_\be(X,Y)\ \ \mbox{(by property (2))}\\
&=\sig_\al(J_\ga X,Y)
\end{split}
\end{equation}

As $\sigma_\al$ is nondegenerate on $\mathcal D$,
$\Ll_{\xi_\al}J_\be=J_\ga$.
Similarly,
\begin{equation}\begin{split}
&\sigma_\ga((\Ll_{\xi_\al}J_\al) X,Y)
=\sigma_\ga(\Ll_{\xi_\al}(J_\al X),Y)-\sigma_\ga(J_\al\Ll_{\xi_\al}(X),Y)\\
&=-(\Ll_{\xi_\al}\sig_\ga)(J_\al X,Y)+
\Ll_{\xi_\al}(\sigma_\ga(J_\al X,Y))\\
&\ -\sigma_\ga(J_\al X,\Ll_{\xi_\al} Y)
-\sigma_\ga(J_\al\Ll_{\xi_\al}X,Y)\\
&=\sig_\be(J_\al X,Y)+\Ll_{\xi_\al}(\sigma_\be(X,Y))-\sigma_\be(X,\Ll_{\xi_\al} Y)-\sigma_\be(\Ll_{\xi_\al}X,Y)\\
&=\sig_\be(J_\al X,Y)+(\Ll_{\xi_\al}\sigma_\be)(X,Y)\\
&=-\sigma_\ga(X,Y)+\sigma_\ga(X,Y)=0,
\end{split}
\end{equation}it follows that $\Ll_{\xi_\al}J_\al=0$.
\smallskip

\noindent (4).
Recall from Lemma \ref{symmetric-bi}
that $g^{\mathcal D}(X,Y)=\sig_\al(J_\al X,Y)
=\rho_\al(J_\al X,Y)$ $(X,Y\in\mathcal D)$ 
for each $\al$.
Then
\begin{equation}\label{ghmetric}
\begin{split}
(\Ll_{\xi_\al}g^{\mathcal D})&(X,Y)=\xi_\al(g^{\mathcal D}(X,Y))-g^{\mathcal D}(\Ll_{\xi_\al}X,Y)-g^{\mathcal D}(X,\Ll_{\xi_\al}Y)\\
&=\xi_\al(\rho_\be(J_\be X,Y))-\rho_\be(J_\be\Ll_{\xi_\al}X,Y)
-\rho_\be(J_\be X,\Ll_{\xi_\al}Y).
\end{split}
\end{equation}On the other hand, 
$\Ll_{\xi_\al}\rho_\be=\rho_\ga$ by property (2) and
so 
\[\xi_\al(\rho_\be(J_\be X,Y))=\rho_\be(\Ll_{\xi_\al}J_\be X,Y)
+\rho_\be(J_\be X,\Ll_{\xi_\al}Y)+\rho_\ga(J_\be X,Y).\]
Substitute this
into the equation \eqref{ghmetric}.
\begin{equation*}\label{metric0}
\begin{split}
(\Ll_{\xi_\al}g^{\mathcal D})(X,Y)&=\rho_\be(\Ll_{\xi_\al}J_\be X,Y)
+\rho_\be(J_\be X,\Ll_{\xi_\al}Y)\\
&+\rho_\ga(J_\be X,Y)-\rho_\be(J_\be\Ll_{\xi_\al}X,Y)
-\rho_\be(J_\be X,\Ll_{\xi_\al}Y)\\
&=\rho_\be((\Ll_{\xi_\al}J_\be)X,Y)+
\rho_\ga(J_\be X,Y)\ \ (\mbox{by property}\ (3))\\
&=\rho_\be(J_\ga X,Y)+\rho_\ga(J_\be X,Y) =0,
\end{split}
\end{equation*}
hence, $\Ll_{\xi_\al}g^{\mathcal D}=0$.
 
\end{proof} 

\subsection{Three $CR$ structures}\label{threeintegrable} 
Let $(\{\om_\al\},\{J_\al\},\{\xi_\al\};\ {\al =1,2,3})$ be a
nondegenerate quaternionic $CR$ structure on
$U\subset M$ such that
$\displaystyle \mathcal D|U=\mathop{\cap}_{\al=1}^3 \mathop{\Null}\om_\al$.
We can extend the almost complex structure $J_\al$ to an almost complex
structure $\bar J_\al$ on $\mathop{\Null} \om_\al=\mathcal D\oplus
\{\xi_\be,\xi_\ga\}$ by
setting:
\begin{equation}\label{extending comp}
\begin{split}
&\bar J_\al |\mathcal D=J_\al,\\
&\bar J_\al\xi_\be=\xi_\ga, \bar J_\al \xi_\ga=-\xi_\be.
\end{split}
\end{equation}$(\al,\be,\ga)$ is a cyclic permutation of $(1,2,3)$.
 First of all, note the following formula (cf. \cite{KONO}):
\begin{equation}\label{formula1}
\begin{split}
\Ll_X(\iota_Yd\om_a)=\iota_{(\Ll_XY)}d\om_a+\iota_Y\Ll_Xd\om_a
=\iota_{[X,Y]}d\om_a+\iota_Y\Ll_Xd\om_a \ \ (\forall X,Y\in TU).
\end{split}
\end{equation}
Secondly, we remark the following.
\begin{lemma}\label{mutual-relation}
For $X\in \mathcal D$,
\[
\iota_Xd\om_a=\iota_{J_c X}d\om_b
\ \ (a,b,c)\sim (1,2,3).
\]
\end{lemma}

\begin{proof}
Let $TU=\mathcal D\oplus V$ where $V= \{\xi_1,\xi_2,\xi_3\}$.
If $X\in \mathcal D$, then 
$d\om_a(X,\xi)=0$ for 
$\forall\ \xi\in V$. As $d\om_b(J_c X,\xi)=0$ similarly,
it follows that $\iota_Xd\om_a=\iota_{J_c X}d\om_b=0$ on $V$.
If $Y\in \mathcal D$, calculate
\begin{equation*}\begin{split}
d\om_a(X,Y)&=-d\om_a(J_a(J_a X),Y)
=-d\om_b(J_b(J_a X),Y) \ (\mbox{from Lemma}\ \ref{symmetric-bi})\\
&=d\om_b(J_c X,Y),\ \mbox{hence}\ 
\iota_Xd\om_a=\iota_{J_c X}d\om_b\ \mbox{on}\ U.
\end{split}
\end{equation*}
\end{proof}
In particular, we have 
\begin{equation}\label{2-3}
\iota_Xd\om_2=\iota_{J_1 X}d\om_3\ \mbox{for}\ \forall\ X\in \mathcal D.
\end{equation}
There is the decomposition with respect to the almost complex structure
 $\bar J_1$: 
\begin{equation}
\begin{split}
{\Null}\ \om_1\otimes{\CC}&=\mathcal D\otimes{\CC}\oplus \{\xi_2,\xi_3\}\otimes
{\CC}\\
&=T^{1,0}\oplus T^{0,1}
\end{split}
\end{equation}where $T^{1,0}=
\mathcal D^{1,0}\oplus \{\xi_2-\mbox{\boldmath$i$}\xi_3\}$.
We shall observe that the same formula as in Lemma 6.8
of Hitchin \cite{HIT} can be also obtained for $\mathcal D$.
(We found Lemma 6.8 when we saw a key lemma to
the Kashiwada's theorem \cite{KASHIWA}.)

\begin{lemma}\label{mutual-relation1}
If $X,Y\in \mathcal D^{1,0}$, then
\[
\iota_{[X,Y]}d\om_2=\mbox{\boldmath$i$}\iota_{[X,Y]}d\om_3.
\]
\end{lemma}
\begin{proof}
Let $X\in \mathcal D^{1,0}$ so that $J_1X=\mbox{\boldmath$i$}X$,
then
\begin{equation}\label{left}
\begin{split}
\Ll_Xd\om_2&= (d\iota_X+\iota_X d)d\om_2=d(\iota_Xd\om_2)
=d(\iota_{J_1 X}d\om_3)\ (\mbox{by}\ \eqref{2-3})\\
&=\mbox{\boldmath$i$}(d\iota_{X})d\om_3=
\mbox{\boldmath$i$}(\Ll_X-\iota_{X}d)d\om_3
=\mbox{\boldmath$i$}\Ll_Xd\om_3.
\end{split}
\end{equation}

Applying $Y\in \mathcal D^{1,0}$ to the equation
\eqref{2-3} and
using the formula \eqref{formula1} 
(extended to a ${\CC}$-valued one),
\begin{equation}\label{right}
\begin{split}
\Ll_X(\iota_Yd\om_2)&=\Ll_X(\iota_{J_1 Y}d\om_3)
=\mbox{\boldmath$i$}\Ll_X(\iota_{Y}d\om_3) \ \mbox{(from \eqref{2-3})}\\
&=\mbox{\boldmath$i$}\iota_{[X,Y]}d\om_3+\iota_Y\mbox{\boldmath$i$}
\Ll_Xd\om_3\\
&=\mbox{\boldmath$i$}\iota_{[X,Y]}d\om_3+\iota_Y\Ll_Xd\om_2\ (\mbox{by}\ 
\eqref{left}).
\end{split}
\end{equation}
Compared this with the equality \eqref{formula1} to $\om_a=\om_2$,
we obtain $\mbox{\boldmath$i$}\iota_{[X,Y]}d\om_3=\iota_{[X,Y]}d\om_2$.

\end{proof}

We prove the following equation (which is used to show the existence of
 a complex contact structure on the quotient of the quaternionic $CR$ manifold 
by $S^1$ \cite {AK}.)

\begin{pro}\label{complex-contact}
For any $X,Y\in \mathcal D^{1,0}$, there exsist
$a\in {\RR}$ and $u\in \mathcal D^{1,0}$ such that
\[
[X,Y]=a(\xi_2-\mbox{\boldmath$i$}\xi_3)+u.
\]
Conversely, given an arbitrary $a\in {\RR}$, 
we can choose such $X,Y\in \mathcal D^{1,0}$ and some $u\in \mathcal D^{1,0}$.

\end{pro}

\begin{proof}
As $g(J_\al\cdot , J_\al\cdot )=g(\cdot,\cdot)$ (\cf Lemma \ref{symmetric-bi}),
 we note that $d\om_1|(\mathcal D^{1,0},\mathcal D^{0,1})$,
$d\om_2|(\mathcal D^{1,0},\mathcal D^{1,0})$, $d\om_3|(\mathcal D^{1,0},\mathcal D^{1,0})$
are nondegenerate. 
Given $X,Y\in \mathcal D^{1,0}$, put
$\displaystyle d\om_2(X,Y)=g(X,J_2Y)=-\frac12a$
for some $a\in {\RR}$. (Note that conversely 
for any $a\in {\RR}$, we can choose
$X,Y\in \mathcal D^{1,0}$
such that $\displaystyle d\om_2(X,Y)=g(X,J_2Y)=-\frac12a$.)
Then $\om_2([X,Y])=a$ so that there is an element
$v\in \mathop{\Null} \om_2\otimes{\CC}$
such that
$[X,Y]-a\cdot\xi_2=v$.
As $\displaystyle d\om_3(X,Y)=g(X,J_1J_2Y)=
-g(X,J_2(J_1Y))=-\mbox{\boldmath$i$}g(X,J_2Y)=
-\frac {\mbox{\boldmath$i$}}2a$,
it follows that $\om_3([X,Y])=-\mbox{\boldmath$i$}a$.
Since $\om_3(v)=\om_3([X,Y]-\xi_2)=\om_3([X,Y])$,
$v=-\mbox{\boldmath$i$}a\cdot \xi_3+u$ for some
 $u\in \mathop{\Null} \om_3\otimes{\CC}$.
Then we have that $\displaystyle [X,Y]=a(\xi_2-\mbox{\boldmath$i$}\xi_3)+u$.
Obviously, $\om_2(u)=0$.
As $X,Y\in \mathcal D^{1,0}$,
$\om_1(u)=\om_1([X,Y])=-2d\om_1(X,Y)=0$ for which
$u\in \mathcal D\otimes{\CC}$.
We now prove that $u\in \mathcal D^{1,0}$.
First we note that
\begin{equation}\label{braeq}
\iota_{[X,Y]}d\om_2=a\iota_{(\xi_2-\mbox{\boldmath$i$}\xi_3)}
d\om_2+\iota_{u}d\om_2.
\end{equation}
As $\xi_2$ (respectively $\xi_3$) is characteristic for $\om_2$ 
(respectively $\om_3$)
from Lemma \ref{properties}, 
$\iota_{\xi_2}d\om_2=0$ (respectively $\iota_{\xi_3}d\om_3=0$).
Using \eqref{localflow1}, the function satisfies $d\iota_{\xi_3}\om_2=0$
(respectively $d\iota_{\xi_2}\om_3=0$). 
It follows that $\displaystyle \iota_{\xi_3}d\om_2=
(\Ll_{\xi_3}-d\iota_{\xi_3})\om_2=\Ll_{\xi_3}\om_2= -\om_1$.
Then $\displaystyle
\iota_{(\xi_2-\mbox{\boldmath$i$}\xi_3)}d\om_2=
(\iota_{\xi_2}d\om_2-\mbox{\boldmath$i$}\iota_{\xi_3}d\om_2)
=\mbox{\boldmath$i$}\om_1$ so
\eqref{braeq} becomes
\begin{equation}\label{form-f}
\iota_{[X,Y]}d\om_2=a\mbox{\boldmath$i$}\om_1+\iota_{u}d\om_2.
\end{equation}
As $\Ll_{\xi_2}\om_3=\om_1$, it follows
$\iota _{\xi_2}d\om_3=\om_1$.
Similarly
\begin{equation}\label{form-ff}
\iota_{[X,Y]}d\om_3=a\iota_{(\xi_2-\mbox{\boldmath$i$}\xi_3)}
d\om_3+\iota_{u}d\om_3=
a\om_1+\iota_{u}d\om_3.
\end{equation}
Substitute \eqref{form-f}, \eqref{form-ff}
into the equlaity $\iota_{[X,Y]}d\om_2=\mbox{\boldmath$i$}\iota_{[X,Y]}d\om_3$
of Lemma \ref{mutual-relation1}, which concludes
that
\begin{equation}\label{reduce H}
\iota_{u}d\om_2=\mbox{\boldmath$i$}\iota_{u}d\om_3.
\end{equation}
Since $d\om_2(u,X)=d\om_3(J_1u,X)$ for any $X\in \mathcal D\otimes {\CC}$,
\eqref{reduce H} implies that
\begin{equation}
d\om_3(J_1u,X)=\iota_{u}d\om_2(X)=d\om_3(\mbox{\boldmath$i$}u,X).
\end{equation}
As $d\om_3$ is nondegenerate on $\mathcal D$ (and so is on $\mathcal D\otimes {\CC}$), 
we obtain that $J_1u=\mbox{\boldmath$i$}u$.
Hence, $u\in \mathcal D^{1,0}$. 

\end{proof}

\par\noindent 
By definition a $CR$ structure
on an odd dimensional  manifold  consists
of the pair $(\mathop{\Null} \om,J)$ where 
$\om$ is a contact structure and $J$
is a complex structure on the contact subbundle $\mathop{\Null} \om$
(\ie $J$ is integrable).
In addition, the characteristic vector field 
$\xi$ for $\om$ is said to be a {\em characteristic $CR$-vector field}
if $\Ll_\xi J=0$. 
Consider $(\mathop{\Null} \om_a, \bar J_a)$ on $U$ $(a=1,2,3).$ 
By Lemma \ref{properties},
each $\xi_a$ is a characteristic
vector field for $\om_a$ on $U$.
From (3) of Lemma \ref{properties},
$\Ll_{\xi_\al}J_\al=0$. It is easy to check that
$\Ll_{\xi_a}\bar J_a=0$.

\begin{theorem}\label{crstructure1}
Each $\bar J_\al$ is integrable on $\mathcal{\Null}\om_\al$.
As a consequence, a nondegenerate quaternionic $CR$ structure 
$\{\om_\al,J_\al\}_{\al=1,2,3}$ on a neighborhood $U$
of $M^{4n+3}$ induces three nondegenerate
$CR$ structures
$(\mathop{\Null} \om_\al, \bar J_\al)$
 equipped with characteristic $CR$-vector field
$\xi_\al$ for each $\om_\al$ $(\al =1,2,3)$.
In fact, $\om_\al(\xi_\al)=1$ and $d\om_\al(\xi_\al,X)=0$ $(\forall\ X\in TM)$
$(\al=1,2,3)$.
\end{theorem}

\begin{proof}
Consider the case for $(\mathop{\Null} \om_1, \bar J_1)$.
Let $\displaystyle \mathop{\Null} \om_1\otimes{\CC}=T^{1,0}\oplus T^{0,1}$
where $\displaystyle T^{1,0}=
\mathcal D^{1,0}\oplus \{\xi_2-\mbox{\boldmath$i$}\xi_3\}$.
By Proposition \ref{complex-contact},
if $X,Y\in \mathcal D^{1,0}$, then there exist elements
$a\in\RR$ and $u\in \mathcal D^{1,0}$
such that
$[X,Y]=a(\xi_2-\mbox{\boldmath$i$}\xi_3)+u$.
By definition, 
\begin{equation}
\bar J_1[X,Y]=a\bar J_1(\xi_2-\mbox{\boldmath$i$}\xi_3)
+J_1u=a\mbox{\boldmath$i$}(\xi_2-\mbox{\boldmath$i$}\xi_3)
+\mbox{\boldmath$i$}u =\mbox{\boldmath$i$}[X,Y],
\end{equation}it follows $[X,Y]\in T^{1,0}$.
It suffices to show that the element
$[\xi_2-\mbox{\boldmath$i$}\xi_3,v]\in T^{1,0}$ for $v\in \mathcal D^{1,0}$.
As $\Ll_{\xi_2}J_1=-J_3$ and
$-J_3v=(\Ll_{\xi_2}J_1)v=\Ll_{\xi_2}(J_1v)-J_1(\Ll_{\xi_2}v)$,
\begin{equation}\label{J2}
\begin{split}
J_1(\Ll_{\xi_2}v)=J_3v+\mbox{\boldmath$i$}\Ll_{\xi_2}v.
\end{split}
\end{equation}
Note that $[\xi_2-\mbox{\boldmath$i$}\xi_3,v]=
\Ll_{\xi_2}v-\mbox{\boldmath$i$}\Ll_{\xi_3}v\in \mathcal D\otimes {\CC}$
on which $\bar J_a=J_a$.
Then $\bar J_1[\xi_2-\mbox{\boldmath$i$}\xi_3,v]
=J_1(\Ll_{\xi_2}v)-\mbox{\boldmath$i$}J_1(\Ll_{\xi_3}v)$.
Moreover, as
$J_2v=(\Ll_{\xi_3}J_1)v=\mbox{\boldmath$i$}\Ll_{\xi_3}(v)-J_1(\Ll_{\xi_3}v)$
and $J_2v=J_3J_1v=\mbox{\boldmath$i$}J_3v$,
it follows that $J_1(\Ll_{\xi_3}v)=
-\mbox{\boldmath$i$}J_3v+\mbox{\boldmath$i$}\Ll_{\xi_3}v$.
Using this equality and \eqref{J2}, it follows that
\begin{equation}
\begin{split}
\bar J_1[\xi_2-\mbox{\boldmath$i$}\xi_3,v]&=
J_1(\Ll_{\xi_2}v)-\mbox{\boldmath$i$}J_1(\Ll_{\xi_3}v)
=\mbox{\boldmath$i$}\Ll_{\xi_2}v+\Ll_{\xi_3}v\\
&=\mbox{\boldmath$i$}(\Ll_{\xi_2}v-\mbox{\boldmath$i$}\Ll_{\xi_3}v)
=\mbox{\boldmath$i$}[\xi_2-\mbox{\boldmath$i$}\xi_3,v].
\end{split}
\end{equation}
Therefore, $[T^{1,0},T^{1,0}]\subset T^{1,0}$ so that
$\bar J_1$ is a complex structure
on $\mathop{\Null} \om_1$, \ie
$(\mathop{\Null}\om_1, \bar J_1)$ is a $CR$ structure
on $U$. The same holds for 
$(\mathop{\Null} \om_b, \bar J_b)$ $(b=2,3)$.

\end{proof}

\section{Model of QCR space forms with type $(4p+3,4q)$}
\label{MCR}
Suppose that $p+q=n$.
Let ${\HH}^{n+1}$ be the quaternionc number space 
in quaternionic dimension $n+1$
with nondegenerate quaternionic Hermitian form
\begin{equation}
\langle x,y\rangle={\bar x}_1 y\sb 1+ \cdots +
 {\bar x}\sb {p+1} y \sb{p+1}-
{\bar x}_{p+2} y\sb{p+2}- \cdots -
 {\bar x}\sb {n+1} y \sb{n+1}.
\end{equation}
If we denote ${\rm Re}\langle x,y\rangle$ the real part of
$\langle x,y\rangle$, then it is noted that
${\rm Re}\langle \ ,\ \rangle$ is a nondegenerate symmetric bilinear form
on ${\HH}^{n+1}$.
In the quaternion case, the group of all invertible
matrices $\mathop{\GL}(n+1,{\HH})$
is acting from the left and ${\HH}^*=\mathop{\GL}(1,{\HH})$ 
acting as the scalar
multiplications from the right on ${\HH}^{n+1}$, which forms the
group $\displaystyle\mathop{\GL}(n+1,{\HH})\cdot \mathop{\GL}(1,{\HH})
=\mathop{\GL}(n+1,{\HH})\mathop{\times}_{\RR^*}^{}\mathop{\GL}(1,{\HH})$.
Let ${\rm Sp}(p+1,q)\cdot {\rm Sp}(1)$ 
be the subgroup of $\mathop{\GL}(n+1,{\HH})\cdot \mathop{\GL}(1,{\HH})$
whose elements preserve the 
nondegenerate bilinear form
${\rm Re}\langle \ ,\ \rangle$.
Denote by $\Sigma^{3+4p,4q}_{\HH}$ the $(4n+3)$-dimensional quadric space:
\begin{equation}
\begin{split}
&\{(z_1,\cdots,z_{p+1},w_1,\cdots,w_q)\in {\HH}^{n+1}\ | \\
& \ \ \ \ \ \ \ ||(z,w)||^2=|z_1|^2+\cdots+|z_{p+1}|^2-|w_1|^2-\cdots-|w_q|^2=1\}.
\end{split}\end{equation}In particular, the group 
${\rm Sp}(p+1,q)\cdot {\rm Sp}(1)$ leaves
$\Sigma^{3+4p,4q}_{\HH}$ invariant.
Let $\langle \ ,\ \rangle_x$ be the 
nondegenerate quaternionic inner product on the tangent space 
$T_x{\HH}^{n+1}$ obtained from the parallel translation of
$\langle \ ,\ \rangle$ to the point $x\in {\HH}^{n+1}$.
Recall that $\{I,\ J,\ K\}$ is the standard quaternionic structure
on ${\HH}^{n+1}$ which operates as $Iz=z\mbox{\boldmath$i$}$,
$Jz=z\mbox{\boldmath$j$}$, or $Kz=z\mbox{\boldmath$k$}$.
As usual, $\{I_x,\ J_x,\ K_x\}$ acts on $T_x{\HH}^{n+1}$ at each point $x$.
Then it is easy to see that $g^{\HH}_x(X,Y)={\rm Re}\langle X,Y\rangle_x$
 $(\forall\ X,Y\in T_x{\HH}^{n+1})$ is the standard pseudo-euclidean metric
of type $(p+1,q)$ 
on ${\HH}^{n+1}$ which is invariant under
$\{I,\ J,\ K\}$.
Restricted $g^{\HH}$ to the
quadric $\Sigma^{3+4p,4q}_{\HH}$ in ${\HH}^{n+1}$,
we obtain a nondegenerate pseudo-Riemannian
metric $g$ of type $(3+4p,4q)$ where $p+q=n$.

\begin{definition}(\cite{KU1})\label{fqcc}
 The quadric $\Sigma^{3+4p,4q}_{\HH}$ is referred to
the quaternionic pseudo-\\
Riemannian space form of type $(3+4p,4q)$ 
with constant curvature $1$ 
endowed with a transitive group
of isometries ${\rm Sp}(p+1,q)\cdot {\rm Sp}(1)$
for which $\displaystyle
\Sigma^{3+4p,4q}_{\HH}={\rm Sp}(p+1,q)\cdot {\rm Sp}(1)/
{\rm Sp}(p,q)\cdot {\rm Sp}(1)$
where ${\rm Sp}(p,q)\cdot {\rm Sp}(1)$ 
is the stabilizer at $(1,0,\cdots,0)$. 
\end{definition}Compare \cite{WO}, \cite{KU1}. 
When $(\Sigma^{3+4p,4q}_{\HH},g^{\HH})$ is viewed 
as a real pseudo-Riemannian space form,
the full group of isometries is ${\rm O}(4p+4,4q)$.
It is noted that the intersection of ${\rm O}(4p+4,4q)$
with $\mathop{\GL}(n+1,{\HH})\cdot \mathop{\GL}(1,{\HH})$ is
${\rm Sp}(p+1,q)\cdot {\rm Sp}(1)$.
When $N_x$ is the
 normal vector at $x\in \Sigma^{3+4p,4q}_{\HH}$,
$T_x\Sigma^{3+4p,4q}_{\HH}=N_x^{\perp}$ with respect to
$g^{\HH}$. If $N$ is
a normal vector field on $\Sigma^{3+4p,4q}_{\HH}$, then
$IN,JN,KN\in T\Sigma^{3+4p,4q}_{\HH}$ such that
there is the decomposition $\displaystyle
T\Sigma^{3+4p,4q}_{\HH}=\{IN,JN,KN\}\oplus \{IN,JN,KN\}^{\perp}$.
Let $\mathcal D=\{IN,JN,KN\}^{\perp}$ which is the
$4n$-dimensional subbundle.
As $g^{\HH}$ is a $\{I,\ J,\ K\}$ -invariant metric, $(\mathcal D,g|\mathcal D)$ is also invariant under
$\{I,\ J,\ K\}$.
Now, ${\rm Sp}(1)$ acts
freely on $\Sigma^{3+4p,4q}_{\HH}$
as (right translations): $(\lam\in {\rm Sp}(1))$
\begin{equation*}
(\lam,(z_1,\cdots,z_{p+1},w_1,\cdots,w_q))=
(z_1\cdot\bar\lam,\cdots,z_{p+1}\cdot\bar\lam,
w_1\cdot\bar\lam,\cdots,w_q\cdot\bar\lam).
\end{equation*}

\begin{definition}\label{p-projective}
Let ${\HH}{\PP}^{p,q}$ be the orbit space
$\Sigma^{3+4p,4q}_{\HH}/{\rm Sp}(1)$
which is said to be the 
quaternionic pseudo-K\"ahler projective space of type $(4p,4q)$.
\end{definition}
See Definition \ref{p-Kahler}
for the definition of quaternionic pseudo-K\"ahler  
manifold in general.
${\HH}{\PP}^{p,q}$ is shown to be a quaternionic 
pseudo-K\"ahler manifold in Theorem \ref{q-K-s}
provided that $4n\geq 8$.
When $p=n,q=0$, ${\HH}{\PP}^{n,0}$ is the standard
quaternionic projective space ${\HH}{\PP}^n$.
When $p=0,q=n$, ${\HH}{\PP}^{0,n}$ is the quaternionic 
hyperbolic space ${\HH}_{\HH}^n$.
It is easy to see that
${\HH}{\PP}^{p,q}$ is homotopic to the canonical quaternionic
line bundle over the quaternionic K\"ahler
projective space ${\HH}{\PP}^{p}$.
There is the equivariant principal bundle:
\begin{equation}\label{sp1-bundle}
{\rm Sp}(1)\ra  ({\rm Sp}(p+1,q)\cdot {\rm Sp}(1),\Sigma^{3+4p,4q}_{\HH})
\stackrel{\pi}\lra ({\rm PSp}(p+1,q),{\HH}{\PP}^{p,q})
\end{equation}

On the other hand, let
\begin{equation}\label{sphere- form}
\om_0=-(\bar z_1dz_1 +\cdots + \bar z_{p+1}dz_{p+1}
-\bar w_1dw_1 -\cdots - \bar w_{q}dw_{q}).
\end{equation}Then it is easy to check that
$\om_0$ is an ${\mathfrak s}{\mathfrak p}(1)$-valued $1$-form
on $\Sigma^{3+4p,4q}_{\HH}$.
Let $\xi_1,\xi_2,\xi_3$ be the vector fields
on $\Sigma^{3+4p,4q}_{\HH}$ induced by the one-parameter subgroups 
$\{e^{{\mbox{\boldmath$i$}}\theta}\}_{\theta\in {\RR}}$,
$\{e^{{\mbox{\boldmath$j$}}\theta}\}_{\theta\in {\RR}}$,
$\{e^{{\mbox{\boldmath$k$}}\theta}\}_{\theta\in {\RR}}$ respectively,
which is equivalent to that $\xi_1=IN, \xi_2=JN, \xi_3=KN$.
A calculation shows that
\begin{equation}\label{ijk}
\om_0(\xi_1)={\mbox{\boldmath$i$}},\
\om_0(\xi_2)={\mbox{\boldmath$j$}},\ 
\om_0(\xi_3)={\mbox{\boldmath$k$}}.
\end{equation}
By the formula of $\om_0$, if $a\in {\rm Sp}(1)$,
then the right translation ${\rm R}_a$ on $\Sigma^{3+4p,4q}_{\HH}$
satisfies that
\begin{equation}\label{rightaction}
{\rm R}_a^*\om_0=a\cdot \om_0\cdot \bar a.
\end{equation}
Therefore, $\om_0$ is a connection form of the above bundle
\eqref{sp1-bundle}. Note that
${\rm Sp}(p+1,q)$ leaves $\om_0$ invariant.
We shall check the conditions (i), (ii), (iii), (iv) and \eqref{localflow}
so that $(\Sigma^{3+4p,4q}_{\HH},\{I,J,K\},g,\om_0)$ will be a 
quaternionic $CR$ manifold.
First of all, it follows that
\begin{equation*}
\omega_0\wedge\omega_0\wedge\omega_0\wedge
\overbrace{(d\omega_0\wedge d\omega_0)\wedge 
\cdots\wedge(d\omega_0\wedge d\omega_0)}^{n\ \mbox{\scriptsize{times}}}
\neq 0 \ \ \mbox{at any point of}\ \Sigma^{3+4p,4q}_{\HH}.
\end{equation*}(Compare \cite{KA},\cite{SA} for example).
In fact, letting $\om_0=\omega_1{\mbox{\boldmath$i$}}+
\omega_2{\mbox{\boldmath$j$}}+\omega_3{\mbox{\boldmath$k$}}$ as before,
\[
{\omega_0}^3\wedge d{\omega_0}^{2n}=
6\omega_1\wedge\omega_2\wedge\omega_3\wedge({d\omega_1}^{2}+
{d\omega_2}^{2}+{d\omega_3}^{2})^n.
\]This calculation shows (iii).
In particular, each $\om_a$ is a nondegenerate
contact form on $\Sigma^{3+4p,4q}_{\HH}$. 
Using \eqref{rightaction} and as
$\xi_1$ generates $\{e^{{\mbox{\boldmath$i$}}\theta}\}_{\theta\in \RR}
\subset{\rm Sp}(1)$,
$\Ll_{\xi_1}\om_1=0$. (Similarly we have
$\Ll_{\xi_2}\om_2=\Ll_{\xi_3}\om_3=0$.)
Noting that $\om_a(\xi_a)=1$
and $0=\Ll_{\xi_a}\om_a=\iota_{\xi_a}\circ d\om_a$
from \eqref{ijk},
each $\xi_a$ is the characteristic vector field
for $\om_a$.
Moreover, note that
$\{\xi_1,\xi_2,\xi_3\}$ generates 
the fields of Lie algebra of ${\rm Sp}(1)$.
It follows that
 $\displaystyle \mathcal D=\mathop{\cap}_{a=1}^{3}{\Null} \om_a$
for which there is the decomposition
$T\Sigma^{3+4p,4q}_{\HH}=\{\xi_1,\xi_2,\xi_3\}\oplus \mathcal D$.
If $\{e_i\}_{i=1,\cdots,4n}$ is the orthonormal basis
of $\mathcal D$, then the dual frame $\se^i$ is obtained as 
$\se^i(e_j)=\delta^i_j$ and 
$\se^i(\xi_1)=\se^i(\xi_2)=\se^i(\xi_3)=0$. 
In order to prove that the distribution uniquely determined 
by \eqref{localflow} are 
$\{\xi_1,\xi_2,\xi_3\}$ (\cf \eqref{3-qcr} also),
we need the following lemma.
\begin{lemma}
\begin{equation*}
d\om_1(X,Y)=g(X,IY),\ d\om_2(X,Y)=g(X,JY),\
d\om_3(X,Y)=g(X,KY)
\end{equation*}where $X,Y\in \mathcal D$.
\end{lemma}
\begin{proof}
Given $X,Y\in \mathcal D_x$,
let $u,v$ be the vectors at the origin by parallel translation of $X,Y$ at
$x\in \Sigma^{3+4p,4q}_{\HH}$ respectively.
Then by definition, $g(X,Y)={\rm Re}\langle u,v\rangle$.
Furthermore, 
\begin{equation}\label{g-part}
g(X,IY)={\rm Re}(\langle u,v\cdot{\mbox{\boldmath$i$}}\rangle)
={\rm Re}(\langle u,v\rangle\cdot{\mbox{\boldmath$i$}}).
\end{equation}
From \eqref{sphere- form},
if $X,Y\in \mathcal D_x$, then 
\begin{equation}
\begin{split}
d\om_0(X,Y)=-(d\bar z_1\we dz_1 +\cdots + d\bar z_{p+1}\we dz_{p+1}
-d\bar w_1\we dw_1 -\cdots - d\bar w_{q}\we dw_{q})(u,v).
\end{split}\end{equation}
Then a calculation shows that
$\displaystyle d\om_0(X,Y)
=-\frac 12(\langle u,v\rangle-\overline{\langle u,v\rangle})$.
It is easy to check that
the ${\mbox{\boldmath$i$}}$-part of 
$\displaystyle-\frac 12(\langle u,v\rangle-\overline{\langle u,v\rangle})$
is  ${\rm Re}(\langle u,v\rangle\cdot{\mbox{\boldmath$i$}})$. 
Since $d\om_1(X,Y)$ is
the ${\mbox{\boldmath$i$}}$-part of $d\om(X,Y)$ and by \eqref{g-part},
we obtain the equality $\displaystyle g(X,IY)=d\om_1(X,Y)$.
Similarly, we have that
$\displaystyle g(X,JY)=d\om_2(X,Y),\ g(X,KY)=d\om_3(X,Y)$.
\end{proof}

From this lemma, 
$d\om_a(e_i,e_j)=
g(e_i,J_a e_j)=-{\bf J^a}_{ij}$.
Since $\{\xi_1,\xi_2,\xi_3\}$ generates 
${\rm Sp}(1)$ of the bundle \eqref{sp1-bundle},
we obtain $\displaystyle
d\omega_a+2\omega_b\wedge\omega_c=
-{\bf J^a}_{ij}\theta^i\wedge\theta^j$.
Applying to $J,K$ similarly,
we obtain the following structure equation of the bundle \eqref{sp1-bundle}:
\begin{equation}\label{localflow1}
d\omega_0+\omega_0\wedge\omega_0=
-({\bf I}_{ij}\mbox{\boldmath$i$}+{\bf J}_{ij}\mbox{\boldmath$j$}
+{\bf K}_{ij}\mbox{\boldmath$k$})\theta^i\wedge\theta^j.
\end{equation}From this equation,
the condition \eqref{localflow} is easily checked so that
$\mathop{\Null}\om_\al=\{\xi_1,\xi_2,\xi_3\}$.
We summarize that
\begin{theorem}\label{spaceform}
$(\Sigma^{3+4p,4q}_{\HH},\{\om_a\}_{a=1,2,3},\{I,J,K\},g)$
is a $(4n+3)$-dimensional homogeneous quaternionic $CR$
manifold of type $(3+4p,4q)$
where $p+q=n\geq 0$. Moreover,
there exists the equivariant principal bundle
\eqref{sp1-bundle} of the pseudo-Riemannian submersion over
the homogeneous quaternionic pseudo-K\"ahler projective space 
${\HH}{\PP}^{p,q}$ of type $(4p,4q)$:
$\displaystyle
{\rm Sp}(1)\ra ({\rm Sp}(p+1,q)\cdot{\rm Sp}(1),
\Sigma^{3+4p,4q}_{\HH},g)\stackrel{\pi}{\lra} ({\rm PSp}(p+1,q),
{\HH}{\PP}^{p,q},\hat g)$.
\end{theorem}
We shall prove more generally in Theorem \ref{q-K-s}
that $({\rm PSp}(p+1,q),{\HH}{\PP}^{4p,4q})$
supports an invariant quaternionic pseudo-K\"ahler 
metric $\hat g$ of type $(4p,4q)$. 

\begin{remark}
${\bf (a)}$\ 
In \cite{AK}, it is shown 
that $(\Sigma^{3+4p,4q}_{\HH},\{I,J,K\},g)$
is a pseudo-Sasakian space form of constant positive curvature
with type $(4p+3,4q)$.
\vskip0.1cm
\noindent ${\bf (b)}$\ 
When $q=0$ or $p=0$,
we can find discrete cocompact subgroups
from ${\rm Sp}(n+1)\cdot{\rm Sp}(1)$
or ${\rm Sp}(1,n)\cdot{\rm Sp}(1)$
that act properly and freely
on $\Sigma^{3+4n,0}_{\HH}=S^{4n+3}$ or
$\Sigma^{3,4n}_{\HH}=V_{-1}^{4n+3}$ respectively.
Thus, we obtain compact nondegenerate
quaternionic $CR$ manifolds.
In fact, 
${\bf (i)}$ The spherical space form
$S^{4n+3}/F$
which is ${\rm Sp}(1)$ or ${\rm SO}(3)$-bundle over 
the quaternionic K\"ahler projective orbifold 
${\HH}{\PP}^{n}/F^*$ of positive scalar curvature.
$(F\subset{\rm Sp}(n+1)\cdot {\rm Sp}(1)$ is 
a finite group$.)$ \
${\bf (ii)}$ The pseudo-Riemannian standard space form 
$V_{-1}^{4n+3}/\Gamma$ of type $(4n,3)$
with constant sectional curvature $-1$
which is an ${\rm Sp}(1)$-bundle
over the quaternionic K\"ahler hyperbolic
 orbifold ${\HH}_{\HH}^n/\Gamma^*$
of negative scalar curvature.
As we know, there exists no compact pseudo-Sasakian 
manifold (or quaternionic $CR$ manifold)
whose pseudo-K\"ahler orbifold is of zero Ricci curvature.
However in our case an indefnite Heisenberg nilmanifold is a compact {\em pseudo-conformal quaternionic
$CR$} manifold whose
pseudo-K\"ahler orbifold is of zero Ricci curvature, see
$\S \ref{pcHei}$.

\end{remark}

\section{Local Principal bundle}\label{principal}
Let $\{e_i\}_{i=1,\cdots,4n}$ be the basis of $\mathcal D|U$
such that $g^{\mathcal D}(e_i,e_j)=g_{ij}$.
We choose a local coframe $\se^i$ for which 
\begin{equation}\label{coframing}
\se^i|V=0\ \ \mbox{and}\ \ \theta^i(e_j)=\delta_{ij}.
\end{equation}
As usual the quaternionic
structure $\{J_\al\}_{\al=1,2,3}$ can be represented locally by the matrix 
${\bf J^\al}_{i}^{j}$ such as $\displaystyle J_\al e_i={\bf J^\al}_{i}^{j}e_j$.
Using \eqref{1trans}, note that
$\displaystyle \rho_\al(e_j,e_i)={\bf J^\al}_{i}^{k}g_{jk}={\bf J^\al}_{ij}$.
Here the matrix $(g_{ij})$ lowers and raises the indices. Then
we can write the structure equation \eqref{integraleq1}
by using $\se^i$:

\begin{equation}\label{newform}
d\om_\al+2{\om_\be\we\om_\ga}=-{\bf J^\al}_{ij}{\theta}^i\wedge{\theta}^j\
 \ (\al=1,2,3).
\end{equation}

Using $\om$ of $(iii)$, we have the following formula
corresponding to \eqref{newform}:
\begin{equation}\label{3-qcr}
d\om+\om\we\om=-({\bf J^1}_{ij}\mbox{\boldmath$i$}+{\bf J^2}_{ij}\mbox{\boldmath$j$}+{\bf J^3}_{ij}\mbox{\boldmath$k$})
\se^i\we\se^j.
\end{equation} 

Denote by $\mathcal E$ the local transformation groups generated by
$V$ acting on a small neighborhood $U'$ of $U$.
As $\mathcal E$ is locally isomorphic to the compact Lie
group ${\rm SO}(3)$ by Lemma \ref{so3},
it acts properly on $U'$. (See for example \cite{PA}.)
If we note that each $\xi_a$ is a nonzero vector field everywhere
on $U$, then the stabilizer of $\mathcal E$ is finite at every point.
By the slice theorem of compact Lie groups \cite{BR},
choosing a sufficiently small neighborhood $\mathcal E'$ 
of the identity from $\mathcal E$, $\mathcal E'$ 
acts properly and freely on $U'$.
We choose such $U'$ (respectively $\mathcal E'$)
from the beginning and replace it by $U$ (respectively $\mathcal E$).
Then there is a principal local fibration:
\begin{equation}\label{prin}
\mathcal E\ra U\stackrel{\pi}{\lra} U/\mathcal E.
\end{equation}
If we note that $V\oplus \mathcal D=TM|U$, $\pi$ maps $\mathcal D$
 isomorphically onto  $T(U/\mathcal E)$ at each point of $U$.
So $\{\pi_*e_i \ |\ i=1,\cdots,4n\}$ is a basis of
$T(U/\mathcal E)$ at each point of $U/\mathcal E$.
Let $\hat \se^i$ be the dual frame on
$U/\mathcal E$ such that 
\begin{equation}\label{dual-co}
\hat\se^i(\pi_*e_j)=\delta_{ij}\ \ \mbox{on}\  U/\mathcal E.
\end{equation}Since $\se^i$ is the coframe
of $\{e_i\}$ and  $\pi^*\hat\se^i|V=\se^i|V=0$,
it follows that 
\begin{equation}\label{pull}
\pi^*\hat\se^i=\se^i\ \mbox{on}\  U\ (i=1,\cdots,4n).
\end{equation}

\begin{lemma}\label{intersec}
Put $J_1=I,\ J_2=J,\ J_3=K$ respectively.
Let $\{\f_\se\}_{-\varepsilon<\se<\varepsilon}$ be a 
local one-parameter subgroup
of the local group $\mathcal E$.
Then there exists an element $G_{\se}\in {\rm SO}(3)$ 
satisfying the following:
\begin{equation}\label{3xi}
\begin{split}
&(1)\ \ \ (\f_\se)_*\left(\begin{array}{c}
\xi_1\\
\xi_2\\
\xi_3\end{array}\right)=
G_{\se}\left(\begin{array}{c}
\xi_1\\
\xi_2\\
\xi_3\end{array}\right).\\
& \\
&(2)\ \  \
\left(\begin{array}{c}
I_{\f_\se y}\\
J_{\f_\se y}\\
K_{\f_\se y}\end{array}\right)\circ {\f_\se}_*=
{\f_\se}_*\circ {}^tG(\se)
\left(\begin{array}{c}
I_y\\
J_y\\
K_y\end{array}\right).
\end{split}
\end{equation}
\end{lemma}

\begin{proof}
Since every leaf of $V$ is locally isomorphic to 
${\rm SO}(3)$, $\xi_a$ is viewed as the fundamental vector field
to the principal fibration $\pi:U\ra U/\mathcal E$.
Thus we may assume that $\xi_1,\xi_2,\xi_3$
 correspond to $\mbox{\boldmath$i$},
\mbox{\boldmath$j$}, \mbox{\boldmath$k$}$ respectively
so that $\displaystyle \f_\se^1=e^{{\mbox{\boldmath$i$}}\se},\
\f_\se^2=e^{{\mbox{\boldmath$j$}}\se},\ \f_\se^3=
e^{{\mbox{\boldmath$k$}}\se}$ up to conjugacy by an element of ${\rm SO}(3)$,
A calculation shows that
$(\f_\se^1)_*((\xi_2)_x)=
\cos 2\se\cdot(\xi_2)_{\f_\se^1x}+\sin 2\se\cdot(\xi_3)_{\f_\se^1x}$.
Similarly, $(\f_\se^1)_*((\xi_3)_x)=
-\sin 2\se\cdot(\xi_2)_{\f_\se^1x}+\cos 2\se\cdot(\xi_3)_{\f_\se^1x}$,
$(\f_\se^1)_*((\xi_1)_x)=(\xi_1)_{\f_\se^1x}$.
This holds similarly for $\f_\se^1,\f_\se^2$.
It turns out that if $\f_\se\in \mathcal E$, then
there exists an element $G_{\se}\in {\rm SO}(3)$ which shows
the above formula $(1)$.
Since $\f_t$ preserves $\mathcal D$ $(-\varepsilon<t<\varepsilon)$,
using $(1)$ we see that
\begin{equation}\label{omega-con}
\f_t^*(\om_1,\om_2,\om_3)=(\om_1,\om_2,\om_3)G_{t}.
\end{equation}
Since there exists an element
$g_t\in {\rm Sp}(1)$ such that
$\displaystyle g_t\left(\begin{array}{c}
\mbox{\boldmath$i$}\\
\mbox{\boldmath$j$}\\
\mbox{\boldmath$k$}
\end{array}\right)\bar g_t=
G_{t}\left(\begin{array}{c}
\mbox{\boldmath$i$}\\
\mbox{\boldmath$j$}\\
\mbox{\boldmath$k$}
\end{array}\right)$ ($\bar g_t$ is the quaternion conjugate of
$g_t$), \eqref{omega-con} is equivalent with
\begin{equation}\label{conj}
\f_t^*\om=g_t\cdot \om\cdot \bar g_t.
\end{equation}
Differentiate this equation which yields that
\begin{equation}\label{conj1}
\f_t^*(d\om+\om\we\om)\equiv  g_t(d\om+\om\we\om)\bar g_t\ \ \mbox{mod}\ \om.
\end{equation}
Using the equation \eqref{newform}, it follows that
\begin{equation*}
\begin{split}
\f^*_t(({I}_{ij}, {J}_{ij}, {K}_{ij})\left(
\begin{array}{c}
\mbox{\boldmath$i$}\\
\mbox{\boldmath$j$}\\
\mbox{\boldmath$k$}\end{array}\right)
\theta^i\wedge\theta^j)&\equiv
({I}_{ij}, {J}_{ij}, {K}_{ij})g_t\left(
\begin{array}{c}
\mbox{\boldmath$i$}\\
\mbox{\boldmath$j$}\\
\mbox{\boldmath$k$}\end{array}\right)\bar g_t
\theta^i\wedge\theta^j\\
&=
({I}_{ij}, {J}_{ij}, {K}_{ij})G_{t}\left(
\begin{array}{c}
\mbox{\boldmath$i$}\\
\mbox{\boldmath$j$}\\
\mbox{\boldmath$k$}\end{array}\right)
\theta^i\wedge\theta^j.
\end{split}
\end{equation*}
Noting that $\f^*_t\se^i=\f^*_t(\pi^*\hat \se^i)=\se^i$,
the above equation implies that
\begin{equation}\label{intsec1}
\begin{split}
({I}_{ij}(\f_t(x)), {J}_{ij}(\f_t(x)), {K}_{ij}(\f_t(x)))
\equiv ({I}_{ij}(x), {J}_{ij}(x), {K}_{ij}(x))G_t(x)
 \ \mbox{mod} \ \om.
\end{split}
\end{equation}Since 
$\pi_*{\f_t}_*((e_i)_x)={\pi}_*((e_i)_{\f_t x})$ $(x\in U)$, it follows
${\f_t}_*((e_i)_x)=(e_i)_{\f_tx}$.
Letting  $G_{t}=(s_{ij})\in{\rm SO}(3)$ and using \eqref{intsec1},
\begin{equation*}
\begin{split}
&I_{\f_tx}(\f_t)_*((e_i)_x)=I_{\f_tx}((e_i)_{\f_tx})=
{I}_{i}^{j}(\f_tx)((e_j)_{\f_tx})\\
&=({I}_{i}^{j}(x)\cdot s_{11}+{J}_{i}^{j}(x)\cdot s_{21}+
 {K}_{i}^{j}(x)\cdot s_{31}))((\f_t)_*((e_j)_x))\\
&=(\f_t)_*( s_{11}\cdot I_{x}((e_i)_x)+s_{21}\cdot J_{x}((e_i)_x)+
s_{31}\cdot K_{x}((e_i)_x))\\
&=(\f_t)_*((s_{11}, s_{21}, s_{31})\left(
\begin{array}{c}
I_x\\
J_x\\
K_x
\end{array}\right)(e_i)_x).
\end{split}
\end{equation*}
The  same argument applies to $J_{\f_tx},K_{\f_tx}$ to
conclude that
$\displaystyle
\left(\begin{array}{c}
I_{\f_t x}\\
J_{\f_t x}\\
K_{\f_t x}\end{array}\right)\circ {\f_t}_*=
{\f_t}_*\circ {}^tG_t
\left(\begin{array}{c}
I_x\\
J_x\\
K_x\end{array}\right)$.
This proves (2).
\end{proof}

\begin{lemma}\label{quat} 
The quaternionic structure
$\{I,J,K\}$ on $\mathcal D|U$ induces a family of
quaternionic structures
$\{\hat I_i,\hat J_i, \hat K_i\}_{i\in \Lambda}$ on $U/\mathcal E$.
\end{lemma}
 
\begin{proof}
Choose a small neighborhood $V_i\subset U/\mathcal E$ and a
section $s_i:V_i\ra U$ for the principal bundle $\pi:
U\ra U/\mathcal E$. Let $\hat x\in V_i$ and a vector 
$\hat X_{\hat x}\in T V_i$. Choose a vector $X_{s_i(\hat x)}\in
\mathcal D_{s_i(\hat x)}$
such that $\pi_*(X_{s_i(\hat x)})=\hat X_{\hat x}$.
Define endomorphisms 
$\hat I_i,\hat J_i, \hat K_i$ on $V_i$ to be
\begin{equation}\label{3complex}
\begin{split}
(\hat I_i)_{\hat x}(\hat X_{\hat x})&=
\pi_*I_{s_i(\hat x)}X_{s_i(\hat x)},\\
(\hat J_i)_{\hat x}(\hat X_{\hat x})&=
\pi_*J_{s_i(\hat x)}X_{s_i(\hat x)},\\
(\hat K_i)_{\hat x}(\hat X_{\hat x})&=
\pi_*K_{s_i(\hat x)}X_{s_i(\hat x)}.
\end{split}
\end{equation}Since
$\pi_*:\mathcal D_{s_i(\hat x)}\ra T_{\hat x}(U/\mathcal E)$
is an isomorphism, $\hat I_i,\hat J_i, \hat K_i$ are well-defined
almost complex structures on $V_i$.
So we have a family $\{\hat I_i,\hat J_i, \hat K_i\}_{i\in \Lambda}$ 
of almost complex structures associated to an open
 cover $\{V_i\}_{i\in\Lambda}$ of $U/\mathcal E$.
Suppose that $V_i\cap V_j\neq \emptyset$.
If $\hat x\in V_i\cap V_j$, then there is an element $\f_{\se}\in \mathcal E$
such that $s_j(\hat x)=\f_{\se}\cdot s_i(\hat x)$.
As $\f_\se$ preserves $\mathcal D$, ${\f_\se}_*X_{s_i(\hat x)}\in 
\mathcal D_{s_j(\hat x)}$
and $\pi_*({\f_\se}_*X_{s_i(\hat x)})=\hat X_{\hat x}$. Then
\begin{equation}\label{proj}
X_{s_j(\hat x)}={\f_\se}_*X_{s_i(\hat x)}.
\end{equation}

Let  $\{{\hat I}_j,{\hat J}_j, {\hat K}_j\}$ 
be almost complex structures on $V_j$ obtained from \eqref{3complex}.
Using Lemma \ref{intersec} and \eqref{proj}, calculate
at $s_j(\hat x)\ (\hat x\in V_i\cap V_j)$,
\begin{equation*}
\begin{split}
\left(\begin{array}{c}
({\hat I}_j)_{\hat x}\\
({\hat J}_j)_{\hat x}\\
({\hat K}_j)_{\hat x}\end{array}\right)
\hat X_{\hat x}&=
\pi_*\left(\begin{array}{c}
I_{s_j(\hat x)}\\
J_{s_j(\hat x)}\\
K_{s_j(\hat x)}\end{array}\right)
X_{s_j(\hat x)}=
\pi_*\left(\begin{array}{c}
I_{\f_{\se}\cdot s_i(\hat x)}\\
J_{\f_{\se}\cdot s_i(\hat x)}\\
K_{\f_{\se}\cdot s_i(\hat x)}\end{array}\right)
{\f_\se}_*X_{s_i(\hat x)}\\
&=
\pi_*{\f_\se}_*\circ {}^t G_{\se}\left(\begin{array}{c}
I_{s_i(\hat x)}\\
J_{s_i(\hat x)}\\
K_{s_i(\hat x)}\end{array}\right)
X_{s_i(\hat x)}\\
&={}^tG(\se)\pi_*\left(\begin{array}{c}
I_{s_i(\hat x)}\\
J_{s_i(\hat x)}\\
K_{s_i(\hat x)}\end{array}\right)
X_{s_i(\hat x)}
={}^tG_{\se}
\left(\begin{array}{c}
(\hat I_i)_{\hat x}\\
(\hat J_i)_{\hat x}\\
(\hat K_i)_{\hat x}\end{array}\right)
\hat X_{\hat x},
\end{split}
\end{equation*}hence
$\displaystyle
\left(\begin{array}{c}
(\hat I_j)_{\hat x}\\
(\hat J_j)_{\hat x}\\
(\hat K_j)_{\hat x}\end{array}\right)
={}^tG_{\se}
\left(\begin{array}{c}
(\hat I_i)_{\hat x}\\
(\hat J_i)_{\hat x}\\
(\hat K_i)_{\hat x}\end{array}\right)$
on $\hat x\in V_i\cap V_j$.
Thus,
$\{\hat I_i,\hat J_i, \hat K_i\}_{i\in \Lambda}$ 
defines a quaternionic structure on $U/\mathcal E$.
\end{proof}

\subsection{Pseudo-Sasakian $3$-structure and Pseudo-K\"ahler structure}
\label{PS3}
We now take $\{e_i\}_{i=1,\cdots,4n}$ of $\mathcal D|U$ as the orthonormal
 basis, \ie $g_{ij}=\delta_{ij}$.
Then the bilinear form 
$\displaystyle g^{\mathcal D}=\mathop{\sum}_{i=1}^{4p}\se^i\cdot\se^i
-\mathop{\sum}_{i=4p+1}^{4n}\se^i\cdot\se^i$ defined on $\mathcal D$
induces a pseudo-Riemannian
metric on $U/\mathcal E$:
\begin{equation}\label{b-form}
\hat g=\mathop{\sum}_{i=1}^{4p}\hat\se^i\cdot\hat\se^i
-\mathop{\sum}_{i=4p+1}^{4n}\hat\se^i\cdot\hat\se^i
\end{equation}such that $g^{\mathcal D}=\pi^*\hat g$.
Let $\hat\nabla$ be the covariant derivative on $U/\mathcal E$.
If $\hat\om_j^i$ is the Levi-Civita
connection with respect to $\hat g$, then
it follows that
\begin{equation}\label{nab}
\hat\nabla\hat e_i=\hat\om^j_i\hat e_j.
\end{equation}
Then $\hat\om_j^i$ satisfies that
\begin{equation}\label{levi}
d\hat\se^i=\hat\se^j\we\hat\om^i_j,\ \  \hat\om_{ij}+\hat\om_{ji}=0.
\end{equation}
Put
\begin{equation}\label{Rcur}
\hat\Omega_j^i=d\hat\om^i_j-\hat\om^\sigma_j\we\hat\om^i_\sigma=
\frac 12\hat R^i_{jkl}\hat\se^k\we\hat\se^\ell.
\end{equation}

Consider the  following pseudo-Riemannian metric on $U$:
\begin{equation}\label{3sasakian meric}
\begin{split}
{\tilde g}_x(X,Y)&=\mathop{\sum}_{a=1}^3\om_a(X)\cdot\om_a(Y)+
\hat g_{\pi(x)}(\pi_*X,\pi_*Y)\ \ \ (X,Y\in T_xU). \\
(\mbox{Equivalently} &\ 
\tilde g=\mathop{\sum}_{a=1}^3\om_a\cdot\om_a+\mathop{\sum}_{i=1}^{4p}\se^i\cdot\se^i-\mathop{\sum}_{i=4p+1}^{4n}\se^i\cdot\se^i.)
\end{split}
\end{equation}
 Then we have shown in \cite{AK_1}
that the local principal fibration 
$\mathcal E\ra (U,\tilde g)\stackrel{\pi}\lra (U/\mathcal E,\hat g)$
is a pseudo-Sasakian $3$-structure.
In fact the following equation 
\eqref{sasaki-condition}
is equivalent with 
the normality condition
of the pseudo-Sasakian $3$-structure.
(Compare \cite{TA}, \cite{BL}.)

\begin{pro}\label{normal-sasaki}
Let $(\{\om_\al\},\{J_\al\},\{\xi_\al\})_{\al =1,2,3}$ be a
nondegenerate quaternionic $CR$ structure on
$U$ of a $(4n+3)$-manifold $M$.
If $\nabla$ is the Levi-Civita connection on $(U,\tilde g)$,
then,
\begin{equation}\label{sasaki-condition}
(\nabla_X\bar J_\al)Y={\tilde g}(X,Y)\xi_\al-\om_\al(Y)X\ \ (\al=1,2,3).
\end{equation}
\end{pro}

\begin{proof} For $X,Y\in TU$, consider the following tensor
\begin{equation}\label{sa-nije}
N^{\om_\al}(X,Y)=N(X,Y)+(X\om_\al(Y)-Y\om_\al(X))\xi_{\al}
\end{equation}where $N(X,Y)=
[\bar J_\al X,\bar J_\al Y]-[X,Y]-\bar J_\al[\bar J_\al X,Y]-
\bar J_\al[X,\bar J_\al Y]$ is the Nijenhuis torsion of 
$\bar J_\al$ $(\al=1,2,3)$.
A direct calculation for a contact metric structure
$\tilde g$ (\cf \cite{BL}) shows that
\begin{equation*}\label{normalityformula}
\begin{split}
  2{\tilde g}((\nabla_X{\bar J}_\al)Y,Z) &
  ={\tilde g}(N^{\om_\al}(Y,Z),\bar J_\al X)+
(\Ll_{{\bar J_\al}X}\om_\al)(Y)\\
& \ \ -(\Ll_{{\bar J_\al}Y}\om_\al)(X)+2{\tilde g}(X,Y)\om_\al(Z)
-2{\tilde g}(X,Z)\om_\al(Y).
\end{split}
\end{equation*}
We have shown in Theorem \ref{crstructure1}
that each $\bar J_\al$ is integrable on $\mathcal{\Null}\om_\al$.
It follows from the famous theorem that
the Nijenhuis torsion of $\bar J_\al$,
$N(X,Y)=0$ $(\forall\ X,Y\in \mathop{\Null}\om_\al)$.
By the formula \eqref{sa-nije},
$N^{\om_\al}(X,Y)=0$ for $\forall\ X,Y\in \mathop{\Null}\om_\al$.
To obtain \eqref{sasaki-condition},
it has to show that $N^{\om_\al}(X,Y)$
vanishes for all $X,Y\in TU$.
Noting the decomposition
$TU=\{\xi_1\}\oplus \mathop{\Null} \om_1$,
it suffices to show that $N^{\om_1}(\xi_1,X)=0$ (similarly for $\al=2,3$).
Since $\xi_\al$ is a characteristic $CR$-vector field 
for $(\om_\al,\bar J_\al)$ $(\al=1,2,3)$, \ie $\Ll_{\xi_1} \bar J_1=0$,
it follows that
$\bar J_1[\xi_1,Y]=[\xi_1,\bar J_1Y]$ $(\forall\ Y\in 
\mathop{\Null}\om_1)$.
In particular, $\displaystyle \bar J_1[\xi_1,\bar J_1 X]=-[\xi_1,X]$.
Hence, $N^{\om_\al}(\xi_1,X)=0$.
As a consequence, we see that
$N^{\om_\al}(X,Y)=0$  $(\forall\ X,Y\in TU)$.
On the other hand, 
if $N^{\om_\al}(X,Y)=0$  $(\forall\ X,Y\in TU)$,
then it is easy to see that
$(\Ll_{{\bar J_\al}X}\om_\al)(Y)-(\Ll_{{\bar J_\al}Y}\om_\al)(X)=0$.
 (See \cite{BL}.) From \eqref{3sasakian meric},
 note that $\om_\al(X)={\tilde g}(\xi_\al,X)$.
 The above equation \eqref{sasaki-condition} follows.
\end{proof}

As $\{\om_\al,\se^i\}_{\al=1,2,3;i=1\cdots 4n}$
are orthonormal coframes for the pseudo-Sasakian metric
$\tilde g$ (\cf \eqref{3sasakian meric}), the structure equation says
that there exist unique $1$-forms
$\f^i_j$, $\tau^i_\al$ $(i,j=1,\cdots,4n; \al=1,2,3)$ satisfying:
\begin{equation}\label{sasaki structure eq}
d\se^i=\se^j\wedge\f^i_j\ +\mathop{\sum}_{\al=1}^{3} \om_\al\we\tau^i_\al\
\ \ \  (\f_{ij}+\f_{ji}=0).
\end{equation} 
Then the normality condition 
for the pseudo-Sasakian $3$-structure is reinterpreted as
the following structure equation.

\begin{theorem}\label{normal-sasaki1}
There exsists a connection form $\{\om_j^i\}$
such that
\begin{equation}\label{3-structure eqs}
\begin{split}
d\bar {\bf J}^a_{ij}-\om^{\sigma}_{i}\bar {\bf J}^a_{\sigma j}-
\bar {\bf J}^a_{i\sigma}\om^{\sigma}_j=
2\bar {\bf J}^b_{ij}\cdot\om_c-2\bar {\bf J}^c_{ij}\cdot\om_b
\ \ \ ((a,b,c)\sim (1,2,3)).
\end{split}
\end{equation} 
\end{theorem}

\begin{proof}
It follows from Proposition \ref{normal-sasaki}
that $(\nabla_X{\bar J}_a)e_i=\tilde g(X,e_i)\xi_a$ 
for the orthonormal basis $e_i\in B$.
From the structure equation \eqref{sasaki structure eq}, 
let $\displaystyle
\nabla _Xe_i=\f_i^j(X)e_j+\mathop{\sum}_{b=1}^{3}{(\tau_b)}_i\xi_b$
which is substituted into the equation $\displaystyle (\nabla_X\bar J_a)e_i=
\nabla_X(\bar J_a e_i)-\bar J_{a}(\nabla_X e_i)$:
\begin{equation}
\begin{split}
(\nabla_X\bar J_a)e_i&=
(d(\bar {\bf J}^a)_i^\ell(X)-\f^{\sigma}_{i}(X)
(\bar {\bf J}^a)_{\sigma}^{\ell}
+(\bar {\bf J}^a)_i^{\sigma}\f_{\sigma}^{\ell}(X))e_{\ell}\\
&\ \ \ \ +\mathop{\sum}_{b=1}^{3}(\bar {\bf J}^a)_i^{\ell}
{(\tau_b)}_{\ell}(X)\xi_b
-\mathop{\sum}_{c\neq a}{(\tau_b)}_i(X)\xi_c\\
&=\tilde g(X,e_i)\xi_a.\ \ \ (\mbox{Here}\  \bar J_a\xi_b=\xi_c.)
\end{split}
\end{equation}
As $\tilde g(X,e_i)=\tilde g_{ki}\se^k(X)$
(\cf \eqref{3sasakian meric}),
this implies that
$d(\bar {\bf J}^a)_i^{\ell}-\f^{\sigma}_{i}(\bar {\bf J}^a)_{\sigma}^{\ell}
+(\bar {\bf J}^a)_i^{\sigma}\f_{\sigma}^\ell=0$ and
$\displaystyle (\bar {\bf J}^a)_i^{\ell}{(\tau_a)}_{\ell}(X)\xi_a=
\tilde g_{ki}\se^k(X)\xi_a$. It follows
$-{(\tau_a)}_{i}=(\bar {\bf J}^a)_{ij}\se^j$.
Then $\displaystyle {(\tau_a)}_{i}\tilde g^{ik}=-(\bar {\bf J}^a)_{ij}
\tilde g^{ik}\se^j=(\bar {\bf J}^a)_{ji}\tilde g^{ik}\se^j$, so that
$\displaystyle {(\tau_a)}^{i}=(\bar {\bf J}^a)_{j}^{i}\se^j$.
As $\tilde g_{ij}=\pm\delta_{ij}$, use $\tilde g^{ij}$ to lower the above equations:
\begin{equation}\label{upper-equations}
\begin{split}
&d(\bar {\bf J}^a)_{ij}-\f^{\sigma}_{i}(\bar {\bf J}^a)_{\sigma j}
-(\bar {\bf J}^a)_{i\sigma}\f_j^{\sigma}=0.\\
&{(\tau_a)}^i=(\bar {\bf J}^a)^i_{j}\se^j.
\end{split}
\end{equation}
Putting 
\begin{equation}\label{new-con}
\om_j^i=\f^i_j-\mathop{\sum}_{a=1}^{3}({\bar{\bf J}}^a)_{j}^i\om_a,
\end{equation}
the equation \eqref{sasaki structure eq} reduces to 
\begin{equation}\label{reduce}
d\se^i=\se^j\wedge\om^i_j\ \ \ \  (\om_{ij}+\om_{ji}=0).
\end{equation}

Differentiate our equation \eqref{newform}
$\displaystyle d\om_a+2\om_b\we\om_c=-{\bar {\bf J}}^a_{ij}\se^i\we\se^j$
\ ($(a,b,c)\sim (1,2,3)$) and substitute \eqref{reduce}.
Then it becomes (after alternation):
\[
(d{\bar {\bf J}}^a_{ij}-\om^{\sigma}_{i}{\bar {\bf J}}^a_{\sigma j}-
{\bar {\bf J}}^a_{i\sigma}\om^{\sigma}_j
+\om_b\cdot 2{\bar {\bf J}}^c_{ij}-\om_c\cdot 2{\bar {\bf J}}^b_{ij})
\we\se^i\we\se^j=0.
\]If we note \eqref{upper-equations} and \eqref{new-con},
$d{\bar {\bf J}}^a_{ij}-\om^{\sigma}_{i}{\bar {\bf J}}^a_{\sigma j}-
{\bar {\bf J}}^a_{i\sigma}\om^{\sigma}_j\equiv 0$ \mbox{mod}\ 
$\om_1,\om_2,\om_3$.
As the forms $\om_a\we\se^i\we\se^j$ are linearly independent,
the result follows.
\end{proof}

\begin{definition}\label{p-Kahler}
Let $\hat\nabla$ be
the Levi-Civita connection on an
 almost quaternionic pseudo-Riemannian manifold $(X,\hat g)$
 of type $(4p,4q)$ $(p+q=n)$.
 Then  $X$ is said to be a
quaternionic pseudo-K\"ahler manifold if
 for each quaternionic structure $\{\hat J_a; a=1,2,3\}$
 defined locally on a neighborhood of $X$, there exists a
 smooth local function $A\in{\mathfrak{so}}(3)$ such that
 \begin{equation*}
\hat\nabla\left(\begin{array}{c}
{\hat J}_1\\
{\hat J}_2\\
{\hat J}_3\end{array}\right)=
A\cdot
\left(\begin{array}{c}
{\hat J}_1\\
{\hat J}_2\\
{\hat J}_3\end{array}\right)
\end{equation*}provided that ${\rm dim}\  X=4n\geq 8$. 
Equivalently 
if $\hat \Omega$
is the fundamental $4$-form globally defined on $X$,
then $\hat\nabla\hat \Omega=0$.

\end{definition}
 
We have shown the following result in \cite{AK}
when  ${\rm dim}\  U/\mathcal E=4n\geq 12$
by  Swann's  method.
\begin{theorem}\label{q-K-s}
The set $(U/\mathcal E,\hat g,
\{{\hat {I}}_i,{\hat {J}}_i,{\hat {K}}_i\}_{i\in\Lambda})$
is a quaternionic pseudo-K\"ahler manifold
of type $(4p,4q)$  provided that ${\rm dim}\  U/\mathcal E=4n\geq 8$.
Moreover, $(U/\mathcal E,\hat g)$ is an Einstein manifold of positive
scalar curvature $(4n\geq 4)$ such that
\begin{equation}\label{ein}
\hat R_{j\ell}=4(n+2)\hat g_{j\ell}.
\end{equation}
\end{theorem}

\begin{proof}
As we put $\se^i=\pi^*\hat \se ^i$, 
the equation \eqref{levi}
implies that
$\displaystyle d\se^i=\se^j\we\pi^*\hat\om^i_j$,
$\displaystyle \pi^*\hat\om_{ij}+\pi^*\hat\om_{ji}=0$.
Compared this with the equation \eqref{reduce} and by skew-symmetry,
it is easy to check that
\begin{equation}\label{pullback}
\pi^*\hat\om^i_j=\om^i_j.
\end{equation}
Put $\hat V=V_i$ and
${\hat J}_1={\hat I}_{i},\
{\hat J}_2={\hat J}_{i},\ {\hat J}_3={\hat K}_{i}$
on $\hat V$.
Let $s=s_i: \hat V \ra U$ be the section as before.
Since $\pi_*s_*((\hat e_j)_x)=(\hat e_j)_x=\pi_*((e_j)_{s(\hat x)})$,
$s_*((\hat e_j)_x)-(e_j)_{s(\hat x)}\in V=\{\xi_1,\xi_2,\xi_3\}$.
Then
$\se^i(s_*((\hat e_j)_x))=\se^i((e_j)_{s(\hat x)})$
from \eqref{coframing}. 
A calculation shows that
$(\hat J_a)_{\hat x}\hat e_i=
\pi_*(J_a)_{s(\hat x)}e_i=\pi_*(({\bar {\bf J}}^a)_{i}^{j}(s(\hat x))e_j)
=({\bar {\bf J}}^a)_{i}^{j}(s(\hat x))\hat e_j$\ (\cf \eqref{3complex}).
As we put ${\hat J}^a_{\hat x}\hat e_i=
({\hat {\bf J}}^a)_{i}^{j}(\hat x)\hat e_j$, note that
\begin{equation}\label{Ipro}
{\bar {\bf J}}^a_{ij}(s(\hat x))={\hat {\bf J}}^a_{ij}(\hat x)\ \
(a=1,2,3).
\end{equation}
In particular,
\begin{equation}\label{IJK-pro}
d({\bar {\bf J}}^a)_{ij}\circ s_*(\hat X_{\hat x})=
d({\hat {\bf J}}^a)_{ij}(\hat X_{\hat x})\ \
(\forall\ \hat X_{\hat x}\in T_{\hat x}(\hat V))\ (a=1,2,3).
\end{equation}
Since $\pi_*s_*(\hat X_{\hat x})=
\hat X_{\hat x}$ $(\hat X_{\hat x}\in T_{\hat x}(\hat V))$,

\begin{equation}\label{induce-cone}
\hat\om^{\sigma}_j(\hat X_{\hat x})=\om^{\sigma}_j(s_*(\hat X_{\hat x})).
\end{equation}

Plug these equations \eqref{induce-cone}, \eqref{Ipro}
and \eqref{IJK-pro}
into \eqref{3-structure eqs}:
\begin{equation}\label{inter}
\begin{split}
&d({\bar {\bf J}}^a)_{ij}(s_*\hat X)
-\om^{\sigma}_{i}(s_*\hat X)\cdot ({\bar {\bf J}}^a)_{\sigma j}(s(\hat x))-
({\bar {\bf J}}^a)_{i\sigma}(s(\hat x))\cdot\om^{\sigma}_j(s_*\hat X)\\
&=d(({\hat {\bf J}}^a)_{ij})_{\hat x}(\hat X)
-{\hat\om}^{\sigma}_i(\hat X)\cdot({\hat {\bf J}}^a)_{\sigma j}(\hat x)
-({\hat {\bf J}}^a)_{i\sigma}(\hat x)\cdot {\hat\om}^{\sigma}_j(\hat X)\\
&=2({\bar {\bf J}}^b)_{ij}(s(\hat x))\cdot\om_c(s_*\hat X)-
2({\bar {\bf J}}^c)_{ij}(s(\hat x))\cdot\om_b(s_*\hat X)\\
&=2({\hat {\bf J}}^b)_{ij}(\hat x)\cdot\om_c(s_*\hat X)-
2({\hat {\bf J}}^c)_{ij}(\hat x)\cdot\om_b(s_*\hat X).
\end{split}
\end{equation}
Using these,
\begin{equation}
\begin{split}
&(\hat\nabla_{\hat X}({\hat J}_a)((\hat e_i)_{\hat x})=
\hat\nabla_{\hat X}({\hat J}_a)\hat e_i
-({\hat J}_a)(\hat\nabla_{\hat X}\hat e_i)\\
&=(d({\hat {\bf J}}^a)_{ij}(\hat X)
-({\hat {\bf J}}^a)_{i\sigma}(\hat x)\cdot\hat \om^{\sigma}_j(\hat X)
-\hat \om^{\sigma}_i(\hat X)\cdot({\hat {\bf J}}^a)_{\sigma j}(\hat x))
(\hat e_j)_{\hat x}\\
&=2({\hat {\bf J}}^b)_{ij}(\hat x)(\hat e_j)_{\hat x}
\cdot s^*\om_c(\hat X)-
2({\hat {\bf J}}^c)_{ij}(\hat x)(\hat e_j)_{\hat x}\cdot
 s^*\om_b(\hat X)\\
&=\Bigl(2({\hat J}_b)_{\hat x}\cdot s^*\om_c(\hat X)-
2({\hat J}_c)_{\hat x}\cdot 
s^*\om_b(\hat X)\Bigr)(\hat e_i)_{\hat x}.
\end{split}
\end{equation}Therefore,
$\hat\nabla_{\hat X}({\hat J}_a)=
2({\hat J}_b)_{\hat x}\cdot s^*\om_c(\hat X)-
2({\hat J}_c)_{\hat x}\cdot s^*\om_b(\hat X)$.
This concludes that
\begin{equation}\label{quat}
\hat\nabla\left(\begin{array}{c}
{\hat J}_1\\
{\hat J}_2\\
{\hat J}_3\end{array}\right)=
2\left(\begin{array}{lcr}
0&s^*\om_3&-s^*\om_2\\
-s^*\om_3&0&s^*\om_1\\
s^*\om_2&-s^*\om_1&0
\end{array}\right)
\left(\begin{array}{c}
{\hat J}_1\\
{\hat J}_2\\
{\hat J}_3\end{array}\right).
\end{equation}
As we put ${\hat J}_1={\hat I}_{i},\
{\hat J}_2={\hat J}_{i},\ {\hat J}_3={\hat K}_{i}$
on $\hat V$, $(U/\mathcal E,\hat g,
\{{\hat {I}}_i,{\hat {J}}_i,{\hat {K}}_i\}_{i\in\Lambda})$
is a quaternionic pseudo-K\"ahler manifold
for ${\rm dim}\ U/\mathcal E\geq 8$.
Using the Ricci identity $($cf. $(2.11), (2.12)$ of \cite{IS}, \cite{TA1}$)$,
a calculation shows that
\begin{equation}\label{3ricci}(n>1)\ \ \
\begin{split}
\hat R_{jl}&=-4(n+2)\Bigl(s^*(d\om_1+2\om_2\we\om_3)\Bigr)(\hat e_j,\hat e_k)
\hat {I}_{\ell}^{k}(\hat x).\\
\hat R_{jl}&=-4(n+2)\Bigl(s^*(d\om_2+2\om_3\we\om_1)\Bigr)(\hat e_j,\hat e_k)
\hat {J}_{\ell}^{k}(\hat x).\\
\hat R_{jl}&=-4(n+2)\Bigl(s^*(d\om_3+2\om_1\we\om_2)\Bigr)(\hat e_j,\hat e_k)
\hat {K}_{\ell}^{k}(\hat x).
\end{split}
\end{equation}
\begin{equation}\label{31ricci}(n=1)\ \ \ 
\begin{split}
\hat R_{jl}&=-4\Bigl(s^*(d\om_1+2\om_2\we\om_3)\Bigl)(\hat e_j,\hat e_k)
\hat {I}_{\ell}^{k}(\hat x)+\\
&\ \ \ \ \ \ \ \ \  -4\Bigl(s^*(d\om_2+2\om_3\we\om_1)\Bigr)
(\hat e_j,\hat e_k)\hat {J}_{\ell}^{k}(\hat x)+\\ 
&\ \ \ \ \ \ \ \ \ \ \ -4\Bigl(s^*(d\om_3+2\om_1\we\om_2)\Bigr)
(\hat e_j,\hat e_k)\hat {K}_{\ell}^{k}(\hat x).
\end{split}
\end{equation}
Using $d\omega_a+2\omega_b\wedge\omega_c=-{\bf J}^a_{ij}\theta^i\we\theta^j$
and \eqref{Ipro},
it follows that
$\Bigr(s^*(d\om_a+2\om_b\we\om_c)\Bigl)(\hat e_j,\hat e_k)=
-{\bf J}^a_{jk}(s(\hat x))=-\hat {\bf J}^a_{jk}(\hat x)$.
Since $(\hat {\bf J}^a)_i^j\cdot (\hat {\bf J}^a)_j^k=-\delta_{i}^{k}$,
$\hat R_{jl}=+4(n+2)(\hat {\bf J}^a)_{jk}(\hat x)\cdot
 (\hat {\bf J}^a)_{\ell}^{k}(\hat x)
=4(n+2)g_{j\ell}$ when $n>1$ and
$\hat R_{jl}=+4(\hat {I}_{jk}(\hat x)\cdot 
{\hat {I}_{\ell}}^{k}(\hat x)
+\hat {J}_{jk}(\hat x)\cdot {\hat {J}_{\ell}}^{k}(\hat x)
+\hat {K}_{jk}(\hat x)\cdot {\hat {K}_{\ell}}^{k}(\hat x))=
4\cdot 3g_{j\ell}$ when $n=1$.
\end{proof}

\section{Quaternionic $CR$ curvature tensor}\label{EP}

Recall from \eqref{reduce}
that $\displaystyle
d\se^i=\se^j\we\om^i_j,\ \  \om_{ij}+\om_{ji}=0$
where
$\pi^*\hat\om^i_j=\om^i_j$, $\pi^*\hat\se^i=\se^i$ from
\eqref{sasaki structure eq}, \eqref{pull} respectively 
 $(i,j=1,\cdots,4n)$.
Define the fourth-order tensor
$R^i_{jk\ell}$ on $U$ by putting
\begin{equation}\label{M-curvature}
d\om^i_j-\om_j^\sigma\we\om^i_{\sigma}\equiv\frac12
R^i_{jk\ell}\se^k\we\se^\ell\ \ \mbox{mod}\ \om_1,\om_2,\om_3.
\end{equation}
By \eqref{Rcur}, it follows that
\begin{equation}\label{updown-curvature}
R^i_{jk\ell}=\pi^*\hat R^i_{jk\ell}.
\end{equation}The equality \eqref{ein} implies that
\begin{equation}\label{Ein2}
R_{j\ell}=4(n+2)g_{j\ell}.
\end{equation}

Differentiate the structure equation \eqref{sasaki structure eq}.
\begin{equation}\label{diff}
0=d\se^j\we \f^i_j-\se^j\we d\f^i_j+\mathop{\sum}_a^{} d\om_a\we \tau^i_a-
\mathop{\sum}_a^{}\om_a\we d\tau^i_a.
\end{equation}Substitute \eqref{newform} and  \eqref{sasaki structure eq}
into \eqref{diff};
\begin{equation*}\label{1diff}
\begin{split}
&\se^j\we (d\f^i_j-\f^k_j\we\f^i_k -
\mathop{\sum}_a {\bf J}^a_{kj}\se^k\we\tau^i_a)
+\mathop{\sum}_a\om_a\we (d\tau^i_a-\tau^k_a\we\f^i_k)\\
&\ \ +2\om_2\we\om_3\we\tau^i_1+2\om_3\we\om_1\we\tau^i_2+
2\om_1\we\om_2\we\tau^i_3=0.
\end{split}
\end{equation*}
This implies  that
\begin{equation}\label{1eq}
\se^j\we(d\f^i_j-\f^k_j\we\f^i_k - \mathop{\sum}_a {\bf J}^a_{kj}
\se^k\we\tau^i_a)
\equiv 0\ \mbox{mod}\ \ \om_1,\om_2,\om_3.
\end{equation}
We use \eqref{1eq} to define the curvature form:
\begin{equation}\label{cur}
\Phi^i_j=d\f^i_j-\f^k_j\we\f^i_k+\mathop{\sum}_{a=1}^{3}
\se^k\we {\bf J}_{jk}^a\tau^i_a -\se^i\we \se_j.
\end{equation}
Set
\begin{equation}\label{bian}
\begin{split}
{}_1\Phi^i&=d\tau^i_1-\tau^k_1\we\f^i_k+\om_2\we\tau^i_3-\om_3\we\tau_2^i,\\
{}_2\Phi^i&=d\tau^i_2-\tau^k_2\we\f^i_k+\om_3\we\tau^i_1-\om_1\we\tau_3^i,\\
{}_3\Phi^i&=d\tau^i_3-\tau^k_3\we\f^i_k+\om_1\we\tau^i_2-\om_2\we\tau_1^i
\end{split}
\end{equation}
which satisfy the following relation.
\begin{equation}\label{1bian}
\se^j\we\Phi^i_j+\om_1\we{}_1\Phi^i
+\om_2\we{}_2\Phi^i+\om_3\we{}_3\Phi^i=0.
\end{equation}
We may define the fourth-order curvature tensor $T^i_{jkl}$ from $\Phi^i_j$:
\begin{equation}\label{ps-cur}
\Phi^i_j\equiv\frac 12T^i_{jkl}\se^k\we 
\se^{\ell}\ \ \mbox{mod}\ \om_1,\om_2, \om_3.
\end{equation}

\begin{remark} {\em In view of \eqref{ps-cur},
there exist the fourth-order curvature tensors
$W^i_{jka}$ $(a=1,2,3)$ and $V^i_{jbc}$ $(1\leq b<c\leq 3)$
for which we can describe: 
\begin{equation}\label{ps-curv}
\Phi^i_j=\frac 12T^i_{jkl}\se^k\we\se^{\ell}+\frac 12\mathop{\sum}_a^{}
W^i_{jka}\se^k\we\om_a+\frac 12\mathop{\sum}_{b<c}
V^i_{jbc}\om_b\we\om_c. 
\end{equation}
}
\end{remark}

\section{Transformation of pseudo-conformal $QCR$ structure}
\label{transform1}

\subsection{$G$-structure}\label{g-struc}
When $\{\se^i\}_{i=1,\cdots,4n}$ are the $1$-forms locally defined
on a neighborhood $U$ of $M$, we form the $\HH$-valued $1$-form 
$\{\om^i\}_{i=1,\cdots,n}$
such as
\begin{equation}\label{q-form}
\om^i=\se^i+\se^{n+i}\mbox{\boldmath$i$}+\se^{2n+i}\mbox{\boldmath$j$}+
\se^{3n+i}\mbox{\boldmath$k$}.
\end{equation}
We shall consider the transformations $f:U\ra U$
 of the following form:
\begin{equation}\label{pq-tran}
\begin{split}
f^*\om&=\lam\cdot\om\cdot\bar \lam\ (=u^2 a\cdot\om\cdot\bar a),\\
f^*(\om^j)&={U'}^j_{\ell}\om^{\ell}\cdot\bar \lam
+\lam\tilde v^j\om\bar\lam
\end{split}
\end{equation}such that $\lam=u\cdot a$ for 
some smooth functions $u>0$,
$a\in{\rm Sp}(1)$ and $U'\in {\rm Sp}(p,q)$ with
$p+q=n$.
Let $G$ be the subgroup of
$\mathop{\GL}(n+1,\HH)\cdot \HH^*$ consisting of matrices

\begin{equation}\label{Gmat}
\left(\begin{array}{c|ccc}
\lam &   &\mbox{\Large $0$}& \\
\hline
     &  &  & \\
\lam\cdot{\tilde v}^i &  & U' & \\
   & & &
 \end{array}\right)\cdot \lam.
\end{equation}
Recall that $\mathop{\Sim}(\HH^n)=
\HH^n\rtimes ({\rm Sp}(p,q)\cdot\HH^*)$ 
is the quaternionic affine similarity group where
$\HH^*={\rm Sp}(1)\times \RR^+$.
Then note that $G$ is anti-isomorphic to 
$\mathop{\Sim}(\HH^n)$ given by the map 
\begin{equation}
\label{similarity-tr}
{{}^{}}^t\left(\begin{array}{c|ccc}
\lambda &  & x^j & \\
\hline
                &  & & \\
\mbox{\Large $0$}& & X & \\
                 & &   &
  \end{array}\right)\cdot \lambda
\lra (X{x^j}^*,X\cdot \lambda)\in \HH^n\rtimes
({\rm Sp}(p,q)\cdot \HH^*).
\end{equation}(Here $x^*={}^t\bar x$.)
We represent $G$ as the real matrices.
Let ${\tilde v}$ be a vector of the quaternionic vector space
$\HH^n$.
The group ${\rm Sp}(p,q)\cdot {\HH}^*$ is the subgroup
of $\mathop{\GL}(4n,\RR)$ acting on $\HH^n$ by
\begin{equation}
(U'\cdot \lam)\tilde v=U'\tilde v\cdot \bar \lam
\end{equation}where $U'\in{\rm Sp}(p,q)$, $\lam\in{\HH}^*$.
Write $\lambda=u\cdot a\in \RR^+\times{\rm Sp}(1)$
so that ${\rm Sp}(p,q)\cdot {\HH}^*$
is embedded into 
$\RR^+\times {\rm SO}(4p,4q)$ in the following manner:
\begin{equation}\label{actionAD}
\begin{split}
U'\cdot \lam(\tilde v)=uU'\tilde v\bar a=
uU'\bar a\circ (a \tilde v \bar a)=u(U'\bar a)\circ {\rm Ad}_{a}(\tilde v)
=u\cdot U\tilde v \quad (\tilde v\in \HH^n=\RR^{4n})
\end{split}
\end{equation}in which  
\begin{equation}\label{SPSO}
U=U'\bar a\circ {\rm Ad}_{a}\in {\rm SO}(4p,4q),
\end{equation}

\begin{equation}\label{barA}
{\rm Ad}_{a}\left(\begin{array}{c}
\mbox{\boldmath$i$}\\
\mbox{\boldmath$j$}\\
\mbox{\boldmath$k$}\end{array}\right)
=a\left(\begin{array}{c}
\mbox{\boldmath$i$}\\
\mbox{\boldmath$j$}\\
\mbox{\boldmath$k$}\end{array}\right)\bar a
=A\left(\begin{array}{c}
\mbox{\boldmath$i$}\\
\mbox{\boldmath$j$}\\
\mbox{\boldmath$k$}\end{array}\right)
\ \mbox{for some}\  
A\in {\rm SO}(3).
\end{equation}

We put the vector ${\tilde v}^j\in \HH^n$ in
 such a way that
$\displaystyle {\tilde v}^j=v^j+v^{n+j}\mbox{\boldmath$i$}+
v^{2n+j}\mbox{\boldmath$j$}+v^{3n+j}\mbox{\boldmath$k$}$
$(j=1,\cdots,n)$.
Form the real $(4\times 3)$-matrix
\begin{equation}
V^j=\left(\begin{array}{ccc}
 -v^{j+n}   & -v^{j+2n} &-v^{j+3n}\\
  v^j       & -v^{j+3n} &v^{j+2n}\\
v^{j+3n} & v^{j}   &-v^{j+n}\\  
-v^{j+2n} & v^{j+n}   &v^{j}
  \end{array}\right).
\end{equation}It is easy to check that
\begin{equation}
\begin{split}
\lam{\tilde v}^j\left(\begin{array}{c}
 \om_1\\
 \om_2\\
  \om_3\end{array}\right)\bar\lam=
\lam((1\ \mbox{\boldmath$i$}\ \mbox{\boldmath$j$}\ \mbox{\boldmath$k$})V^j
\left(\begin{array}{c}
 \om_1\\
 \om_2\\
  \om_3\end{array}\right))\bar \lam 
=(1\ \mbox{\boldmath$i$}\ \mbox{\boldmath$j$}\ \mbox{\boldmath$k$})
u^2 \left(\begin{array}{cc}
 1   &0\\
 0   &{}^tA
  \end{array}\right)V^j
\left(\begin{array}{c}
 \om_1\\
 \om_2\\
  \om_3\end{array}\right).
\end{split}\end{equation}

Then $G$ is isomorphic to the subgroup of
$\mathop{\GL}(4n+3,\RR)$ consisting of matrices
\begin{equation}\label{G}
\left(\begin{array}{c|ccc}
                &  & &\\
u^2\cdot {}^tA  &  &\mbox{\Large $0$} & \\
                &  & & \\
\hline
                &      & &           \\
u^2 \left(\begin{array}{cc}
 1   &0\\
 0   &{}^tA
  \end{array}\right)V^1  &  & & \\
 \vdots            &   & \mbox{\Large $u\cdot U$} &\\
u^2 \left(\begin{array}{cc}
 1   &0\\
 0   &{}^tA
  \end{array}\right)V^n    & & & 
\end{array}\right).
\end{equation}
Here $A\in {\rm SO}(3), U=(U^i_j)\in{\rm SO}(4p,4q)$.\\
Using the coframe field $\{\om_1,\om_2,\om_3,\se^1,\cdots,\se^{4n}\}$, $f$
is represented by 
\begin{equation}\label{1se}
\begin{split}
&f^*(\om_1,\om_2,\om_3)=u^2(\om_1,\om_2,\om_3)A,\\
&f^*\se^i=u\se^kU^i_k+\mathop{\sum}_{\al=1}^{3}\om_\al {v}^i_\al,\\
&\mbox{where}\ \left(\begin{array}{lcr}
{v}^{4j-3}_1& {v}^{4j-3}_2& {v}^{4j-3}_3\\
{v}^{4j-2}_1& {v}^{4j-2}_2& {v}^{4j-2}_3\\
{v}^{4j-1}_1& {v}^{4j-1}_2& {v}^{4j-1}_3\\
{v}^{4j}_1& {v}^{4j}_2& {v}^{4j}_3
\end{array}\right)
=u^2 \left(\begin{array}{cc}
 1   &0\\
 0   &{}^tA
  \end{array}\right)V^j\ \ (j=1,\cdots,n).
\end{split}
\end{equation}
Let $\mathcal F(M)$ be the principal coframe bundle over $M$.
A subbundle $P$ of
$\mathcal F(M)$ is said
to be a {\em bundle of the nondegenerate integrable $G$-structure}
if $P$ is the total space of the
principal bundle $G\ra P\ra M$ whose points 
consist of such coframe
fields $\{\om_1,\om_2,\om_3,\se^1,\cdots,\se^{4n}\}$
satisfying the condition \eqref{qcc}, \eqref{integraleq1}, \eqref{localflow}. 
A diffeomorphism $f:M\ra M$ is
a $G$-automorphism 
if the derivative $f^*: \mathcal F(M)\ra \mathcal F(M)$
induces a bundle map $f^*:P\ra P$ in which 
$f^*$ has the form locally as in \eqref{pq-tran}
(equivalently \eqref{1se}).
\begin{definition}\label{auto-pqcr}
Let $\mathop{\Aut}_{QCR}(M)$ be the group of 
all $G$-automorphisms of $M$.
\end{definition}

\subsection{Automorphism group  $\mathop{\Aut}(M)$}\label{review}
Let $W$ be the $(n+2)$-dimensional arithmetic vector space $\HH^{p+1,q+1}$
 over $\mathbb H$ equipped with the standard
Hermitian metric $\mathcal B$ of signature $(p+1,q+1)$ where $p+q=n$.
Then note that the isometry group
$\mathcal G =\mathop{\Aut}(W,\mathcal B)={\rm Sp}(p+1,q+1)$ 
and $W$ has the gradation 
$\displaystyle W = W^{-1} + W^{0} + W^{+1}$, 
where $W^{\pm 1}$ are dual 1-dimensional
isotropic subspaces and $W^0$ is ($h$-non-degenerate ) orthogonal
 complement to $W^{-1} + W^{+1}$.
The gradation $W$ induces 
the gradation of the Lie algebra $\mathfrak g$ of
depth two, \ie
\[\mathfrak g = \mathfrak g^{-2}+\mathfrak g^{-1}+
             \mathfrak g^0 + \mathfrak g^1 + \mathfrak g^2.
\]Here $\mathfrak g^0 = \mathbb R + \mathfrak {sp}(1) + \mathfrak{sp}(n)$.

In \cite{AK1} we introduced a notion of
pseudo-conformal quaternionic structure.
This geometry is defined  by a 
codimension three distribution $\mathcal H$ on a $(4n +3)$-dimensional
manifold $M$, 
which satisfies the only one condition
that the associated graded tangent space 
$\displaystyle {}^{gr}T_xM=T_xM/\mathcal H_x+\mathcal H_x$
at any point is isomorphic to the quaternionic Heisenberg Lie algebra 
$\mathfrak M(p,q)\cong \mathfrak g^{-}=\mathfrak g^{-2}+
\mathfrak g^{-1}$, \ie the Iwasawa subalgebra of 
${\rm Sp}(p+1,q+1)$. Put $\displaystyle {\rm Sp}(W)={\rm Sp}(p+1,q+1)$.
We proved that such a geometry is a parabolic geometry so that
it admits a canonical Cartan connection and 
its automorphism group $\mathop{\Aut}(M)$ is a Lie group.
 More precisely,
if $P^{+}(\HH)$ is the parabolic connected subgroup of the symplectic group
${\rm Sp}(W)$ corresponding to the
dual parabolic subalgebra $\fp^{+}(\HH) = \fg^{+} + \fg^0$
of $\fs \fp (W)$, then there is a $P^+(\HH)$-principal bundle 
$\pi : B \rightarrow M$ with a normal Cartan connection 
$\kappa : TB \rightarrow \fs \fp(W) $ of type ${\rm Sp}(W)/P^+(\HH)$.
There exists the canonical pseudo-conformal quaternionic structure
$\mathcal H^{\rm can}$ on ${\rm  Sp}(p+1, q+1)/P^+(\HH)$ with all vanishing
curvature tensors (\cf $\S \ref{D-H}$). 
A pseudo-conformal quaternionic manifold $(M, \mathcal H)$ 
is locally isomorphic to a 
$({\rm  Sp}(p+1, q+1)/P^+(\HH), \mathcal H^{\rm can})$ if and only if 
the associated Cartan connection $\kappa$ is flat (\ie has zero curvature).
Put $\displaystyle S^{4p+3,4q}={\rm  Sp}(p+1,q+1)/P^+(\HH)$.
Then $S^{4p+3,4q}$ is the flat homogeneous model diffeomorphic to 
$\displaystyle S^{4p+3}\times S^{4q+3}/{\rm Sp}(1)$
where the product of spheres $\displaystyle S^{4p+3}\times S^{4q+3}
=\{(z^+,z^-)\in \HH^{p+1,q+1}\ |\ \mathcal B(z^+,z^+)=1,
 \mathcal B(z^-,z^-)=-1\}$ is the subspace of $W=\HH^{p+1,q+1}$
and the action of ${\rm Sp}(1)$ is induced by the diagonal right 
action on $W$.
The group of all automorphisms $\displaystyle\mathop{\Aut}(S^{4p+3,4q})$
preserving this flat structure is ${\rm PSp}(p+1, q+1)$.

Suppose that $M$ is a pseudo-conformal quaternionic $CR$ manifold.
By definition,
$T_xM\cong T_xM/\mathcal D_x+\mathcal D_x
={\rm Im}\HH+\HH^n\cong {\mathfrak M}(p,q)$ at $\forall\ x\in M$.
Then each $G$-automorphism of 
$\mathop{\Aut}_{QCR}(M)$ 
preserves ${\mathfrak M}(p,q)$ by the above formula \eqref{1se}.
Since a pseudo-conformal quaternionic $CR$ structure
is a pseudo-conformal quaternionic structure
by Definition \ref{pc-qcr}, note that
$\mathop{\Aut}_{QCR}(M)$ is a (closed) subgroup of
$\mathop{\Aut}(M)$ which is a Lie group as above.

\begin{corollary}\label{Lie-trans}
The group $\mathop{\Aut}_{QCR}(M)$ is a finite dimensional Lie group
for a pseudo-conformal quaternionic $CR$ manifold $M$.
\end{corollary}

\section{Pseudo-conformal $QCR$ structure on $S^{3+4p,4q}$}
\label{existPQCR} 
We shall prove that the quaternionic $CR$  homogeneous model
$\Sigma^{3+4p,4q}_{\HH}$ induces a pseudo-conformal quaternionic $CR$ structure on $S^{3+4p,4q}$ which coincides with
the flat pseudo-conformal quaternionic structure.

\subsection{Quaternionic pseudo-hyperbolic geometry}\label{review2}
Let 
\[
\mathcal B(z,w) = {\bar z}\sb 1w\sb 1+ {\bar z}\sb 2w \sb 2+ \cdots +
 {\bar z}\sb{p+1} w\sb{p+1}-{\bar z}\sb{p+2} w\sb{p+2}-\cdots-
 {\bar z}\sb{n+2} w\sb{n+2}
\]be the above Hermitian form on $\HH^{n+2}=\HH^{p+1,q+1}$ $(p+q=n)$.
We consider the following subspaces in $\HH^{n+2}-\{0\}$:
\begin{equation*}
\begin{split}
V \sb 0 \sp {4n+7} &=\{z \in \HH^{n+2} \vert \ \mathcal B(z,z)=0\}, \\
V \sb-\sp{4n+8}  &= \{z \in \HH^{n+2} \vert \ \mathcal B(z,z)<0\}.
\end{split}
\end{equation*}

Let $\displaystyle \HH^*\ra (({\rm Sp}(p+1,q+1)\cdot \HH^*,\HH^{n+2}-\{0\})
\stackrel P{\lra} ({\rm PSp}(p+1,q+1), \HH\PP^{n+1})$
be the equivariant projection.
The quaternionic pseudo-hyperbolic space $\HH \sb{\HH}\sp{p+1,q}$ 
is defined to be $P(V\sb{-}\sp {4n+8})$ (\cf \cite{CG}).
Let $\mathop{\GL}(n+2,\HH)$ be the group of all invertible
$(n+2)\times (n+2)$-matrices with
quaternion entries. Denote by ${\rm Sp}(p+1,q+1)$ the subgroup
consisting of
\begin{equation}
\begin{split}
\{A\in \mathop{\GL}(n+2, \HH)\ |\ 
\mathcal B(Az,Aw)=\mathcal B(z,w), z,w\in \HH^{n+2}\}.
\end{split}
\end{equation}
The action ${\rm Sp}(p+1,q+1)$ on  $V\sb{-}\sp {4n+8}$ induces an action on
$\HH\ \sb{\HH}\sp {p+1,q}$.
The kernel of this action is the center $\ZZ/2=\{\pm 1\}$
whose quotient is
the pseudo-quaternionic hyperbolic group
$\mathop{\PSp}(p+1,q+1)$.
It is known that $\HH \sb{\HH}\sp{p+1,q}$ is 
a complete simply connected pseudo-Riemannian manifold
of negative sectional curvature from $-1$ to $-\frac14$,
and with the group of isometries $\mathop{\PSp}(p+1,q+1)$
(\cf \cite{KONO}).
Remark that when $q=0, p=n$, $P(V\sb{-}\sp {4n+8})=\HH \sb{\HH}\sp{n+1}$
is the quaternionic K\"ahler hyperbolic space with the group of 
isometries $\mathop{\PSp}(n+1,1)$.
The projective compactification of
$\HH_{\HH}^{p+1,q}$ is
obtained by taking the closure
$\bar{\HH}_{\HH}^{p+1,q}$ in $\HH\PP^{n+1}$.
Then it is easy to check that
$\bar{\HH}_{\HH}\sp{p+1,q}=
\HH_{\HH}\sp{p+1,q}\cup P(V \sb 0 \sp {4n+7})$.
The boundary $P(V \sb 0 \sp {4n+7})$ of $\HH \sb{\HH}\sp{p+1,q}$ is 
identified with the quadric $S^{3+4p,4q}$
by the correspondence:
\begin{equation}
\begin{split}
\left[z_+,z_-\right]
\mapsto
\left[\frac{z_+}{|z_-|},\frac{z_-}{|z_-|}\right].
\end{split}
\end{equation}
Since the pseudo-hyperbolic action of
$\mathop{\PSp}(p+1,q+1)$ on $\HH^{p+1,q}_{\HH}$ extends to
a smooth action on $S^{3+4p,4q}=P(V \sb 0 \sp {4n+7})$
as projective transformations because
the projective compactification $\bar \HH^{p+1,q}_{\HH}$ is
an invariant domain of $\HH\PP^{n+1}$.

\subsection{Existence of pseudo-conformal $QCR$ structure on $S^{3+4p,4q}$}
\label{D-H} 
Recall that
$\displaystyle\Sigma^{3+4p,4q}_{\HH}
=\{(z_1,\cdots,z_{p+1},w_1,\cdots,w_q)\in {\HH}^{n+1}\ |\
|z_1|^2+\cdots+|z_{p+1}|^2-|w_1|^2-\cdots-|w_q|^2=1\}$ equipped
with the quaternionic $CR$ structure $\om_0$ (\cf $\S \ref{MCR}$).
The embedding $\iota$ of $\Sigma^{3+4p,4q}_{\HH}$
into $S^{4p+3,4q}$ is defined by
$\displaystyle (z_1,\cdots,z_{p+1},w_1,\cdots,w_q)\mapsto 
\left[(z_1,\cdots,z_{p+1},w_1,\cdots,w_q,1)\right]$.
Then $\iota(\Sigma^{3+4p,4q}_{\HH})$ is an open dense submanifold 
of $S^{4p+3,4q}$ since it misses the subspace
$\displaystyle S^{4p+3,4(q-1)}=S^{4p+3}\times S^{4q-1}/
{\rm  Sp}(1)$ in $S^{4p+3,4q}$.
We know that 
$\Sigma^{3+4p,4q}_{\HH}$ has the transitive isometry
group ${\rm Sp}(p+1,q)\cdot {\rm Sp}(1)$ (\cf Definition \ref{fqcc}).
Then this embedding implies that
${\rm Sp}(p+1,q)\cdot {\rm Sp}(1)$ is identified with the subgroup
$P({\rm Sp}(p+1,q)\times {\rm Sp}(1))$ of
${\rm PSp}(p+1,q+1)$ leaving the last component $z_{n+2}$ invariant
in $V\sb 0\sp {4n+7}\subset\HH^{n+2}$.

By pullback, each element $h$ of ${\rm PSp}(p+1,q+1)$ gives a
 quaternionic $CR$ structure $h^{-1*}\om_0$ on the open subset
$h(\Sigma^{3+4p,4q}_{\HH})$ of $S^{3+4p,4q}$.
Noting that $h^{-1*}\mathcal H^{can}=\mathcal H^{can}$
and Definition \ref{pc-qcr},
we shall prove that $(S^{3+4p,4q},\mathcal H^{can})$ admits a pseudo-conformal
quaternionic $CR$ structure by showing that
$\mathcal{\Null}h^{-1*}\om_0$
coincides with the restriction of
$\mathcal H^{can}|h(\Sigma^{3+4p,4q}_{\HH})$.

\begin{theorem}\label{pcq=pcqcr}
The $(4n+3)$-dimensional 
pseudo-conformal quaternionic manifold
$\displaystyle (S^{4p+3,4q},\mathcal H^{can})$
supports a pseudo-conformal quaternionic $CR$ structure, \ie
there exists locally a quaternionic $CR$ strucrure $\om$ on a neighborhood $U$ such that
\[
\mathcal H^{can}|U=\mathcal{\Null}\om.
\]Moreover, the automorphism group
$\mathop{\Aut}_{QCR}(S^{4p+3,4q})$ with respect to this
pseudo-conformal quaternionic $CR$ structure
is ${\rm PSp}(p+1,q+1)$. 
\end{theorem}

\begin{proof}
First we describe the canonical pseudo-conformal quaternionic
structure $\mathcal H^{can}$ on $S^{3+4p,4q}$ explicitly.
Choose isotropic vectors $x,y\in V_0$ 
such that $\mathcal B(x,y)=1$ and denote by $V$ 
the orthogonal complement to $\{x,y\}$ in $\HH^{p+1,q+1}$.
Then it follows that
$\displaystyle T_xV_0 = {\fs}{\fp}(W)x = y{\rm Im}\HH + V + x\HH$
where $T_x(x\HH^*)=x\HH$. Then 
\[
T_{[x]}S^{4k+3,4q}=
P_*(T_xV_0) = (y{\rm Im}\HH + V + x\HH )/x\HH.
\]
We associate to each $[x]\in S^{4k+3,4q}$ the orthogonal complement 
$x^{\perp}= V + x\HH$. It does not depend on the choice of
points from $[x]$. In fact, if $x'\in [x]$, then $x'=x\cdot \lam$ 
for some $\lam \in\HH^*$.  By the definition choosing $y'$
such that $\displaystyle T_{x'}V_0 = y'{\rm Im}\HH + V' + {x'}\HH$
where the orthogonal complement $V'$ to
$\{x',y'\}$ in $\HH^{p+1,q+1}$ is uniquely determined.
Let $v'$ be any vector of $V'$ which is described as
$v'=y\cdot a+v+x\cdot b$ for some $a,b\in \HH$.
Then 
\begin{equation*}\begin{split}
0=\mathcal B(x',v')&=\mathcal B(x',y)a+\mathcal B(x',v)+\mathcal B(x',x)b\\
&=\bar \lam\mathcal B(x,y)a+\bar\lam \mathcal B(x,v)+\bar\lam
\mathcal B(x,x)b=\bar\lam\cdot a.
\end{split}\end{equation*}Since $\lam\neq 0$, $a=0$ and so
$v'=v+x\cdot b$. Hence ${x'}^{\perp}=V'+x'\HH=V+x\HH$.
Therefore the orthogonal complement $x^{\perp} = V + x\HH$ in $\HH^{p+1,q+1}$
determines  a codimension three subbundle
\begin{equation}\label{pcqbdle}
\begin{split}
\mathcal H^{can}&=\mathop{\cup}_{[x]\in S^{4p+3,4q}}^{}P_*( x^{\perp}).\\
\displaystyle P_*( x^{\perp})&=V + x\HH/x\HH\subset TS^{4p+3,4q}.
\end{split}\end{equation}On the other hand, recall that if $N_p$ is the
 normal vector at $p\in \Sigma^{3+4p,4q}_{\HH}$,
then  $(\mathop{\Null}\om_0)_p=\mathcal D_p=\{IN_p,JN_p,KN_p\}^{\perp}$ by the definition
 (\cf $\S$ \ref{MCR}).  Since $T_p\Sigma^{3+4p,4q}_{\HH}=N_p^{\perp}$ with respect to
$g^{\HH}$, it follows that $T_p\HH^{n+1}|\Sigma^{3+4p,4q}_{\HH}=
\{N_p,IN_p,JN_p,KN_p\}\oplus
\mathcal D_p$. If we note that
$\{N_p,IN_p,JN_p,KN_p\}=p\HH$, then we have
$\displaystyle \mathcal D_p=p\HH^{\perp}$.
It is easy to see that the orthogonal complement to $p\HH$ with respect 
to $g^{\HH}$ coincides with the orthogonal complement to $p$ with respect to 
the inner product $\mathcal B=\langle, \rangle$.
Hence, $\mathcal D_p=p^{\perp}$.
As the tangent subspace $\iota_*(\mathcal D_p)$
at $\iota(p)$ in $T_{\iota(p)}V_0$ is $(\mathcal D_p,0)$ which is parallel to 
$\mathcal D_p$ in $T_pV_0$, it implies that
$\displaystyle \mathcal B(\iota_*(\mathcal D_p), \iota(p))=
\mathcal B((\mathcal D_p,0),(p,1))=\langle \mathcal D_p, p\rangle -\langle  0,1\rangle =0$.
Hence $\iota_*(\mathcal D_p)\subset {\iota(p)}^{\perp}$ (with respect to $\mathcal B$). As ${\iota(p)}^{\perp} = V + {\iota(p)}\HH$,
$\iota_*(\mathcal D_p)\subset V + {\iota(p)}\HH$.
As above $\iota_*(\mathcal D_p)=(\mathcal D_p,0)$ at $\iota(p)$, but
$\iota(p)\HH=(p,1)\cdot H$. The intersection
$\iota_*(\mathcal D_p)\cap \iota(p)\HH=\{0\}$.
It implies that
$\displaystyle \iota_*(\mathcal D_p)=
\iota_*(\mathcal D_p)/\iota(p)\HH\subset  V + {\iota(p)}\HH/{\iota(p)}\HH$.
By \eqref{pcqbdle},
$\displaystyle \iota_*((\mathop{\Null}\om_0)_p)= P_*( \iota(p)^{\perp})=
\mathcal H^{can}_{\iota (p)}$.
Therefore $S^{4p+3,4q}$ admits a pseudo-conformal quaternionic $CR$ structure.
Then $\mathop{\Aut}_{QCR}(S^{4p+3,4q})$
is a subgroup of $\mathop{\Aut}(S^{4p+3,4q})={\rm  PSp}(p+1, q+1)$ 
from $\S \ref{review}$. \end{proof}
To prove $\mathop{\Aut}_{QCR}(S^{4p+3,4q})={\rm PSp}(p+1, q+1)$,
 we recall the quaternionic Heisenberg Lie group.

\subsection{Pseudo-conformal quaternionic Heisenberg geometry}\label{pcHei}
Let ${\rm  PSp}(p+1, q+1)$ 
be the group of all automorphisms 
 preserving the flat pseudo-conformal quaternionic structure
of $S^{4p+3,4q}={\rm  PSp}(p+1, q+1)/P^+(\HH)$ (\cf \ref{review}.)
We consider the stabilizer of the point at infinity 
$\{\infty\}=[1,0,\cdots,0,1]\in \Sigma^{3+4p,4q}_{\HH}\subset S^{4p+3,4q}$.
Recall the (indefinite) Heisenberg nilpotent Lie group $\mathcal M=\mathcal M(p,q)$
from \cite{KA}.
It is the product $\RR^3\times \HH^n$ with group law:
\[(a,y)\cdot (b,z)=(a+b-\mbox{Im}\langle y,z\rangle, y+z).\]
Here  $\langle  \,\ \rangle$ is 
the Hermitian inner product of signature $(p,q)$ on $\HH^n$
and $\mbox{Im}\langle  \,\ \rangle$ is the imaginary part\ $(p+q=n)$.
It is nilpotent because the commutator subgroup
 $[\mathcal M,\mathcal M]=\RR^3$
which is the center consisting of the form $(a,0)$.
In particular, there is the central extension:
\begin{equation}\label{nil-center}
1\ra \RR^3\ra \mathcal M\lra \HH^n\ra 1.
\end{equation}
Denote by $\mathop{\Sim}({\mathcal M})$ the semidirect product
${\mathcal M}\rtimes ({\rm Sp}(p,q)\cdot {\rm Sp}(1)\times \RR\sp {+})$
where the action $(A\cdot g,t)\in {\rm Sp}(p,q)\cdot {\rm Sp}(1)
\times \RR\sp{+}$ on $(a,y)\in\mathcal M$
is given  by:
\begin{equation}
(A\cdot g,t)\circ (a,y) =(t^{2}\cdot gag^{-1},\ t\cdot Ayg^{-1}).
\end{equation}
Denote the origin by $O=[1,0,\cdots,0,-1]\in 
\Sigma^{3+4p,4q}_{\HH}-\{\infty\}$.
Then, the stabilizer $\mathop{\Aut}(S^{3+4p,4q})_{\infty}$
is isomorphic
to $\mathop{\Sim}({\mathcal M})$
(\cf \cite{KA2}).
The orbit $\mathcal M\cdot O$ is a dense open subset
of $S^{4p+3,4q}$.  The embedding $\iota$ is defined by:
\begin{equation}\label{embeddheisen}
((a,b,c),(z_+,z_-))\in\mathcal M
\stackrel{\iota}\ra \left[\begin{array}{c}
\frac {|z_+|^2-|z_-|^2}{2} -1+ {\bf i}a+{\bf j}b+{\bf k}c\\
\sqrt 2 z_+\\
\sqrt 2 z_-\\
\frac {|z_+|^2-|z_-|^2}{2} +1+ {\bf i}a+{\bf j}b+{\bf k}c
\end{array}\right]
\end{equation}
Then the pair 
$(\mathop{\Sim}({\mathcal M}),\mathcal M)$ is 
said to be {\em pseudo-conformal quaternionic Heisenberg geometry}
which is a subgeometry of
the flat pseudo-conformal quaternionic geometry\\
$(\mathop{\Aut}(S^{3+4p,4q}),S^{3+4p,4q})$.
We prove the rest of Theorem \ref{pcq=pcqcr}.
\begin{pro}\label{flat2}
$\mathop{{\Aut}}_{QCR}(S^{4p+3,4q})={\rm PSp}(p+1,q+1)$.
\end{pro}

\begin{proof}
First note that $\mathop{\PSp}(p+1,q+1)$ decomposes
into $\mathop{\Sim}({\mathcal M})\cdot 
({\rm Sp}(p+1,q)\cdot{\rm Sp}(1))$. We know from $\S \ref{MCR}$ that
each element $f=(A,a)\in{\rm Sp}(p+1,q)\cdot{\rm Sp}(1)$
satisfies that $f^*\om_0=a\om_0 \bar a$, obviously 
$f\in \mathop{\Aut}_{QCR}(S^{4p+3,4q})$.
On the other hand, it is shown that an
element $h$ of $\mathop{\Sim}({\mathcal M})$
satisfy that $h^*\om_0= \lam\om_0 \bar\lam$ for
 some function $\lam\in {\HH}^*$ by using the explicit formula of $\om_0$.
(See \cite{KA}.) 
When $h\in\mathop{\Sim}({\mathcal M})$, note that
$h(\infty)=\infty$.
Let
$\tau : \mathop{\PSp}(p+1,q+1)_{\infty}\ra {\Aut}(\mbox{T}_{\{\infty\}}
(S^{3+4p,4q}))$ be the tangential representation
at $\{\infty\}$.
Since the elements of
 the center $\RR^3$ of ${\mathcal M}$
are tangentially identity maps at 
$\mbox{T}_{\{\infty\}}(S^{3+4p,4q})$, 
$\tau(\mathop{\PSp}(p+1,q+1)_{\infty})=
{\HH^n}\rtimes ({\rm Sp}(p,q)\cdot {\rm Sp}(1)\times \RR\sp {+})$
which is isomorphic to the structure group $G$
(\cf \eqref{G}).
As $\tau(h)=h_*$,  $h\in\mathop{\Aut}_{QCR}(S^{3+4p,4q})$
by Definition \ref{auto-pqcr}.
We have ${\rm PSp}(p+1,q+1)\subset\mathop{\Aut}_{QCR}(S^{3+4p,4q})$.
\end{proof}

\section{Pseudo-conformal quaternionic $CR$ invariant}\label{EP1}
We shall consider the equivalence of pseudo-conformal quaternionic 
$CR$ structure.
Let $d\omega+\omega\wedge\omega=
-(I_{ij}\mbox{\boldmath$i$}+J_{ij}\mbox{\boldmath$j$}
+K_{ij}\mbox{\boldmath$k$})\theta^i\wedge\theta^j$
be the equation \eqref{3-qcr} as before.
We examine how this equation behaves under the transformation
$f\in \mathop{\Aut}_{QCR}(M)$;
 $f^*\om=\lam\cdot\om\cdot \bar\lam$. Put $\om'=f^*\om$.
By \eqref{1se},
\begin{equation*}
\begin{split}
d\om'+\om'\we\om'&=f^*(d\om+\om\we\om)=
-(I_{ij}\mbox{\boldmath$i$}+J_{ij}\mbox{\boldmath$j$}+
K_{ij}\mbox{\boldmath$k$})f^*\se^i\we f^*\se^j\\
&=-(I_{ij}\mbox{\boldmath$i$}+J_{ij}\mbox{\boldmath$j$}+
K_{ij}\mbox{\boldmath$k$})(u\se^k U^i_k+\mathop{\sum}_a\om_a v^i_a)\we
(u\se^\ell U^j_\ell+\mathop{\sum}_b\om_b v^j_b)\\
&=-(I_{ij}\mbox{\boldmath$i$}+J_{ij}\mbox{\boldmath$j$}+
K_{ij}\mbox{\boldmath$k$})\Bigl(u^2U^i_kU^j_\ell\se^k\we\se^\ell+\\
&\ \ \ 
\mathop{\sum}_a\om_a\we(uv_a^iU_{\ell}^j\se^{\ell}
-uv_a^jU^i_\ell\se^\ell)\Bigr.
\Bigl.+\mathop{\sum}_{a<b}{\om_a\we\om_b}(v_a^iv_b^j-v_b^iv_a^j)
\Bigr)\\
&=-(I_{ij}\mbox{\boldmath$i$}+J_{ij}\mbox{\boldmath$j$}+
K_{ij}\mbox{\boldmath$k$})\Bigl(u^2U^i_kU^j_\ell\se^k\we\se^\ell
+\mathop{\sum}_a\om_a\we 2uv_a^iU_{\ell}^j\se^{\ell}\Bigr.\\
&\ \ \ \ \Bigl.+\mathop{\sum}_{a<b}{\om_a\we\om_b}(2v_a^iv_b^j)\Bigr).
\end{split}
\end{equation*}
Choosing $w^i_a$ $(a=1,2,3)$ such that $U_k^iw_a^k=v_a^i$,
the above equation becomes

\begin{equation*}
\begin{split}
d{\om'}+{\om'}\we{\om'}&=
-(I_{ij}\mbox{\boldmath$i$}+J_{ij}\mbox{\boldmath$j$}+
K_{ij}\mbox{\boldmath$k$})\Bigl(u^2U^i_kU^j_\ell\se^k\we\se^\ell+\Bigr.\\
&\quad \Bigl.\mathop{\sum}_a\om_a\we 2uw_a^kU_k^iU_{\ell}^j\se^{\ell}
+\mathop{\sum}_{a<b}{\om_a\we\om_b}(2U^i_kU^j_\ell w_a^kw_b^\ell)
\Bigr).
\end{split}
\end{equation*}
Let $U=U'\bar a\circ {\rm Ad}_{a}\in {\rm SO}(4p,4q)$
be the matrix as in \eqref{SPSO} so that
$Uz=U'z\bar a$ $(z\in \HH^n)$ (\cf \eqref{actionAD}).
If $\{I,J,K\}$ is the set of
the standard quaternionic structure, then
\begin{equation*}
\begin{split}
IU(z)&=I(U'z\bar a)=U'z\bar a\mbox{\boldmath$i$}=
U'z(\bar a\mbox{\boldmath$i$}a) \bar a\\
&=U'z(a_{11}\mbox{\boldmath$i$}+a_{21}\mbox{\boldmath$j$}
+a_{31}\mbox{\boldmath$k$})\bar a=
a_{11}U'z\mbox{\boldmath$i$}\bar a+a_{21}U'z\mbox{\boldmath$j$}\bar a
+a_{31}U'z\mbox{\boldmath$k$}\bar a\\
&=a_{11}U(z\mbox{\boldmath$i$})+a_{21}U(z\mbox{\boldmath$j$})
+a_{31}U(z\mbox{\boldmath$k$})=a_{11}UI(z)+a_{21}UJ(z)+a_{31}UK(z).
\end{split}
\end{equation*}This follows that
$\displaystyle
IU=a_{11}UI+a_{21}UJ+a_{31}UK$.
Since $IU(e_i)=U_{j}^jI_j^\ell e_\ell$,
a calculation shows that
$U_{i}^{j}I_{j}^{\ell}=a_{11}I_{i}^{j}U_j^\ell+a_{21}J_{i}^jU_j^{\ell}+
a_{31}K_{i}^jU_j^{\ell}$, similarly for $J,K$.
As 
\begin{equation}\label{newqu}
\left(\begin{array}{c}
I'\\
J'\\
K'\end{array}\right)=
{}^tA\left(\begin{array}{c}
I\\
J\\
K\end{array}\right)
\end{equation} is a new quaternionic structure (\cf \eqref{two-quaternionic}),
it follows that
\begin{equation}\label{1nqu}
\begin{split}
I_{ij}U^i_kU^j_{\ell}=
a_{11}I_{k\ell}+a_{21}J_{k\ell}+a_{31}K_{k\ell}&={I'}_{k\ell}.\\
J_{ij}U^i_kU^j_{\ell}=
a_{12}I_{k\ell}+a_{22}J_{k\ell}+a_{32}K_{k\ell}&={J'}_{k\ell}.\\
K_{ij}U^i_kU^j_{\ell}=
a_{13}I_{k\ell}+a_{23}J_{k\ell}+a_{33}K_{k\ell}&={K'}_{k\ell}.
\end{split}
\end{equation}
Then we obtain that

\begin{equation}\label{1integene}
\begin{split}
d{\om'}+{\om'}\we{\om'}&=
-({I'}_{ij}\mbox{\boldmath$i$}+{J'}_{ij}\mbox{\boldmath$j$}+
{K'}_{ij}\mbox{\boldmath$k$})\Bigl(u^2\se^i\we\se^j \Bigr.\\
&\Bigl. \  +\mathop{\sum}_a\om_a\we 2uw_a^i\se^j+\mathop{\sum}_{a<b}{\om_a\we\om_b}(2w_a^iw_b^j)\Bigr).
\end{split}
\end{equation}
We shall derive an invariant under the change
$\om'=\lam\cdot \om\cdot\bar \lam$.
Recall from \eqref{1se} that
\begin{equation}\label{qcc-equiva}
(\om'_1,\om'_2,\om'_3)=
(\om_1,\om_2,\om_3)u^2\cdot A.
\end{equation}
Let $\displaystyle
d\se^i=\se^j\wedge\f^i_j\ +\mathop{\sum}_a \om_a\we\tau^i_a$
be the structure equation \eqref{sasaki structure eq}.
We define $1$-forms ${\nu'}^i_a$ by
setting
\begin{equation}\label{tau1}
\left(\begin{array}{c}
{\nu'}^i_1\\
{\nu'}^i_2\\
{\nu'}^i_3
\end{array}\right)
=u^{-2}\cdot {}^tA\left(\begin{array}{c}
{\tau}^i_1\\
{\tau}^i_2\\
{\tau}^i_3
\end{array}\right).
\end{equation}Since
${\tau}^i_a\equiv 0\ \ \mbox{mod}\ \se^k\ (k=1,\cdots 4n)$
by \eqref{upper-equations},
note that
\begin{equation}\label{cond1}
{\nu'}^i_a\equiv 0\ \mbox{mod}\ \se^k.
\end{equation}
Using \eqref{qcc-equiva} and \eqref{tau1},
\begin{equation*}
\begin{split}
\mathop{\sum}_a
\om_a\we{\tau}^i_a
=({\om'}_1,{\om'}_2,{\om'}_3)\we
\left(\begin{array}{c}
{\nu'}^i_1\\
{\nu'}^i_2\\
{\nu'}^i_3
\end{array}\right)
=\mathop{\sum}_a
{\om'}_a\we{\nu'}^i_a,
\end{split}
\end{equation*}
the equation \eqref{sasaki structure eq} becomes
\begin{equation}\label{gene3}
d\se^i=\se^j\we{\f}^i_j+\mathop{\sum}_a{\om'}_a\we{\nu'}^i_a.
\end{equation}
Differentiate \eqref{gene3}, and then 
substitute \eqref{1integene}, we obtain
that
\begin{equation}\label{1eq'}
\begin{split}
\se^j\we(d{\f}^i_j-{\f}^{\sigma}_j\we{\f}^i_{\sigma}+
u^2{\bf J'}^{1}_{jk}\se^k\we{\nu'}^i_1+u^2{\bf J'}^{2}_{jk}\se^k\we{\nu'}^i_2
+u^2{\bf J'}^{3}_{jk}\se^k\we{\nu'}^i_3) \equiv 0 \ \mbox{mod} \ \om_\al.
\end{split}
\end{equation}
Taking into account \eqref{1eq'} (which corresponds to \eqref{1eq}),
we have the fourth-order tensor up to the terms $\om_1,\om_2,\om_3$:
\begin{equation}\label{vani}
\frac 12{T'}^i_{jk\ell}\se^k\we\se^\ell\equiv d{\f}_j^i-
{\f}_j^\sigma\we{\f}_\sigma^i
+\mathop{\sum}_au^2\cdot {\bf J'}^{a}_{jk}\se^k\we{\nu'}^i_a
-{\se}^i\we\se_j.
\end{equation}
Since $({I'}_{ij},{J'}_{ij},{K'}_{ij})=(I_{ij},J_{ij},K_{ij})A$
from \eqref{newqu} and \eqref{tau1},
\begin{equation*}
\begin{split}
&\mathop{\sum}_au^2\cdot {\bf J'}^{a}_{jk}\se^k\we{\nu'}^i_a
=\se^k\we u^2({I'}_{jk},{J'}_{jk},{K'}_{jk})\left(
\begin{array}{c}
{\nu'}^i_1\\
{\nu'}^i_2\\
{\nu'}^i_3
\end{array}\right)\\
&=\se^k\we ({I}_{jk},{J}_{jk},{K}_{jk})\left(
\begin{array}{c}
{\tau}^i_1\\
{\tau}^i_2\\
{\tau}^i_3
\end{array}\right)
=\se^k\we\mathop{\sum}_a{\bf J}^a_{jk}{\tau}^i_a.
\end{split}
\end{equation*}
The equation \eqref{vani} can be reduced to the following:
\begin{equation}\label{Nvani}
\begin{split}
{T'}^i_{jk\ell}\se^k\we\se^\ell&\equiv 
d{\f}_j^i-{\f}_j^\sigma\we{\f}_\sigma^i
+\se^k\we\mathop{\sum}_a{\bf J}^a_{jk}{\tau}^i_a
-\se^i\we\se_j.
\end{split}
\end{equation}From \eqref{ps-cur} and \eqref{cur}, we have shown

\begin{pro}\label{invarinat tensor}
If $\om'=\lam\cdot\om\cdot \bar\lam$ for which
$\om$ is a quaternionic $CR$ structure, then
the curvature tensor $T'$ satisfies that
${T'}^i_{jk\ell}={T}^i_{jk\ell}$.
In particular, $T=({T}^i_{jk\ell})$ 
is an invariant tensor under the pseudo-conformal quaternionic
$CR$ structure.
\end{pro}

\begin{remark}\label{extend of I,J,K}
{\em
\ {\bf 1.}\ \ Similarly, the quaternionic structures $\{I',J',K'\}$
extends to almost complex structures $\{\bar I',\bar J',\bar K'\}$
respectively. 

\par\noindent\quad {\bf 2.}\ \ Let $f\in \mathop{\Aut}_{QCR}(M)$ be an element 
satisfying \eqref{1se}. Then,\\
$f_*e_i=uU_i^ke_k$. Using \eqref{1nqu},
\begin{equation*}
\begin{split}
&If_*e_i=uU_i^kI_{k}^{j}e_j=
u(a_{11}I_{i}^{m}+a_{21}J_{i}^{m}+a_{31}K_{i}^{m}) U_m^j e_j\\
&=f_*((a_{11}I_{i}^{m}+a_{21}J_{i}^{m}+a_{31}K_{i}^{m})e_m)\\
&=f_*((a_{11}I+a_{21}J+a_{31}K)e_i).
\end{split}
\end{equation*}The similar argument to $J,K$ yields that
\begin{equation}\label{f-com}
\left(\begin{array}{c}
f^{-1}_*If_*\\
f^{-1}_*Jf_*\\
f^{-1}_*Kf_*
\end{array}\right)=
{}^tA
\left(\begin{array}{c}
I\\
J\\
K
\end{array}\right)\ \ \mbox{on}\ \mathcal D.
\end{equation}
}
\end{remark}
\smallskip

\vskip0.2cm

\subsection{Formula of Curvature tensor}\label{VC}
We shall find a formula of the tensor $T$.
Substitute \eqref{new-con},
\eqref{upper-equations} into \eqref{Nvani}:
\begin{equation*}
\begin{split}
{T}^i_{jk\ell}\se^k\we\se^{\ell}&=
d(\om^i_j+\mathop{\sum}_a{({\bf J}^a)^{i}_{j}}\om_a)
-(\om^{\sigma}_j+\mathop{\sum}_a{({\bf J}^a)^{\sigma}_{j}}\om_a)\we
(\om^i_{\sigma}+\mathop{\sum}_a{({\bf J}^a)^{i}_{\sigma}}\om_a)\\
&\ \ \ \ \ +\se^k\we({I}_{jk}\cdot{{I}}^{i}_{\ell}\se^\ell+
{J}_{jk}\cdot{J}^{i}_{\ell}\se^{\ell}+{K}_{jk}\cdot
{K}^{i}_{\ell}\se^{\ell}) -\se^i\we\se_j \ \ \ \mbox{mod}\ \om_a\\
&=d\om^i_j+\mathop{\sum}_a{({\bf J}^a)^{i}_{j}}d\om_a
-\om_j^{\sigma}\we\om^i_{\sigma}+\mathop{\sum}_a({{\bf J}^a_{jk}}{({\bf J}^a)^{i}_{\ell}})\se^k\we\se^{\ell}
-\se^i\we\se_j \ \ \mbox{mod}\ \om_a\\
&=(d\om^i_j-\om_j^{\sigma}\we\om^i_{\sigma})\\
&\ \ \ \ \ +\mathop{\sum}_a{({\bf J}^a)^{i}_{j}}
({-{\bf J}^a_{k\ell}})\se^k\we\se^{\ell}+\mathop{\sum}_a({{\bf J}^a_{jk}}
{({\bf J}^a)^{i}_{\ell}})\se^k\we\se^{\ell}-\se^i\we\se_j \ \ 
\mbox{mod}\ \om_a\\
&=\left(\frac 12R^i_{jk\ell}-\mathop{\sum}_a{({\bf J}^a)^{i}_{j}}
{{\bf J}^a_{k\ell}}+
\mathop{\sum}_a {{\bf J}^a_{jk}}{({\bf J}^a)^{i}_{\ell}}-
g_{j\ell}\cdot\delta^{i}_{k}\right)\se^k\we\se^{\ell}\ \ \mbox{mod}\ \om_a.
\end{split}
\end{equation*}By alternation, we have

\begin{equation} \label{Sten}
\begin{split}
{T}^i_{jk\ell}=R^i_{jk\ell}- \left(2\mathop{\sum}_a{({\bf J}^a)^{i}_{j}}
{{\bf J}^a_{k\ell}}-\mathop{\sum}_a{{\bf J}^a_{jk}}{({\bf J}^a)^{i}_{\ell}}+
\mathop{\sum}_a{{\bf J}^a_{j\ell}}{({\bf J}^a)^{i}_{k}} 
+(g_{j\ell}\delta^{i}_{k}-g_{jk}\delta^{i}_{\ell})\right).
\end{split}
\end{equation}

Recall the space of all curvature tensors $\mathcal R({\rm Sp}(p,q)\cdot
{\rm Sp}(1))$. (See \cite{AL} for example.)
It decomposes into the direct sum
$\mathcal R_0({\rm Sp}(p,q)\cdot{\rm Sp}(1))\oplus
\mathcal R_{\HH\PP}({\rm Sp}(p,q)\cdot{\rm Sp}(1))$ $(n\geq 2)$.
Here $\mathcal R_0$ is the space of those curvatures with zero Ricci forms
and $\mathcal R_{\HH\PP}\approx \RR$ is the space of curvature tensors of
 the quaternionic pseudo-K\"ahler projective space ${\HH}{\PP}^{p,q}$
 (\cf Definition \ref{p-projective}).
\smallskip
\par\noindent {\bf Case ${\bf n\geq 2}$}.\ Since we know that
$R^i_{ji\ell}=R_{j\ell}=(4n+8)g_{j\ell}$
from \eqref{Ein2},
the curvature tensor $T=({{T}}^i_{jk\ell})$ satisfies
the {\it tracefree} condition:
\begin{equation*}
{{T}}_{j\ell}=({T}^i_{ji\ell})=
(4n+8)g_{j\ell}-\Bigl(3\cdot 3g_{j\ell}+
(4n-1)g_{j\ell}\Bigr)=0.
\end{equation*}
This implies that
our curvature tensor $T$ belongs to 
$\mathcal R_0({\rm Sp}(p,q)\cdot{\rm Sp}(1))$ when $n\geq 2$.
\smallskip
\par\noindent {\bf Case ${\bf n=1}$}.
When ${\rm dim}\ M=7$,
either $p=1,q=0$ or $p=0,q=1$. 
Choose the orthonormal basis $\{e_i\}_{i=1,2,3,4}$
with $e_1=e, e_2=Ie, e_3=Je, e_4=Ke$.
Form another curvature tensor:
\begin{equation}\label{4q-curvature}
\begin{split}
&{R'}^i_{jk\ell}=(g_{j\ell}\delta_k^i-g_{jk}\delta_{\ell}^i)+
\Bigl[I_{j\ell}I^{i}_{k}-I_{jk}I^{i}_{\ell}+2I^{i}_{j}I_{k\ell}\Bigr.\\
&\quad \Bigl.+J_{j\ell}J^{i}_{k}-J_{jk}J^{i}_{\ell}+2J^{i}_{j}J_{k\ell}
+K_{j\ell}K^{i}_{k}-K_{jk}K^{i}_{\ell}+2K^{i}_{j}K_{k\ell}\Bigr].
\end{split}
\end{equation}
For any two distinct $e_i,e_j$,
\begin{equation*}
\begin{split}
{R'}^i_{jij}&=
(g_{jj}\delta_i^i-g_{ji}\delta_{j}^i)+
\Bigl[I_{ji}I^{i}_{i}-I_{ji}I^{i}_{j}+2I_{ij}I^{i}_{j}+
J_{ji}J^{i}_{i}-J_{ji}J^{i}_{j}+
2J^{i}_{j}J_{ij}\Bigr.\\
&\Bigl.+K_{ji}K^{i}_{i}-K_{ji}K^{i}_{j}+2K^{i}_{j}K_{ij}\Bigr]
=g_{jj}+3\Bigl[I_{ij}{I^{i}_{j}}+J_{ij}{J^{i}_{j}}+K_{ij}{K^{i}_{j}}\Bigr].
\end{split}
\end{equation*}
Since $i\neq j$ and 
$e_j$ is either one of $\pm Ie_i,\pm Je_i,\pm Ke_i$,
$I_{ij}^2+J_{ij}^2+K_{ij}^2=1$ (for example, 
if $e_j=Ie_i$, then ${I^{i}_{j}}^2=1$, $J^{i}_{j}=0, K^{i}_{j}=0$ so that
$I_{ij}I_{j}^i=g_{jj}$.)
Thus, $\displaystyle {R'}^i_{jij}=4g_{jj}$. It follows from the Schur's
theorem (cf. \cite{KONO} for example) that
\begin{equation}\label{Schur}
{R'}^i_{jk\ell}=4(g_{j\ell}\delta_k^i-g_{jk}\delta_{\ell}^i).
\end{equation}
When $n=1$,
we conclude that
\begin{equation}\label{ano-cur}
{T}^i_{jk\ell}=R^i_{jk\ell}-{R'}^i_{jk\ell}
=R^i_{jk\ell}-4(g_{j\ell}\delta_k^i-g_{jk}\delta_{\ell}^i).
\end{equation}
As the curvature
$R^i_{jk\ell}$ satisfies the Einstein property
from \eqref{Ein2};
$R_{j\ell}=4\cdot 3g_{j\ell}$, the scalar curvature
$\sigma=4\cdot 12$. On the other hand, the curvature
tensor $R^i_{jk\ell}$ has the decomposition:
\[
R^i_{jk\ell}=W^i_{jk\ell}+\frac{4\cdot 12}{4\cdot 3}
(g_{j\ell}\delta_k^i-g_{jk}\delta_{\ell}^i)
\]in the space $\mathcal R({\rm SO}(4))$ where
 ${\rm SO}(4)={\rm Sp}(1)\cdot{\rm Sp}(1)$.
Hence, 
\begin{equation}\label{4cur}
{T}^i_{jk\ell}=W^i_{jk\ell}\in\mathcal R_0({\rm SO}(4))
\end{equation}for which $W^i_{jk\ell}$ is the Weyl curvature tensor 
(of $(U/\mathcal E,\hat g)$). 

\smallskip
\par\noindent {\bf Case ${\bf n=0}$}.
If ${\rm dim}\ M=3$, then
the above tensor is empty, so we simply set $T=0$.
Define the Riemannian metric on a neighborhood $U$
of a $3$-dimensional
pseudo-conformal quaternionic $CR$ manifold $M$:
\begin{equation}\label{riemann}
g_x(X,Y)=\om_1(X)\cdot \om_1(Y)+\om_2(X)\cdot \om_2(Y)+
\om_3(X)\cdot \om_3(Y)
\end{equation}$(\forall\ X,Y\in T_xU)$.
Suppose that $\om'=\lam\cdot\om\cdot \bar\lam$.
Since $({\om'}_1,{\om'}_2,{\om'}_3)=u^2\cdot({\om}_1,{\om}_2,{\om}_3)A$
for $A\in {\rm SO}(3)$, the metric $g$ changes into
${g'}={\om'}_1\cdot {\om'}_1+{\om'}_2\cdot {\om'}_2+
{\om'}_3\cdot {\om'}_3$ satisfying that
\begin{equation}\label{conf}
{g'}_x(X,Y)=u^4\cdot g_x(X,Y)\ \   (\forall\ X,Y\in T_xU).
\end{equation}Then $g'$ is conformal to $g$ on $U$.
Define $TW(\om)$ to be the Weyl-Schouten tensor
$TW(g)$ of the Riemannian metric $g$ on $U$.
Then, it turns out that
\begin{equation}\label{schou}
TW(\om')=TW(\om).
\end{equation}As a consequence,
$TW(\om)$ is an invariant tensor of $U$ under the change
$\om'=\lam\cdot\om\cdot \bar\lam$.

\section{Uniformization of pseudo-conformal $QCR$ structure} 
\label{CI}
If $\{\omega^{(\alpha)},(I^{(\al)}, J^{(\al)}, 
K^{(\al)}), g_{(\al)}, U_{\al}\}_{\alpha\in\Lambda}$
is a pseudo-conformal quaternionic $CR$ structure
on $M$ where $\displaystyle \mathop{\cup}_{\al\in \Lambda}U_{\al}=M$,
then we have the curvature tensor
$T^{(\al)}=({{}^{(\al)}T}^i_{jk\ell})$ on each $(U_\al,\om^{(\al)})$
$(n\geq 1)$.
Similarly, $TW^{(\al)}=TW(\om^{(\al)})$ on $(U_\al,\om^{(\al)})$
for $3$-dimensional case\ $(n=0)$.
Then it follows from Proposition \ref{invarinat tensor} and 
\eqref{schou} that
if $\om^{(\be)}=\lam_{\al\be}\cdot\om^{(\al)}\cdot\bar\lam_{\al\be}$
on $U_{\alpha}\cap U_{\beta}$,
then
\begin{equation}\label{curvature-int}
\begin{split}
&T^{(\al)}=T^{(\be)},\\
&TW^{(\al)}=TW^{(\be)}.
\end{split}
\end{equation}
By setting $T|U_\al=T^{(\al)}$ (respectively
$TW|U_{\al}=TW^{(\al)}$), the curvature $T$ 
(respectively $TW$) is globally defined on a $(4n+3)$-dimensional
pseudo-conformal quaternionic $CR$ manifold $M$\ $(n\geq 0)$.
This concludes that

\begin{theorem}\label{formula}
Let $M$ be a pseudo-conformal quaternionic $CR$
manifold of dimension $4n+3$ $(n\geq 0)$.
If $n\geq 1$, 
there exists  the fourth-order curvature tensor
$T=({T}^i_{jk\ell})$ on $M$ satisfying that:
\begin{itemize}
\item[(i)] When $n\geq 2$, 
$T=({T}^i_{jk\ell})\in \mathcal R_0({\rm Sp}(p,q)\cdot{\rm Sp}(1))$
which has the formula:
\begin{equation*}\label{Tensor}
\begin{split}
{T}^i_{jk\ell}&=R^i_{jk\ell}-\Bigl\{(g_{j\ell}\delta^{i}_{k}
-g_{jk}\delta^{i}_{\ell})+
\Bigl[I_{j\ell}I^{i}_{k}-I_{jk}I^{i}_{\ell}+
2I^{i}_{j}I_{k\ell}\Bigr.\Bigr.\\
&\ \ \ \ \ \ \Bigl.\Bigl.+J_{j\ell}J^{i}_{k}-J_{jk}J^{i}_{\ell}+
2J^{i}_{j}J_{k\ell}
+K_{j\ell}K^{i}_{k}-K_{jk}K^{i}_{\ell}+
2K^{i}_{j}K_{k\ell} \Bigr]\Bigr\}.
\end{split}
\end{equation*}
\item[(ii)] When $n=1$,
$T=({W}^i_{jk\ell})\in \mathcal R_0({\rm SO}(4))$ 
which has the same formula as the Weyl conformal curvature tensor.
\end{itemize}
If $n=0$, there exists the fourth-order curvature tensor
$TW$ on $M$ which has the same formula as the
Weyl-Schouten curvature tensor.
\end{theorem}

We associated to a pseudo-conformal quaternionic $CR$ structure
$(\{\om_a\},\{J_a\},\{\xi_a\})_{a =1,2,3}$ 
the pseudo-Sasakian metric 
$\displaystyle g=\mathop{\sum}_{a=1}^3\om_a\cdot \om_a+\pi^*\hat g$
on $U$ for which
$\mathcal E\ra (U,g)\stackrel{\pi}\lra (U/\mathcal E,\hat g)$
is a pseudo-Riemannian submersion and
the quotient $(U/\mathcal E,\hat g,
\{{\hat {I}}_i,{\hat {J}}_i,{\hat {K}}_i\}_{i\in\Lambda})$
is a quaternionic pseudo-K\"ahler manifold
by Theorem \ref{q-K-s}.
Let ${}^{(g)}{R}^i_{jk\ell}$ (respectively  $\hat R^i_{jk\ell}$)
 denote the
curvature tensor of $g$ (respectively $\hat g$).
If $R_{\HH\PP}$ is the generator of
$\mathcal R_{\HH\PP}({\rm Sp}(p,q)\cdot{\rm Sp}(1))
\approx \RR$ $(n\geq 2)$, then it can be described as (cf. \cite{AL}):
\begin{equation}\label{q-curvature}
R_{\HH\PP}=(g_{j\ell}g_{ik}-g_{jk}g_{i\ell})+
\mathop{\sum}_{a=1}^3{{\bf J}^a_{j\ell}}{{\bf J}^a_{ik}}
-\mathop{\sum}_{a=1}^3{{\bf J}^a_{jk}}{{\bf J}^a_{i\ell}}+
2\mathop{\sum}_{a=1}^3{{\bf J}^a_{ij}}{{\bf J}^a_{k\ell}}
\end{equation}where
$i,j,k,\ell$ run over $\{1,\cdots,4n\}$.
Then the formula $(12.8)$ of curvature
tensor of $g$ \cite{TA} $(n\geq 1)$ shows the following.

\begin{lemma}\label{tano-result}
\begin{equation}\label{cal-tano}
\begin{split}
\pi^*\hat R_{ijk\ell}&=
{}^{(g)}{R}_{ijk\ell}+
\left(\mathop{\sum}_{a=1}^{3}{{\bf J}^a_{j\ell}}{{\bf J}^a_{ik}}
-\mathop{\sum}_{a=1}^{3}{{\bf J}^a_{jk}}{{\bf J}^a_{i\ell}}+
2\mathop{\sum}_{a=1}^{3}{{\bf J}^a_{ij}}{{\bf J}^a_{k\ell}}\right)\\
&={}^{(g)}{R}_{ijk\ell}-(g_{j\ell}\delta_{ik}-
g_{jk}\delta_{i\ell})+R_{\HH\PP}.
\end{split}
\end{equation}

\end{lemma}

We now state the uniformization theorem.  
\begin{theorem}\label{uni1}
{\bf (1)}\ 
Let $M$ be a $(4n+3)$-dimensional pseudo-conformal quaternionic $CR$ 
manifold $(n\geq 1)$.
If the curvature tensor $T$ vanishes, then
$M$ is locally modelled on $S^{4p+3,4q}$
with respect to the group ${\rm PSp}(p+1,q+1)$.
\par\ \par
{\bf (2)}\ If $M$ is
a $3$-dimensional pseudo-conformal quaternionic $CR$
manifold whose curvature tensor $TW$ vanishes, then
$M$ is conformally flat \ie $($locally modelled on $S^3$
with respect to the group ${\rm PSp}(1,1)$ $)$.
\end{theorem}

\begin{proof}
Using \eqref{updown-curvature} and \eqref{q-curvature},
the formula of Theorem \ref{formula}
becomes 
\begin{equation}\label{HP}
{T}^i_{jk\ell}=\pi^*\hat R^i_{jk\ell}-R_{\HH\PP}.
\end{equation}
Compared this with \eqref{cal-tano},
we obtain that
\begin{equation}\label{T-formula}
{T}^i_{jk\ell}={}^{(g)}{R}^i_{jk\ell}-
(g_{j\ell}\delta^{i}_{k}-g_{jk}\delta^{i}_{\ell}).
\end{equation} The equality
\eqref{T-formula} is also true for $n=1$. In fact,
when $n=1$,
$R_{\HH\PP}=4(g_{j\ell}\delta_k^i-g_{jk}\delta_{\ell}^i)$
(cf. \eqref{4q-curvature}, \eqref{Schur}) and from \eqref{cal-tano},
${}^{(g)}{R}^i_{jk\ell}-(g_{j\ell}\delta^{i}_{k}-g_{jk}\delta^{i}_{\ell})
=\pi^*\hat R^i_{jk\ell}-R_{\HH\PP}={T}^i_{jk\ell}$
by \eqref{ano-cur}.

Suppose that $T$ (respectively $TW$)
vanishes identically on $M$.
First we show that
$M$ is locally isomorphic to $S^{4p+3,4q}$
(respectively $M$ is locally isomorphic to $S^3$.)
As $T|U_\al=({}^{(\al)}{T}^i_{jk\ell})=0$ on $U_\al$,
for brevity, we omit $\al$ so that
$T=({T}^i_{jk\ell})$ vanishes identically on $U$ for $n\geq 2$.
As a consequence, 
\begin{equation}\label{const1}
{}^{(g)}{R}^i_{jk\ell}=g_{j\ell}\delta^{i}_{k}-g_{jk}\delta^{i}_{\ell}
\ \ \mbox{on} \ \mathcal D|U.
\end{equation}
Since  $(U,g)$ is a pseudo-Sasakian $3$-structure
with Killing fields $\{\xi_1,\xi_2,\xi_3\}$,
the normality of \eqref{sasaki-condition}
can be stated as $\displaystyle {}^{(g)}R(X,\xi_a)Y=g(X,Y)\xi_a-g(\xi_a,Y)X$
(cf. \cite{TA}). It turns out that
\begin{equation}\label{const2}
{}^{(g)}R(\xi_a,X,Y,Z)=g(X,Z)g(\xi_a,Y)-g(X,Y)g(\xi_a,Z)
\end{equation} $(\forall X,Y,Z\in TU)$.
Then \eqref{const1} and \eqref{const2} imply
that $(U,g)$ is the space  of positive constant curvature.
As $\hat R^i_{jk\ell}= R_{\HH\PP}$
by \eqref{HP},
the quotient space $(U/\mathcal E,\hat g)$
is locally isometric to the quaternionic pseudo-K\"ahler projective space 
$({\HH}{\PP}^{p,q},\hat g_0)$.
(Note that if ${T}^i_{jk\ell}=0$ for $n=1$, then
$\pi^*\hat R^i_{jk\ell}=R^i_{jk\ell}
=4(\delta_{j\ell}\delta_k^i-\delta_{jk}\delta_{\ell}^i)$
from \eqref{ano-cur}. When $p=1,q=0$, 
the base space 
$(U/\mathcal E,\hat g)$ is locally isometric to the standard sphere
$S^4$ which is identified
with the $1$-dimensional quaternionic projective space 
$\HH\PP^1$.
If $p=0,q=1$, then $(U/\mathcal E,\hat g)$ is locally isometric 
to the quaternionic hyperbolic space $\HH^1_{\HH}=\HH\PP^{0,1}$
in which we remark that the metric $\hat g$ is negative definite.)
Hence, the bundle: 
$\mathcal E\ra (U,g)\stackrel{\pi}\lra (U/\mathcal E,\hat g)$
is locally isometric to
the Hopf bundle as the Riemannian submersion $(n\geq 1)$
(\cf Theorem \ref{spaceform}): 
\[
{\rm Sp}(1)\ra (\Sigma^{4p+3,4q}_{\HH},g_0)\lra (\HH\PP^{p,q},\hat g_0).
\]
This is obviously true for $n=0$.

Let $\varphi: (U,g)\ra (\Sigma^{4p+3,4q}_{\HH},g_0)$
be an isometric immersion preserving the above principal bundle.
If $V_0=\{\xi^0_1,\xi^0_2,\xi^0_3\}$ is the distribution of
Killing vector fields which generates ${\rm Sp}(1)$ of the above 
Hopf bundle, then we can assume that $\f_*\xi_a=\xi^0_a$
$(a=1,2,3)$ (by a composite of some element of ${\rm Sp}(1)$
if necessary). As $\om_a(X)=g(\xi_a,X)$ $(X\in TU)$
and $\om^0_a(X)=g_0(\xi_a,X)$ $(X\in T\Sigma^{4p+3,4q}_{\HH})$
respectively, the equality $g=\f^*g_0$ implies that
\begin{equation}\label{const3}
\om_a=\f^*\om^0_a\  \ (a=1,2,3),\  \ \om=\f^*\om_0.
\end{equation}
If we represent $\displaystyle
\f^*\se^i=\se^kT^i_k+\mathop{\sum}_a\om_a v^i_a$
for some matrix $T^i_j$ and $v^i_a\in \RR$, then
the equality $\f_*\xi_a=\xi^0_a$ shows that
$v^i_a=0$ for $i=1,\cdots,4n$. Thus,
\begin{equation}\label{repre-mat}
\f^*\se^i=\se^kT^i_k.
\end{equation}

For each $\al\in\Lambda$,
we have an immersion $\f_\al:U_\al\ra \Sigma^{4p+3,4q}_{\HH}$ as above
so that there is a collection of charts
$\{U_\al,\f_\al\}_{\al\in\Lambda}$ on $M$.
Put $g_{\al\be}=\f_\be\circ {\f_\al}^{-1} :
\f_\al(U_{\alpha}\cap U_{\beta})\ra \f_\be(U_{\alpha}\cap U_{\beta})$
when $U_{\alpha}\cap U_{\beta}\neq\emptyset$.
It suffices to prove that 
$g_{\al\be}$ extends uniquely to an element of 
${\PSp}(p+1,q+1)=\mathop{\Aut}_{QCR}(S^{4p+3,4q})$.
Suppose that
\begin{equation}\label{inter-al-be}
\omega^{(\beta)}=\lambda\cdot\omega^{(\alpha)}\cdot
\overline{\lambda}=u^2\cdot a\cdot\omega^{(\alpha)}\cdot 
\bar a \ \ \ \mbox{on}\ U_{\alpha}\cap U_{\beta}\neq\emptyset
\end{equation}where $\lambda=u\cdot a$.
The immersions $\f_\al:U_\al\ra \Sigma^{4p+3,4q}_{\HH}$,
$\f_\be:U_\be\ra \Sigma^{4p+3,4q}_{\HH}$ satisfy $\om^{(\al)}=\f_\al^*\om_0$,
$\om^{(\be)}=\f_\be^*\om_0$ as in \eqref{const3}.
If we put $\mu=\lambda\circ {\f_\al}^{-1}$
on $\f_\al(U_{\alpha}\cap U_{\beta})$, then the above relation shows that
\begin{equation}\label{conuj3}
g_{\al\be}^*\om_0=\mu\cdot\om_0\cdot\bar\mu.
\end{equation}
Using
the equation 
that $\displaystyle d\om^{(\al)}_a(X,Y)=
g^{(\al)}(X,J^{(\al)}_aY)$ $(\forall\ X,Y \in \mathcal D, a=1,2,3)$
from \eqref{1trans}, calculate that
\begin{equation*}
\begin{split}
d\om^0_a(\f_{\al*}J_a^{(\al)}X,\f_{\al*}Y)&=
d\om_a(J_a^{(\al)}X,Y)=g^{(\al)}(X,Y)\\
&=g_0(\f_{\al*}X,\f_{\al*}Y)=
d\om^0_a(J_a^{0}\f_{\al*}X,\f_{\al*}Y).
\end{split}
\end{equation*}As $d\om^0_a$ is nondegenerate on $\mathcal D$,
for each $\al\in\Lambda$ we have
\begin{equation}\label{holo}
\f_{\al*}\circ J_a^{(\al)}=J_a^{0}\circ\f_{\al*}\ \ \mbox{on}\ 
\mathcal D\ \ (a=1,2,3).
\end{equation}
Let $\displaystyle {\f_{\al}^*}\se^i={\se_{(\al)}^k}\cdot{}^{(\al)}T^i_k$
for some matrix ${}^{(\al)}T^i_k$ as in \eqref{repre-mat}.
Then \eqref{holo} means that
${}^{(\al)}T_i^k\cdot (J^a)_{k}^j=(J^a)_i^k\cdot {}^{(\al)}T_k^j$,
which implies that $\displaystyle {}^{(\al)}T^i_k\in {\GL}(n,\HH)$.
As $g^{(\al)}(X,Y)=g_{0}(\f_{\al*}X,\f_{\al*}Y)$
from \eqref{const3}, this reduces to
\begin{equation}\label{Spnn}
{}^{(\al)}T^i_k\in {\rm Sp}(p,q).
\end{equation}
Let $\{\om_{(\al)},
\om^i_{(\al)}\}_{i=1,\cdots,n}$, $\{\om_{(\be)},
\om^i_{(\be)}\}_{i=1,\cdots,n}$ be
two coframes on the
intersection $U_\al\cap U_\be$
where $\om_{(\al)}$ is a ${\rm Im}\HH$-valued $1$-form 
and each $\om^i_{(\al)}$ is a $\HH$-valued $1$-form,
simlarly for $\be$.   
Noting \eqref{inter-al-be}, the coordinate change of the fiber $\HH^n$
satisfies that
 \begin{equation}\label{fiberchange}
\left(\begin{array}{c}
\om_{(\be)}\\
\om^1_{(\be)}\\
\vdots  \\
\om^n_{(\be)}
\end{array}\right)
=\left(\begin{array}{c|ccc}
\lam     &   &\mbox{\Large $0$}& \\
              \hline              
         &   &                  & \\
\tilde v^i     &    & U' &        \\ 
               &    &    &   
\end{array}\right)
\left(\begin{array}{c}
\om_{(\al)}\\
\om^1_{(\al)}\\
\vdots  \\
\om^n_{(\al)}
\end{array}\right)\cdot \bar\lam.
\end{equation}
In order to transform them into the real forms,
recall that ${\rm GL}(n,\HH)\cdot{\rm GL}(1,\HH)$ 
is the maximal closed 
subgroup of ${\rm GL}(4n,\RR)$
acting on $\RR^{4n}$ preserving the standard quaternionic structure
$\{I,J,K\}$.
For each fiber of $\mathcal D_\al\ (=\mathcal D_\be)$,
there exists a matrix 
$\tilde U=(\tilde U^i_j)=U'\cdot \lam\in {\rm GL}(n,\HH)\cdot {\rm GL}(1,\HH)$
such that:
\begin{equation}\label{cap}
e^{(\al)}_j=\tilde U^i_je^{(\be)}_i.
\end{equation}with respect to the basis 
$\{e^{(\al)}_i\}_x\in (\mathcal D_\al)_x$, 
$\{e^{(\be)}_i\}_x\in (\mathcal D_\be)_x$.
From Corollary \ref{conf-metric},
\begin{equation*}
\pm u^2\delta_{k\ell}=u^2g_{(\al)}(e^{(\al)}_k,e^{(\al)}_\ell)=
g_{(\be)}(\tilde U^i_ke^{(\be)}_i,\tilde U^j_\ell e^{(\be)}_j)
=\pm\delta_{ij}\tilde U^i_k \tilde U^j_\ell,
\end{equation*}so $(u^{-1}\tilde U^i_k)\in {\rm Sp}(p,q)\cdot {\rm Sp}(1)
={\rm GL}(n,\HH)\cdot {\rm GL}(1,\HH)\cap {\rm SO}(4p,4q)$ up to conjugacy
$(n\geq 1)$. Put $\displaystyle
U=(U^i_k)=(u^{-1}\tilde U^i_k)\in {\rm Sp}(p,q)\cdot {\rm Sp}(1)$,
then
\begin{equation}\label{tildeU}
\tilde U=uU=(uU^i_k)\in {\rm Sp}(p,q)\cdot {\rm Sp}(1)\times \RR^+.
\end{equation}
Using coframes 
$\{\se^i_{(\al)}\}, \{\se^i_{(\be)}\}$
(induced from $\{\om^i_{(\al)},\om^i_{(\be)}\}_{i=1,\cdots,n}$),
the equation  \eqref{cap} translates into
$\se_{(\be)}^i=\se^k_{(\al)}\tilde U^i_k$ on $\mathcal D$.
Using \eqref{fiberchange}, it follows that
\[
\se_{(\be)}^i=\se^k_{(\al)}\tilde U^i_k+\sum_{a=1}^3\om^{(\al)}_a\cdot v^i_a
\ \ \mbox{on}\ U_\al\cap U_\be.
\]Here $v^i_a$ are determined by $\tilde v^i$, see \eqref{1se}.
Then,
\begin{equation}\label{conju4}
\begin{split}
{g_{\al\be}}^*(\se^i)&=
({\f_\al}^{-1})^*{\f_\be}^*(\se^i)=
({\f_\al}^{-1})^*({\se_{(\be)}^j}\cdot{}^{(\be)}T^i_j)\\
&=({\f_\al}^{-1})^*\left((\se^k_{(\al)}\tilde U^j_k+\sum_{a=1}^3
\om^{(\al)}_a\cdot v^j_a)\cdot{}^{(\be)}T^i_j\right)\\
&=\se^\ell\cdot{({}^{(\al)}T^{-1})}_\ell^k\tilde U^j_k\cdot{}^{(\be)}T^i_j
+\sum_{a=1}^3\om^0_a\cdot (v^j_a\cdot{}^{(\be)}T^i_j).
\end{split}\end{equation}
If we put $S=(S_\ell^i)=(({{}^{(\al)}T}^{-1})_\ell^k\cdot
\tilde U^j_k\cdot{}^{(\be)}T^i_j)$, then
\eqref{tildeU} and \eqref{Spnn} imply
$S\in {\rm Sp}(p,q)\cdot {\rm Sp}(1)\times \RR^+$.
Therefore, $g_{\al\be}$ satisfies the conditions 
of \eqref{1se} from \eqref{conuj3}, \eqref{conju4}. 
Then the diffeomorphism
$g_{\al\be}:\f_\al(U_{\alpha}\cap U_{\beta})\ra 
\f_\be(U_{\alpha}\cap U_{\beta})$ is viewed locally 
as an element of  
$\mathop{\Aut}_{QCR}(S^{4p+3,4q})={\rm PSp}(p+1,q+1)$
because 
$\Sigma^{4p+3,4q}_{\HH}\subset S^{4p+3,4q}$.
As ${\rm PSp}(p+1,q+1)$ acts real analytically on $S^{4p+3,4q}$,
$g_{\al\be}$ extends uniquely to an element
of ${\rm PSp}(p+1,q+1)$.
Therefore, the collection of charts
$\{U_\al,\f_\al\}_{\al\in\Lambda}$
gives rise to a uniformization of 
a pseudo-conformal quaternionic $CR$ manifold $M$ with respect to
$({\rm PSp}(p+1,q+1),S^{4p+3,4q})$.

Recall the $3$-dimensional
conformal geometry $({\rm PO}(4,1),S^3)$ for which
the orthogonal Lorentz group ${\rm PO}(4,1)$
is isomorphic as a Lie group to ${\rm PSp}(1,1)$.  
Then the same is true for $({\rm PSp}(1,1),S^3)=
({\rm PO}(4,1),S^3)$ $(n=0)$.

\end{proof}

\section{Quaternionic bundle}\label{q-bundle}
It is known that the first Stiefel-Whitney class  is 
the obstruction to the existence of
a global $1$-form of the contact structure and
the first Chern class is the obstruction to the existence
of a global $1$-form of the complex contact
structure respectively.
It is natural to ask whether the first Pontrjagin class $p_1(M)$
is the obstruction to the existence
of global $1$-form of the pseudo-conformal quaternionic structure
(respectively pseudo-conformal quaternionic $CR$ structure)
on a $(4n+3)$- manifold $M$ $(n\geq 1)$.
In order to see that, we need the elementary properties of the quaternionic
bundle theory. However our structure group is 
${\rm GL}(n,\HH)\cdot{\rm GL}(1,\HH)$ but not ${\rm GL}(n,\HH)$.
To our knowledge the fundamental properties of the bundle
theory in this case are not proved explicitly. So we prepare
the necessary facts.
Let $\mathcal D$ be the $4n$-dimensional bundle defined by 
$\displaystyle \mathcal D=\mathop{\cup}_{\al}\mathcal D_\al$ where
$\mathcal D_\al=\mathcal D|U_{\alpha}=\mathop{\Null\ }\omega^{(\alpha)}$
in which there is the  relation on the intersection $U_\al\cap U_\be$:
\begin{equation}\label{interlast}
\omega^{(\beta)}=\bar \lambda\cdot\omega^{(\alpha)}\cdot
\lambda=u^2\cdot \bar a\om^{(\be)}\cdot a \ \ \mbox{where}\
\lam=u\cdot a\in \HH^*.
\end{equation}
We have already discussed the transition functions on $\mathcal D$
in $\S \ref{fiberchange}$ (\cf \eqref{CI}).
In fact,

The gluing condition of the quaternionic bundle
$\mathcal D$ in $U_\al\cap U_\be$ is given by
\[
\left(\begin{array}{c}
v_1^{(\al)}\\
\vdots \\
v_n^{(\al)}
\end{array}\right)
=uT\left(\begin{array}{c}
v_1^{(\be)}\\
\vdots \\
v_n^{(\be)}
\end{array}\right)\cdot a,
\]in which $\displaystyle u(T\cdot \bar a)\in {\rm Sp}(n)\cdot {\rm Sp}(1)\times \RR^+$.
\begin{definition}\label{qbundle}
A quaternionic $n$-dimensional bundle
is a vector bundle over a paracompact Hausdorff space $M$
with fiber isomorphic to the $n$-dimensional
quaternionic vector space $\HH^n$. For an open cover
 $\{U_\al\}_{\al\in\Lambda}$ of $M$, if
$U_\al\cap U_\be\neq\emptyset$, then there exists a
transition function $g_{\al\be}:U_\al\cap U_\be\ra {\rm GL}(n,\HH)\cdot
{\rm GL}(1,\HH)$.
\end{definition}
As a consequence, $\mathcal D$ is 
a quaternionic $n$-dimensional bundle on $M$.

As ${\rm GL}(1,\HH)\cdot{\rm GL}(1,\HH)\approx {\rm SO}(4)\times \RR^+$,
the quaternionic line bundle is isomorphic to an oriented real 
$4$-dimensional bundle. Define the inner product $\langle ,\rangle$ on $\HH^n$
by
\[
\langle z,w\rangle=\bar z_1w_1+\cdots+\bar z_nw_n.
\]Then, $\langle ,\rangle$ satisfies that
$\langle z,w\cdot\lambda\rangle=\langle z,w\rangle\cdot\lambda,\
\langle z\cdot\lambda,w\rangle =\bar\lam\langle z,w\rangle,\ \langle z,w\rangle=\overline{\langle w,z\rangle}$ for $\lambda\in\HH$, and so on. By a subspace $W$ in $\HH^n$ we
mean a right $\HH$-module.
Choose $v_0\in \HH^n$.
Let $V=\{v_0\cdot \lam\ |\ \lam\in\HH \}$ be a $1$-dimensional subspace
of $\HH^n$. Denote $V^{\perp}=\{v\in\HH^n |\ \langle v,x\rangle =0, \forall\ x\in V\}$.
Then it is easy to check that $V^{\perp}$ is a right $\HH$-module
for which there is a decomposition:
$\HH^n=V\oplus V^{\perp}$ as a right $\HH$-module.
The following is a quaternionic analogue of the splitting theorem.

\begin{pro}\label{splitting}
Given a quaternionic $n$-dimensional bundle $\xi$, there
exists a quaternionic line bundle $\xi_i$ $(i=1,\cdots, n)$
over a paracompact Hausdorff space $N$ and a $($splitting$)$
map $f:N\ra M$ for which:

\begin{equation*}
\begin{split}
&(1)\quad f^*\xi=\xi_1\oplus\cdots \oplus \xi_n.\\
&(2)\quad f^*: H^*(M)\ra H^*(N) \ \mbox{is injective. Moreover,}\\
&(3)\quad \mbox{The bundle isomorphism}\ 
b:\xi_1\oplus\cdots \oplus \xi_n \ra \xi\
\ \mbox{compatible with}\ f\ \mbox{can be} \\
&\ \ \ \ \quad \mbox{chosen to preserve  the inner product}.
\end{split}
\end{equation*}

\end{pro}

\begin{proof}
Let $\HH^n-\{0\}\ra \xi_0\stackrel{\pi}\lra M$ be the subbundle of
$\xi$ consisting of nonzero sections. Noting that
$\HH^n$ is a right $\HH$-module, it induces a fiber bundle
with fiber the quaternionic $n$-dimensional projective space 
$\HH\PP^{n-1}$:
\[
\HH\PP^{n-1}\ra Q\stackrel{q}\lra M.
\]Since 
the cohomology group $H^*(\HH\PP^{n-1};\ZZ)$ is a free abelian group,
$q^*: H^*(M)\ra H^*(Q)$ is injective by the Leray-Hirsch's theorem
(cf. \cite{MI}.)
Put
\[
q^*\xi=\{(\ell,v)\in Q\times \xi\ |\ q(\ell)=\pi(v)\}.
\]Then, $(q^*\xi,{\rm pr},Q)$ is a quaternionic bundle.
Let $\xi_1=\{(\ell,v)\in q^*\xi\ |\ v\in\ell\}$ which is the 
$1$-dimensional quaternionic subbundle of $q^*\xi$. 
Choose a (right) $\HH$-inner product
on $\xi$. Then it induces a (right) $\HH$-inner product
on $q^*\xi$ such that the bundle projection 
${\rm Pr}:q^*\xi\ra \xi$ preserves the inner product obviously.
Moreover, we obtain that
\[
q^*\xi=\xi_1\oplus {\xi_1}^{\perp}.
\]
Since ${\xi_1}^{\perp}$ is an $(n-1)$-dimensional
quaternionic bundle over $Q$, an induction for $n-1$
implies that there exist a paracompact Hausdorff space $N$
and a splitting map
$f_1:N\ra Q$ such that
$\displaystyle f_1^*\xi_1^{\perp}=\xi_2\oplus\cdots\oplus \xi_n$ and
$f_1^* :H^*(Q)\ra H^*(N)$ is injective.
If $\displaystyle b_1:\xi_2\oplus\cdots\oplus \xi_n\ra {\xi_1}^{\perp}$ is the bundle map
compatible with $f_1$, then by induction $b_1$ preserves 
the inner product on the fiber
between $\xi_2\oplus\cdots\oplus \xi_n$ and ${\xi_1}^{\perp}$.
Putting $f=q\circ f_1:N\ra M$, we see that
$f^*:H^*(M)\ra H^*(N)$ is injective and $f^*\xi=
f^*_1\xi_1\oplus\xi_2\oplus\cdots\oplus \xi_n$.
Let $\displaystyle {\rm Pr}_1\times b_1: f^*_1\xi_1\oplus(\xi_2\oplus\cdots\oplus \xi_n)
\ra \xi_1\oplus\xi_1^{\perp}$ be the bundle map.
Then the map $\displaystyle {\rm Pr}\circ ({\rm Pr}_1\times b_1):
f^*_1\xi_1\oplus\xi_2\oplus\cdots\oplus \xi_n\lra \xi$ is compatible
with $f$ and preserves the inner product.
This proves the induction step for $n$.

\end{proof}

Let $\xi$ be a quaternionic line bundle over $M$ with gluing condition
on  $U_\al\cap U_\be$:
\begin{equation}\label{line-glu1} 
z_\al=\bar \lam(x) z_\be\mu(x)=
u(x)\cdot\bar b(x)z_\be a(x)\  \ (u>0,a,b\in {\rm Sp(1)}).
\end{equation}
Consider the tensor
$\displaystyle\bar\xi\mathop{\otimes}_\HH\xi$ so that
the gluing condition on $U_\al\cap U_\be$
is given by          
\begin{equation*}
\begin{split}
&(\bar z_\al\mathop{\otimes}_\HH z_\al)
={u}^2(x)\bar a(x)(\bar z_\be b(x)
\mathop{\otimes}_\HH\bar b(x)
z_\be)a(x)\\
&={u}^2(x)\bar a(x)
(\bar z_\be\mathop{\otimes}_{\HH}z_\be) a(x).
\end{split}
\end{equation*}Then $\displaystyle\bar\xi\mathop{\otimes}_\HH\xi$
is a quaternionic line bundle over $M$
whenever $\xi$ is a quaternionic line bundle.
 
\begin{lemma}\label{pon1}
If $\displaystyle \bar\xi\mathop{\otimes}_{\HH} \xi$ is viewed as a 
$4$-dimensional real vector bundle, then
$\displaystyle p_1(\bar\xi\mathop{\otimes}_{\HH} \xi)=p_1(\bar\xi)+
p_1(\xi)$. Moreover, $p_1(\bar\xi)=p_1(\xi)$ so that
$\displaystyle p_1(\bar\xi\mathop{\otimes}_{\HH} \xi)=2p_1(\xi)$.
\end{lemma}

\begin{proof}
Let $\gamma$  be the canonical $4$-dimensional real 
vector bundle over $B{\rm SO}(4)$ (cf. \cite{MI}).
Then, $\xi$ is determined by a
classifying map $f:M\ra B{\rm SO}(4)$ such that
$f^*\gamma=\xi$.
Let ${\rm pr}_i:B{\rm SO}(4)\times B{\rm SO}(4)\ra B{\rm SO}(4)$
be the projection $(i=1,2)$. As $\gamma$ inherits a quaternionic structure
from $\xi$ through the bundle map, there is a quaternionic line bundle
$\displaystyle
{\rm pr}_1^*\bar \gamma\mathop{\otimes}_{\HH}{\rm pr}_2^*\gamma$
over $B{\rm SO}(4)\times B{\rm SO}(4)$. Now,
let $h:B{\rm SO}(4)\times B{\rm SO}(4)\ra B{\rm SO}(4)$ be
a classifying map of this bundle
so that $\displaystyle h^*\gamma={\rm pr}_1^*\bar \gamma
\mathop{\otimes}_{\HH}{\rm pr}_2^*\gamma$.
When $\iota_i
:B{\rm SO}(4)\ra B{\rm SO}(4)\times B{\rm SO}(4)$ is the inclusion map
on each factor, $\iota_1^*{\rm pr}_2^*\gamma$ is the trivial 
quaternionic line bundle (we simply put $\se_{\bf h}^1$) and so
$\displaystyle
\iota_1^*h^*p_1(\gamma)=
\iota_1^*p_1({\rm pr}_1^*\bar \gamma\mathop{\otimes}_{\HH}
{\rm pr}_2^*\gamma)=p_1(\bar \gamma\mathop{\otimes}_{\HH}
\se_{\bf h}^1)=p_1(\bar\gamma)$.
Similarly, $\iota_2^*h^*p_1(\gamma)=p_1(\gamma)$. Hence we obtain that
\[
h^*p_1(\gamma)=p_1(\bar \gamma)\times 1+1\times p_1(\gamma).
\]
Let $f':M\ra B{\rm SO}(4)$ be a classifying map for $\bar \xi$
such that ${f'}^*\gamma=\bar \xi$.
Then the map $h(f'\times f)d$ composed of the diagonal map $d:M\ra M\times M$
satisfies that
\[
(h(f'\times f)d)^*\gamma={f'}^*\bar\gamma\mathop{\otimes}_{\HH} 
f^*\gamma=\bar\xi\mathop{\otimes}_{\HH}\xi.
\]Therefore, $\displaystyle p_1(\bar\xi\mathop{\otimes}_{\HH}\xi)=
d^*({f'}\times f)^*(p_1(\bar \gamma)\times 1+1\times p_1(\gamma))
=p_1({f'}^*\bar\gamma)+p_1(f^*\gamma)=
p_1(\bar\xi)+p_1(\xi)$.

Next, the conjuagte $\bar \xi$ is isomorphic to $\xi$
as real $4$-dimensional vector bundle without orientation.
But the correspondence
$(1,\mbox{\boldmath$i$},\mbox{\boldmath$j$},
\mbox{\boldmath$k$})\mapsto (1,-\mbox{\boldmath$i$},-\mbox{\boldmath$j$},
-\mbox{\boldmath$k$})$ gives an isomorphism of $\bar \xi$
onto $(-1)^3\xi$. And so, the complexification $\bar\xi_{\CC}$ of $\bar\xi$
(viewed as a real vector bundle)
 is isomorphic to $(-1)^6\xi_{\CC}=\xi_{\CC}$.
By definition, $p_1(\bar\xi)=p_1(\xi)$.

\end{proof}

\subsection{Relation between Pontrjagin classes}\label{q-Pont}
Suppose that $\{\omega^{(\alpha)},(I^{(\al)}, J^{(\al)}, 
K^{(\al)}), g_{(\al)}, U_{\al}\}_{\alpha\in\Lambda}$
represents a pseudo-conformal quaternionic structure
$\mathcal D$ on a $(4n+3)$-dimensional manifold 
$\displaystyle M=\mathop{\cup}_{\al\in \Lambda}U_{\al}$.
Let $L$ be the quotient bundle
$TM/\mathcal D$. Choose the local vector fields
 $\{\xi_1^{(\al)},\xi_2^{(\al)},\xi_3^{(\al)}\}$ on each neighborhood $U_\al$ such that $\om_a^{(\al)}(\xi^{(\al)}_b)=\delta_{ab}$. Then, $L|U_\al$ is spanned
by $\{\xi_1^{(\al)}\}_{i=1,2,3}$ for each $\al\in \Lambda$.
Moreover, the gluing condition between $L|U_\al$ and
$L|U_\be$ is exactly given by
\begin{equation}\label{vectgl}
\left(\begin{array}{c}
\xi_1^{(\al)}\\
\xi_2^{(\al)}\\
\xi_3^{(\al)}\end{array}\right)=
u^2A
\left(\begin{array}{c}
\xi_1^{(\be)}\\
\xi_2^{(\be)}\\
\xi_3^{(\be)}\end{array}\right).
\end{equation} (Compare Definition \ref{pc-qcr}.)
It is easy to see that
$\displaystyle \mathop{\sum}_{a=1}^3
\om_a^{(\al)}\cdot \xi_a^{(\al)}
=\mathop{\sum}_{a=1}^3\om_a^{(\be)}\cdot \xi_a^{(\be)}$.
We can define a section $\se:TM \ra L$ which
is an $L$-valued $1$-form by setting 
\begin{equation}\label{se-om}
\se|U_\al=\om_1^{(\al)}\cdot \xi_1^{(\al)}+\om_2^{(\al)}\cdot \xi_2^{(\al)}+
\om_3^{(\al)}\cdot \xi_3^{(\al)}.
\end{equation}Then there is an exact sequence of bundles:
$\displaystyle 1\ra \mathcal D\ra TM\stackrel{\se}\lra L\ra 1$.

Let $E$ be the $1$-dimensional quaternionic line bundle
obtained from the union $\displaystyle \mathop{\bigcup}_{\al\in\Lambda}^{}
U_\al\times {\HH}$ by identifying 
\begin{equation}\label{L+}
(p,z_\al)\sim (q,z_\be)\ \mbox{if and only if}\
\left\{\begin{array}{l}
p=q\in U_\al\cap U_\be,\\
z_\al=\lam\cdot z_\be\cdot\bar\lam=u^2 a\cdot z_\be\cdot\bar a.
\end{array}\right.
\end{equation}If $L\oplus \se$ is
the Whitney sum composed of the trivial (real) line bundle $\se$
 on $M$, then it is easy to see that $L\oplus \se$
is isomorphic to the $1$-dimensional quaternionic line bundle $E$.
Then the Pontrjagin class $p_1(E)=p_1(L\oplus \se)$.
We prove that
\begin{theorem}\label{pont}
The first Pontrjagin classes of $M$ and the bundle $L$ has the relation:
\[
2p_1(M)=(n+2)p_1(L\oplus \se).
\]
\end{theorem}

\begin{proof}
As $\mathcal D$ is a quaternionic bundle in our sense, there is a splitting map
$f:N\ra M$  such that
$f^*\mathcal D=\xi_1\oplus\cdots\oplus\xi_n$ from Proposition \ref{splitting}.
Let $\Psi:\xi_1\oplus\cdots\oplus\xi_n\ra \mathcal D$ be a bundle map which
is compatible with $f$. Since $\Psi$ is a right $\HH$-linear
map on the fiber at each point $x\in N$,
we can describe
\[
\Psi\left(\begin{array}{c}
v_1\\
\vdots\\
v_n\end{array}\right)_x=
P(x)\left(\begin{array}{c}
v_1\\
\vdots\\
v_n\end{array}\right)_{f(x)}
\]for some function $P:N\ra {\rm GL}(n,\HH)$.
With respect to an appropriate inner product on $\mathcal D$ and
the direct inner product on $\xi_1\oplus\cdots\oplus\xi_n$,
 $\Psi$ preserves the inner product between
them by (3) of Theorem \ref{splitting}, which implies
\begin{equation}\label{unitary}
P(x)\in {\rm  Sp}(n).
\end{equation}
We examine the gluing condition of each $\xi_i$ on $f^{-1}(U_\al)\cap
f^{-1}(U_\be)\neq \emptyset$.
For $x\in f^{-1}(U_\al)\cap
f^{-1}(U_\be)$, let $v_i^{(\al)}\in \xi_i|f^{-1}(U_\al)$.
Suppose that there is an elemt $v_i^{(\be)}\in \xi_i|f^{-1}(U_\be)$
such that $v_i^{(\al)}\sim v_i^{(\be)}$, \ie 
\[
v_i^{(\al)}=\bar\lam_i v_i^{(\be)}\mu_i\ \  (\lam_i, \mu_i\in \HH^*;
 i=1,\cdots, n).
\] Since $\Psi(v_i^{(\al)})\sim \Psi(v_i^{(\be)})$ at $f(x)$,
it follows from \eqref{interlast} that
$\displaystyle \Psi\left(\begin{array}{c}
v_1^{(\al)}\\
\vdots \\
v_n^{(\al)}
\end{array}\right)
=uT\cdot \Psi\left(\begin{array}{c}
v_1^{(\be)}\\
\vdots \\
v_n^{(\be)}
\end{array}\right)\cdot a\ \  \mbox{at}\ f(x)\in U_a\cap U_\be$
where $T\in{\rm Sp}(n)$.
As
\begin{equation*}
\begin{split}
P\left(\begin{array}{c}
v_1^{(\al)}\\
\vdots \\
v_n^{(\al)}
\end{array}\right)
=
\Psi\left(\begin{array}{c}
v_1^{(\al)}\\
\vdots \\
v_n^{(\al)}
\end{array}\right)
=uT\cdot P\left(\begin{array}{c}
v_1^{(\be)}\\
\vdots \\
v_n^{(\be)}
\end{array}\right)\cdot a
=
P(uP^{-1}TP)\left(\begin{array}{c}
v_1^{(\be)}\\
\vdots \\
v_n^{(\be)}
\end{array}\right)\cdot a,
\end{split}
\end{equation*}it follows that
\begin{equation*}
\left(\begin{array}{c}
v_1^{(\al)}\\
\vdots \\
v_n^{(\al)}
\end{array}\right)
=u(P^{-1}TP \left(\begin{array}{c}
v_1^{(\be)}\\
\vdots \\
v_n^{(\be)}
\end{array}\right)\cdot a.
\end{equation*}
Since  $v_i^{(\al)}=\bar\lam_i v_i^{(\be)}\mu_i$ as above, we have
\begin{equation*}
\begin{split}
&{\bf (1)}\quad u(x)P(x)^{-1}T(f(x))P(x)=
\left(\begin{array}{ccc}
\bar\lam_1(x)&\cdots& 0\\
0    &\ddots& 0\\
0    &0&\bar\lam_n(x)
\end{array}\right).\\
&{\bf (2)}\quad \left(\begin{array}{ccc}
\mu_1(x)&\cdots& 0\\
0    &\ddots& 0\\
0    &0&\mu_n(x)
\end{array}\right)
=\left(\begin{array}{ccc}
a(x)&\cdots& 0\\
0    &\ddots& 0\\
0    &0&a(x)
\end{array}\right).
\end{split}
\end{equation*}From \eqref{unitary} and ${\bf (1)}$,
calculate
\begin{equation*}
\begin{split}
&\left(\begin{array}{ccc}
|\lam_1|^2& \hdots &0\\
\vdots    &\ddots  &\vdots\\
0         &\hdots  &|\lam_n|^2
\end{array}\right)
={u}^2P^{-1}TP\circ P^*T^*P={u}^2{\rm I}_n.
\end{split}\end{equation*}Hence,
$\displaystyle 
\lam_i=u\cdot\nu_i\ \mbox{where}\ 
\nu_i={\lam_i}/|\lam_i|\in {\rm Sp}(1)$.
It follows from ${\bf (2)}$ that $\mu_i=a$ for each $i$.
We have that
\begin{equation}\label{line-glu}
v_i^{(\al)}=u(x)\cdot\bar\nu_i(x)v_i^{(\be)}a(x)\ \ \ (i=1,\cdots,n).
\end{equation}
 Therefore each $\xi_i$ is a quaternionic line bundle over $N$ with the above gluing condition \eqref{line-glu} on
 $\displaystyle f^{-1}(U_\al)\cap f^{-1}(U_\be)$.
If we consider the tensor
$\displaystyle\bar\xi_i\mathop{\otimes}_\HH\xi_i$, then
the gluing condition on $\displaystyle f^{-1}(U_\al)\cap f^{-1}(U_\be)$
is given by          
\begin{equation*}
\begin{split}
(\bar v_i^{(\al)}\mathop{\otimes}_\HH v_i^{(\al)})
={u}^2(x)\bar a(x)(\bar v_i^{(\be)}\mathop{\otimes}_{\HH}v_i^{(\be)})a(x).
\end{split}
\end{equation*}
In view of \eqref{L+},
we see that $\displaystyle
f^*(L\oplus\se^1)=\bar\xi_i\mathop{\otimes}_{\HH}\xi_i\ \ (i=1,\cdots,4n)$.
By Lemma \ref{pon1},
$f^*p_1(L\oplus\se^1)=2p_1(\xi_i)$ for each $i$.
Since $f^*p_1(\mathcal D)\equiv p_1(\xi_1)+\cdots+p_1(\xi_n)$
{\rm mod} $2$-torsion in $H^4(N;\ZZ)$,
$\displaystyle f^*(2p_1(\mathcal D))=
2p_1(\xi_1)+\cdots+2p_1(\xi_n)=nf^*p_1(L\oplus\se^1)=nf^*p_1(L)$.
As $f^*$ is injective,
$\displaystyle 2p_1(\mathcal D)=np_1(L)$ in $H^4(M;\ZZ)$.
Noting $TM=\mathcal D\oplus L$, we have
$2p_1(M)=(n+2)p_1(L)$.
\end{proof} 

\begin{corollary}\label{global-obstruction}
Let $(M,\mathcal D)$ be a $(4n+3)$-dimensional simply connected
pseudo-conformal quaternionic manifold associated with the
local forms $\{\omega^{(\alpha)},(I^{(\al)}, J^{(\al)}, 
K^{(\al)}), g_{(\al)}, U_{\al}\}_{\alpha\in\Lambda}$.
Then the following are equivalent.
\begin{enumerate}
\item $2p_1(M)=0$. In particular, the rational Pontrjagin class vanishes. 
\item $L$ is the trivial bundle so that
$\{\xi_\al\}_{\al=1,2,3}$ exists  globally on $M$.                              \item There exists a global ${\rm Im}\HH$-valued
$1$-form $\omega$ on $M$ which represents a pseudo-conformal
quaternionic structure $\mathcal D$. In particular,
there exists a hypercomplex structure $\{I,J,K\}$ on $\mathcal D$.
\end{enumerate}
\end{corollary}

\begin{proof}
First note that the Whitney sum $L\oplus \se^1$ is the quaternionic 
line bundle $E$ with structure group lying in ${\rm SO}(3)\times \RR^+
\subset {\rm Sp}(1)\cdot {\rm Sp}(1)\times \RR^+$.
As above we have the quaternionic line bundle of $\ell$-times tensor
$\displaystyle\mathop{\otimes}_{\HH}^{\ell}E$ with 
structure group ${\rm SO}(3)\times \RR^+$.
Viewed as the $4$-dimensional real vector bundle, it determines 
the classifying map $g:M\ra B({\rm SO}(3)\times \RR^+)=B{\rm SO}(3)$.
Note that $p:B({\rm Sp}(1)\times \RR^+)\ra B({\rm SO}(3)\times \RR^+)$ is
the two-fold covering map. As $M$ is simply connected by the hypothesis, 
the map $g$
lifts to a classifying map
$\tilde g: M\ra B{\rm Sp}(1)$ such that $g=p\circ \tilde g$.
Let $\gamma$ be the $4$-dimensional
universal bundle over $B{\rm SO}(3)$. (Compare \cite{MI}.) Then the pull back
$p^*\gamma$ is the $4$-dimensional canonical 
bundle over $B{\rm Sp}(1)=\HH\PP^{\infty}$ whose 
first Pontrjagin class $p_1(p^*\gamma)$
generates the cohomology ring $H^*(\HH\PP^{\infty};\ZZ)$. 
So the bundle $\displaystyle\mathop{\otimes}_{\HH}^{\ell}E$
is classified by the map $\tilde g$ where $[\tilde g]=
\tilde g^*p_1(p^*\gamma)\in H^4(M;\ZZ)$, which coincides with
$\displaystyle p_1(\mathop{\otimes}_{\HH}^{\ell}E)$.\\

{\em $1\Rightarrow 2$.}\ If $2p_1(M)=0$, then Theorem \ref{pont} shows 
$(n+2)p_1(L)=0$, \ie
 $\displaystyle p_1((\mathop{\otimes}_{\HH}^{n+2}E))=0$.
(See Lemma \ref{pon1}.)
Hence, the classifying map $\tilde g: M\ra B{\rm Sp}(1)$ for
$\displaystyle\mathop{\otimes}_{\HH}^{n+2}E$ is null homotopic so that
$\displaystyle\tilde g^*p^*\gamma=\mathop{\otimes}_{\HH}^{n+2}E$ is trivial.
There exists a family of functions $\{h_\al\}\in{\rm Sp}(1)\times \RR^+$
such that the transition function $g_{\al\be}(x)=\delta^1 h(\al,\be)(x)$
\ ($x\in U_\al\cap U_\be$).
As the gluing relation for $\displaystyle\mathop{\otimes}_{\HH}^{n+2}E$
is given by $u_{\al\be}^{2(n+2)}\bar a_{\al\be}\cdot z\cdot a_{\al\be}$,
letting $h_\al=a_\al\cdot u_\al\in {\rm Sp}(1)\times \RR^+$,
it follows that
\[
u_{\al\be}^{2(n+2)}\cdot \bar a_{\al\be}\cdot z\cdot a_{\al\be}=
(h_\al^{-1}h_\be)z=u_\al^{-1} u_\be a_\al\bar a_\be\cdot z
\cdot a_\be\bar a_\al\ \ \ (z\in \HH).
\]Then, $u_{\al\be}^{2(n+2)}=u_\al^{-1} u_\be\in \RR^+$ and 
$a_{\al\be}=\pm a_\be\bar a_\al$.
As $u_{\al\be}>0$, $\displaystyle u_{\al\be}=
(u_\al^{-1})^{\frac 1{2(n+2)}}\cdot u_\be^{\frac 1{2(n+2)}}$.
Since the gluing relation of $E=L\oplus \se$ is given 
by $z_\al=u_{\al\be}^2\cdot\bar a_{\al\be}\cdot z_\be\cdot a_{\al\be}$,
putting $\displaystyle {u'}_\al=(u_\al)^{\frac 1{(n+2)}}$, 
a calculation shows
$\displaystyle z_\al=
{u_{\al\be}}^2\cdot{\bar a_{\al\be}}\cdot z_\be \cdot a_{\al\be}
={{u'}_\al}^{-1}{u'}_\be\cdot a_\al
{\bar a}_\be\cdot z_\be\cdot a_\be{\bar a}_\al$.
Moreover if $C(\al)\in {\rm SO}(3)$ is the matrix
defined by $\displaystyle \bar a_{\al}\cdot \left(\begin{array}{c}
\mbox{\boldmath$i$}\\
\mbox{\boldmath$j$}\\
\mbox{\boldmath$k$}
\end{array}\right)\cdot a_{\al}=C(\al)
\left(\begin{array}{c}
\mbox{\boldmath$i$}\\
\mbox{\boldmath$j$}\\
\mbox{\boldmath$k$}
\end{array}\right)$ (similarly for $C(\be)$), then

\begin{equation}\label{cobound}
u_{\al\be}^2\cdot A^{\al\be}={u'}_\al^{-1}{u'}_\be\cdot C(\al)^{-1}\circ
C(\be).
\end{equation}
Substitute this into \eqref{vectgl}, it follows that
\begin{equation*}
{u'}_\al\cdot C(\al)\left(\begin{array}{c}
\xi_1^{(\al)}\\
\xi_2^{(\al)}\\
\xi_3^{(\al)}\end{array}\right)
={u'}_\be \cdot C(\be)\left(\begin{array}{c}
\xi_1^{(\be)}\\
\xi_2^{(\be)}\\
\xi_3^{(\be)}\end{array}\right)\ \ \mbox{on}\ U_\al\cap U_\be.
\end{equation*} 
We can define the vector fields $\{\xi_1,\xi_2,\xi_3\}$ on $M$
to be
\begin{equation}
\left(\begin{array}{c}
\xi_1\\
\xi_2\\
\xi_3\end{array}\right)\mbox{\Large $\bigl |$} U_\al=
{u'}_\al\cdot C(\al)\left(\begin{array}{c}
\xi_1^{(\al)}\\
\xi_2^{(\al)}\\
\xi_3^{(\al)}\end{array}\right).
\end{equation}
Then $\{\xi_1,\xi_2,\xi_3\}$ spans $L$, therefore, $L$ is trivial.

{\em $2\Rightarrow 3.$}\ 
Since
$\displaystyle (\om^{(\be)}_1, \om^{(\be)}_2, \om^{(\be)}_3)
=(\om^{(\al)}_1, \om^{(\al)}_2, \om^{(\al)}_3)u_{\al\be}^2\cdot A^{\al\be}$,
\eqref{cobound} implies that
\[
(\om^{(\be)}_1, \om^{(\be)}_2, \om^{(\be)}_3)
{u'}_\be^{-1}\cdot C(\be)^{-1}=
(\om^{(\al)}_1, \om^{(\al)}_2, \om^{(\al)}_3)
{u'}_\al^{-1}\cdot C(\al)^{-1} \  \mbox{on }\ U_\al\cap U_\be.
\]Then, the ${\rm Im}\HH$-valued $1$-form
$\omega$ on $M$ can be defined by
\begin{equation}
\om|U_\al=(\om^{(\al)}_1, \om^{(\al)}_2, \om^{(\al)}_3)
{u'}_\al^{-1}\cdot C(\al)^{-1}
\left(\begin{array}{c}
\mbox{\boldmath$i$}\\
\mbox{\boldmath$j$}\\
\mbox{\boldmath$k$}
\end{array}\right). 
\end{equation}
Note that $\om$ satisfies that
 $\om|U_\al=\bar\lam_\al\cdot\om^{(\al)}\cdot\lam_\al$ for some function
$\lam_\al:U_\al\ra \HH^*$ $(\al\in\Lambda)$.
Recall that two quaternionic structures on $U_\al\cap U_\be$
are related:
\begin{equation*}\label{two-quaternionic2}
\left(\begin{array}{c}
I^{(\al)}\\
J^{(\al)}\\
K^{(\al)}
\end{array}\right)
=A^{\al\be}\left(\begin{array}{c}
I^{(\be)}\\
J^{(\be)}\\
K^{(\be)}
\end{array}\right).
\end{equation*}As $A^{\al\be}=C(\al)^{-1}\circ C(\be)$,
it follows that
\begin{equation}\label{two-quaternionic3}
C(\al)\cdot\left(\begin{array}{c}
I^{(\al)}\\
J^{(\al)}\\
K^{(\al)}
\end{array}\right)
=C(\be)\cdot \left(\begin{array}{c}
I^{(\be)}\\
J^{(\be)}\\
K^{(\be)}
\end{array}\right).
\end{equation}Letting $\displaystyle \left(\begin{array}{c}
I\\
J\\
K
\end{array}\right)|U_\al=
C(\al)\cdot\left(\begin{array}{c}
I^{(\al)}\\
J^{(\al)}\\
K^{(\al)}
\end{array}\right)$, 
there exists a hypercomplex structure $\{I,J,K\}$ on $\mathcal D$.

{\em $3\Rightarrow 1.$}\ If the global ${\rm Im}\HH$-valued 
$1$-form $\om$ exists, then $\om$ defines a three independent
vector fields isomorphic to $L$, \ie $p_1(L)=0$.

\end{proof}

\end{document}